\documentclass{amsart}
\usepackage{amssymb}
\usepackage{amsthm}
\usepackage{amscd}
\usepackage{graphics}
\newtheorem*{reftheorem}{Theorem}
\newtheorem*{reflemma}{Lemma}

\newtheorem*{Theorem 1}{Theorem 1}
\newtheorem*{Theorem 2}{Theorem 2}
\newtheorem*{Theorem 1A}{Theorem 1A}
\newtheorem*{Theorem 2A}{Theorem 2A}
\newtheorem*{theoremcp}{Cup Product Theorem}

\newtheorem*{theoremst}{Steenrod Tower Theorem}

\newtheorem{example}{Example}
\newtheorem{introexample}{Example}
\newtheorem{ML}{Main Lemma}
\newtheorem{lemma}{Lemma}[section]

\newtheorem{prop}{Proposition}[section]

\newtheorem*{conj}{Conjecture}

\newtheorem*{vanishing lemma}{Vanishing Lemma}
\numberwithin{equation}{section}
\title{Minimax problems related to cup powers and Steenrod squares}
\author{Larry Guth}
\address{Department of Mathematics, Stanford, Stanford CA, 94305 USA}
\email{lguth@math.stanford.edu}

\begin{document}
\begin{abstract} If $F$ is a family of mod 2 k-cycles in the unit n-ball,
we lower bound the maximal volume of any cycle in $F$ in terms of the homology
class of $F$ in the space of all cycles.  We give examples to show that these
lower bounds are fairly sharp.
\end{abstract}

\maketitle

This paper is about minimax estimates for the volumes of
cycles in complicated families.  The simplest example of a minimax problem is a
classical result about curves in the unit disk.  First consider the
family of vertical lines in the unit disk.  The longest line in
the family has length 2.  Then consider any other family of curves
that sweeps out the unit disk.  One of the curves in the other
family must still have length at least 2.  We illustrate the
situation in Figure 1.

\vskip10pt

\includegraphics{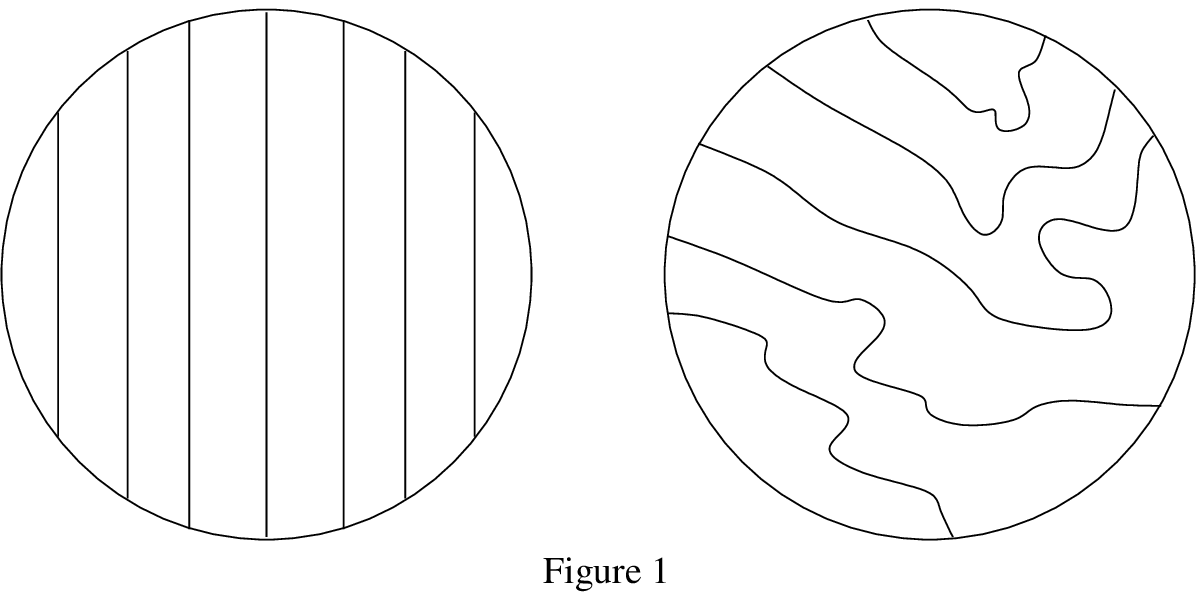}

We can write this result in the following form, which explains
why it is called a minimax estimate.

$$\inf_{F \in \mathbb{F}} \sup_{C \in F} \textrm{ length}(C) =
2.$$

\noindent In this equation, $\mathbb{F}$ denotes the set of all
1-parameter families of curves in the unit disk which sweep out
the disk.  The letter $F$ denotes a particular 1-parameter family
of curves in $\mathbb{F}$, and the letter $C$ denotes a curve in
the family $F$.

In this paper, we study more general minimax problems that may 
involve very high-parameter families.  For example, we will
study a minimax problem for p-parameter
families of planar curves for any integer $p$ and compute the 
asymptotic behavior as $p$ tends to infinity.
More generally, we will study analogous questions for k-dimensional
surfaces in the unit n-ball for any $k < n$.

We let $Z(k,n)$ denote the space of mod 2 relative k-cycles in the unit
n-ball.  (More precisely, we use the space of flat k-cycles.  In Section 1,
we give a self-contained definition of this space.)

By a family of k-cycles, we
mean a continuous map $F$ from a simplicial complex to $Z(k,n)$. 
If $\alpha$ is a cohomology class in
$H^*(Z(k,n), \mathbb{Z}_2)$, we say that $F$ detects
$\alpha$ if $F^*(\alpha) \not= 0$.  Then we define
$\mathbb{F}(\alpha)$ to be the set of all families of cycles that
detect the cohomology class $\alpha$.  We define a minimax volume
for the cohomology class $\alpha$ by the following formula.

$$\mathbb{V}(\alpha) := \inf_{F \in \mathbb{F}(\alpha)} \sup_{C
\in F} \textrm{Volume} (C).$$

This formula defines infinitely many minimax volumes $\mathbb{V}(\alpha)$,
and we will investigate how the minimax volume depends on the cohomology
class.

The first cohomology class that we will use measures whether
a family of cycles sweeps out the ball.  In Section 1, we will
define precisely what it means for a family of cycles to sweep out
the ball.  Suppose that $z$ is a
mod 2 (n-k)-cycle in the space $Z(k,n)$.  There is a cohomology
class $a(k,n)$ in $H^{n-k}(Z(k,n), \mathbb{Z}_2)$ whose pairing 
$<a(k,n), [z] >$ is equal to 1 if
the family $z$ sweeps out the unit ball (mod 2) and is equal to 0 if it
doesn't.  We call $a(k,n)$ the fundamental
cohomology class of $Z(k,n)$.  Determining $\mathbb{V}(a(k,n))$
is the classical minimax problem for families of k-cycles that
sweep out the unit n-ball.  For this problem, there are good
results due to Almgren (\cite{A2}) and Gromov
(\cite{G}, \cite{G3}), which we describe in detail in Section 2.

We get other cohomology classes by applying cohomology operations
to the fundamental cohomology class $a(k,n)$.  The next simplest
cohomology classes are cup powers of $a(k,n)$.  Recently, Gromov studied
the problem of estimating $\mathbb{V}( a(k,n)^p )$ in Section 8
of \cite{G3}.  He proved the following theorem.

\begin{Theorem 1} (Gromov) There are constants $0 < c(n) < C(n)$  so that
the following estimate holds.

$$c(n) p^{\frac{n-k}{n}} \le \mathbb{V}(a(k,n)^p) \le C(n)
p^{\frac{n-k}{n}}.$$

\end{Theorem 1}

In this paper, we will reprove Theorem 1 in detail.  We
will construct an explicit family of cycles that detects $a(k,n)^p$,
check that each cycle in the family has volume at most $C(n) p^{\frac{n-k}{n}}$,
and prove that this value is nearly optimal.

The main goal of the paper is to extend this analysis from cup powers
to towers of Steenrod squares.  Recall that $Sq^i$ 
denotes the Steenrod square cohomology operation that maps 
$H^N(X, \mathbb{Z}_2)$ to $H^{N+i}(X, \mathbb{Z}_2)$ 
for any space $X$.  For background on
Steenrod squares, see the chapter on them in \cite{H}.  If
$\alpha$ is a cohomology class in $H^N(X, \mathbb{Z}_2)$, we
write $Sq_i \alpha$ to denote $Sq^{N-i} \alpha \in H^{2N-i}(X,
\mathbb{Z}_2)$.  We write $Sq_i^2 \alpha$ to denote $Sq_i [Sq_i
\alpha]$, and in a similar way we define $Sq_i^p$ and $Sq_i Sq_j$.
Our second theorem estimates the minimax volume for any cohomology class
of the form $Sq_0^{Q_0} ... Sq_{n-k-1}^{Q_{n-k-1}} a(k,n)$, where
$Q_0, ..., Q_{n-k-1}$ are any non-negative integers.

\begin{Theorem 2} For each $\epsilon > 0$, there is a constant
$c(n, \epsilon) > 0$, and there is a constant $C(n)$ independent
of $\epsilon$, so that the following estimate holds.

$$c(n, \epsilon) \prod_{i=0}^{n-k-1} (2 -
\epsilon)^{\frac{n-k-i}{n-i} Q_i} \le \mathbb{V}(Sq_0^{Q_0} ...
Sq_{n-k-1}^{Q_{n-k-1}} a(k,n)) \le C(n) \prod_{i=0}^{n-k-1}
2^{\frac{n-k-i}{n-i} Q_i}.$$

\end{Theorem 2}

These formulas are pretty complicated, so we make a few comments
about them.  Let $\alpha$ be a cohomology class 
$Sq_0^{Q_0} ... Sq_{n-k-1}^{Q_{n-k-1}} a(k,n)$ lying in $H^d(Z(k,n), \mathbb{Z}_2)$.
At first, we might hope to estimate $V(\alpha)$ in terms of the dimension
$d$.  The formula in Theorem 2 implies that $V(\alpha)$ may be as large
as $d^{\frac{n-k}{n}}$ or as small as $d^{\frac{1}{k+1}}$ depending on
the values of $Q_0, ..., Q_{n-k-1}$.  Knowing $d$ gives us some idea
of $\mathbb{V}(\alpha)$, but our best guess would be subject to an error that
is polynomial in $d$.  Theorem 2 gives us an estimate for $\mathbb{V}(\alpha)$
which is accurate up to an error
of order $d^{\epsilon}$ for any $\epsilon > 0$.

The reader may want to know what fraction of the cohomology ring of $Z(k,n)$
is covered by Theorem 2.  The homotopy groups of spaces of cycles were studied
by Almgren in his thesis \cite{A}.  Almgren proved that the space
of relative integral k-cycles in the unit n-ball has homotopy
groups $\pi_{n-k} = \mathbb{Z}$ and all other homotopy groups
zero.  I believe that Almgren's argument should apply to the
space $Z(k,n)$ of mod 2 flat cycles.  The argument should prove
that $\pi_i(Z(k,n))$ is equal to $\mathbb{Z}_2$ for $i = n-k$ and
zero otherwise.  Unfortunately, this argument is not written down
anywhere as far as I know.  I hope to write an exposition of it in the
future.  If Almgren's argument applies, then $Z(k,n)$ is weak
homotopic to the Eilenberg-Maclane space $K(\mathbb{Z}_2, n-k)$.  
The mod 2 cohomology ring of the space $K(\mathbb{Z}_2, n-k)$ was
determined by Serre in \cite{S}.  The smallest non-zero cohomology
group is $H^{n-k}$ which is equal to $\mathbb{Z}_2$.  We call
the generator of this group $a$.  The entire cohomology ring
of $K(\mathbb{Z}_2, n-k)$ is a free $\mathbb{Z}_2$ algebra
with generators $Sq_1^{Q_1} ... Sq_{n-k-1}^{Q_{n-k-1}} a$, where
$Q_i \ge 0$ are any numbers.  (See \cite{H} for more information.)

To get a sense of what it means for a family of cycles to detect
a certain cohomology class, consider the following topological
properties.  First suppose that a family of cycles detects $a(k,n)$.  In other
words, the family of cycles sweeps out the unit n-ball.  This implies that
one of the cycles must go through the center of the ball.  More
generally, if we pick any point $x$ in the unit ball, one of the cycles
must go through $x$.  Next suppose that a family of cycles detects $a(k,n)^p$.  If
we pick any $p$ points $x_1, ..., x_p$ in the unit n-ball, then one
of the cycles in the family must go through all $p$ points.  The lower
bounds in Theorem 1 exploit this property.  The analogue of this property
for Steenrod squares is a little bit more complicated.  Suppose that a family of 
cycles detects $Sq_i a(k,n)$.  If we pick two particular points, $x_1$ and $x_2$, 
then our family may not contain a cycle that goes through
them both.  But if
we pick a continuous map $f$ from $S^i$ to the unit ball, then we
can find a point $\theta \in S^i$ so that a cycle from our family 
goes through both $f(\theta)$ and $f(- \theta)$.  There is
an analogous property for a family of cycles that detects $Sq_i^Q a(k,n)$,
but it becomes pretty complicated.  Roughly speaking, the lower bounds
in Theorem 2 exploit this property.

To give some context for our theorems, we consider some examples
of families of cycles.  The most interesting examples are families
of algebraic varieties.  These examples are discussed more in
Section 6 of the paper, where we prove the claims made in the
discussion below.  The one thing we don't prove is that the complex
algebraic hypersurfaces are an honest family of flat cycles - therefore,
Example 5 below is not completely rigorous.

\begin{introexample} Vertical lines in the unit disk.  
\end{introexample}

This is the example that began the paper.  The family $F$ is
parametrized by $[-1,1]$, with $F(t)$ being the line $x=t$
intersected with the unit disk.  Since $F(-1)$ and $F(1)$ are
both the empty cycle, we can think of the family as a map from
the circle to $Z(1,2)$.  Because the vertical lines sweep out the
disk, the class $F^*(a(1,2))$ is the non-zero class in $H^1(S^1,
\mathbb{Z}_2)$.  Since each vertical line has length at most 2,
it follows that $\mathbb{V}(a(1,2)) \le 2$, and this inequality
is sharp.

\begin{introexample} Sets of p vertical lines in the disk.
\end{introexample}

Start with the open p-simplex $-1 < t_1 < ... < t_p < 1$.  Define
$F$ on the open simplex by taking $F(t_1, ..., t_p)$ to be the
union of vertical lines $x= t_i$ for $1 \le i \le p$.  By
continuity, $F$ extends to map the closed simplex into $Z(1,2)$. 
The image of this map turns out to be a p-cycle in $Z(1,2)$
that detects $a(1,2)^p$.  Every p-tuple of vertical lines
has total length at most $2p$, and so this example shows that
$\mathbb{V}(a(1,2)^p) \le 2p$.

According to Theorem 1, $\mathbb{V}(a(1,2)^p) \sim p^{1/2}$.  For
large $p$, this family of cycles has much longer curves than
necessary.  Our next example shows how to improve it.

\begin{introexample} Planar real algebraic curves.
\end{introexample}

Let $P(x,y)$ be a real polynomial of degree at most $d$.  If $P$
is not uniformly zero, then we define $F(P)$ to be the
intersection of the real algebraic curve $\{ (x,y) | P(x,y)=0
\}$ with the unit disk.

We should mention that some of these polynomials have no
solutions.  For example, $P(x,y)$ could be $1$, or it could be $1
+ x^2 + y^2$.  In these cases, $F(P)$ is the empty cycle.  Some
of the curves in our family have one connected component, but
some have more than one connected component, while others are
empty.  In spite of the changing topology, the family of real
algebraic curves is continuous in the flat topology.  (We prove
this statement in Section 6.)

The space of all real polynomials of degree at most $d$ is a real
vector space of dimension ${d+2 \choose 2}$.  Two polynomials
define the same curve if one is a constant multiple of the other. 
Therefore, we can think of the space of degree d curves as a
family parameterized by $\mathbb{RP}^{D(d)}$, where $D(d) = {d+2
\choose 2} - 1 = (1/2) (d^2 + 3d)$.  This family detects the class
$a(1,2)^p$ for any $p \le D(d)$.  Also, the length of a degree
$d$ algebraic curve in the unit disk is less than $4 d$. 
Therefore, $\mathbb{V}(a(1,2)^{D(d)}) < 4 d$.  Since $D(d)$ is
roughly $d^2$, it follows that $\mathbb{V}(a(1,2)^p) < 10
p^{1/2}$.

According to Theorem 1, $\mathbb{V}(a(1,2)^p) \sim p^{1/2}$,
and so the family of degree d curves is roughly the optimal way
of detecting $a(1,2)^p$.

\begin{introexample} Real algebraic hypersurfaces.
\end{introexample}

The last construction generalizes to real algebraic hypersurfaces
in any dimension.  The space of degree d real algebraic
hypersurfaces in the unit n-ball can be parametrized by
$\mathbb{RP}^{D(d,n)}$ for a dimension $D(d,n)$ on the order of
$d^n$.  It detects $a(n-1,n)^p$ for any $p \le D(d,n)$.  Each
degree d hypersurface in the unit ball has volume at most $C(n)
d$.  These examples show that $\mathbb{V}(a(n-1,n)^p) \le C(n)
p^{1/n}$, which is roughly optimal according to Theorem 1.

\begin{introexample} Complex algebraic hypersurfaces.
\end{introexample}

The same construction applies to complex algebraic hypersurfaces. 
If $n$ is even, then we can think of the unit n-ball as the unit
ball in $\mathbb{C}^{n/2}$, and we can look at the degree $d$
complex hypersurfaces.  These are parametrized by
$\mathbb{CP}^{D(d,n/2)}$, where the dimension $D(d,n/2)$ is on the
order of $d^{n/2}$.  Each degree d complex hypersurface in the
unit n-ball has volume at most $C(n) d$.  Therefore, we get an
upper bound on minimax volumes $\mathbb{V}(a(n-2,n)^p) \le C(n)
p^{2/n}$ for even $n$.  This upper bound is roughly optimal
according to Theorem 1.

We can also think of complex hypersurfaces as integral cycles. 
We discuss the situation for integral cycles in Appendix 2.

\begin{introexample} Translates of real algebraic planar curves
in $\mathbb{R}^3$.
\end{introexample}

This example is a modification of the family of degree d real
algebraic planar curves.  Recall that the degree d curves formed
a family $F(d,2): \mathbb{RP}^{D(d,2)} \rightarrow Z(1,2)$,
detecting the cohomology class $a(1,2)^{D(d,2)}$.  Now we define
a new family of 1-cycles in the 3-ball by using translates of the
degree d curves.  The new family is parametrized by
$\mathbb{RP}^{D(d,2)} \times [-1, 1]$.  Our new family $F$ is
defined by taking $F(p, t)$ to be the restriction to the unit
3-ball of the product $F(d,2)(p) \times \{ t \} \subset B^2(1)
\times [-1,1]$.  If $t = \pm 1$, then
$F(p,t)$ is the empty cycle, so $F$ extends to a continuous
family parametrized by $\mathbb{RP}^{D(d,2)} \times S^1$.  The
family $F$ detects $Sq_1^Q a(1,3)$ for any $Q$ with $2^Q \le
D(d,2)$.  Each 1-cycle in $F$ has length at most $C d$.  Since
$D(d,2)$ is roughly $d^2$, this example shows that
$\mathbb{V}(Sq_1^Q a(1,3)) \le C 2^{Q/2}$.  According to Theorem
2, this example is approximately sharp in the sense that for any
$\epsilon$, we can choose a constant $c(n, \epsilon) > 0$ so that
$\mathbb{V}(Sq_1^Q a(1,3)) \ge c(n, \epsilon) (2-
\epsilon)^{Q/2}$.

\begin{introexample} Products of previous examples.
\end{introexample}

We can also take products of previous examples.  For example, we
can look at the product $F(d,2) \times F(d,2)$, which defines a
family of 2-cycles in the unit 4-ball parametrized by
$\mathbb{RP}^{D(d,2)} \times \mathbb{RP}^{D(d,2)}$.  This family
detects the class $a(2,4)^p$ for any $p \le D(d,2)$.  Each
surface in the family has area at most $C d^2$.  This
example shows that $\mathbb{V}(a(2,4)^p) \le C p$.  According to
Theorem 1, the actual value of $\mathbb{V}(a(2,4)^p)$ is much
smaller, on the order of $p^{1/2}$, which can be achieved by
looking at complex hypersurfaces.  Similarly, most other products
lead to families of cycles that are far from optimal with respect
to our minimax problem.

\vskip10pt

Now that we have seen some examples of families of cycles, we
discuss the proofs of the lower bounds in Theorems 1 and 2.

Lower bounds for $\mathbb{V}(a(k,n))$ for general $k$ were
first proven by Almgren using his version of Morse theory on the
space of cycles \cite{A2}. He proved that a family of k-cycles
sweeping out the unit n-sphere must contain a cycle with volume
at least equal to that of the unit k-sphere.  This statement
implies a lower bound for $\mathbb{V}(a(k,n))$. 
Almgren's argument involves a lot of geometric measure theory. 
In \cite{G}, Gromov gave a lower bound for $\mathbb{V}(a(k,n))$ 
by using the isoperimetric inequality repeatedly.  We
include Gromov's argument in Section 2.

The lower bounds in Theorem 1 are proven by combining a lower
bound for $\mathbb{V}(a(k,n))$ with Lusternik-Schnirelmann theory.  To give
the idea, we explain how to bound $\mathbb{V}(a(1,2)^p)$.  Let $B_1, ...,
B_p$ be disjoint disks inside the unit disk, each with radius $r
\sim p^{-1/2}$.

Let $S(i) \subset Z(1,2)$ be the subset of 1-cycles $z \in
Z(1,2)$ so that $z \cap B_i$ has length at most $r$.  (We use the
letters $S(i)$ as an abbreviation for ``cycles which are small in
$B_i$''.)  The dashed curve in Figure 2 illustrates a cycle in
$S(1)$.  This cycle does not belong to $S(2)$ or $S(3)$. 

\includegraphics{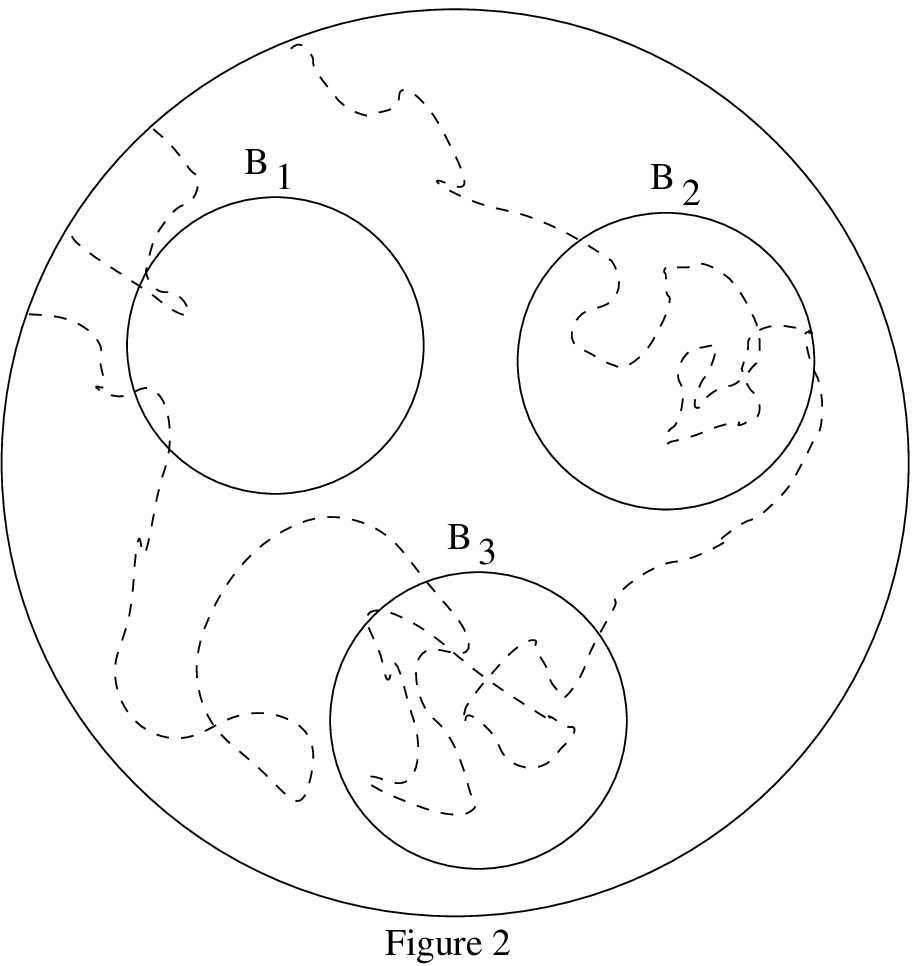}

Any family of curves that sweeps out the unit disk must also
sweep out each $B_i$.  By scaling the classical result at the
beginning of this introduction, we know that any family of curves that
sweeps out $B_i$ must contain a curve of length at least $2r$. 
Therefore, the set $S(i)$ does not contain a family of curves
sweeping out the unit disk.  In other words, the restriction of
$a(1,2)$ to $S(i)$ vanishes.

According to Lusternik-Schnirelmann theory, the class $a(1,2)^p$
vanishes on the union $\cup_{i=1}^p S(i)$.  (For a reference on
Lusternik-Schnirelmann theory, see the first chapter of \cite{LS}.  
The result about the cup powers is proven on pages 2-3.)
Therefore, if $F$ is a
family of cycles that detects $a(1,2)^p$, then $F$ must contain a
cycle $C$ which is not in any $S(i)$.  In other words, for each
$i$, $C \cap B_i$ has length at least $r$.  Since the disks $B_i$
are disjoint, the total length of $C$ must be at least $p r \sim
p^{1/2}$.

To prove the lower bounds in Theorem 2, we need an analogue of
Lusternik-Schnirelmann theory involving Steenrod squares instead
of cup powers.  The fundamental topological fact about cup squares
that we exploited is the following
vanishing result.  If $\alpha$ is a cohomology
class in $H^*(X)$, and $V_1, V_2$ are open sets in $X$ with the
property that $\alpha |_{V_1} = 0$ and $\alpha |_{V_2} = 0$,
then the cup square $\alpha^2$
vanishes on the union $V_1 \cup V_2$.  We prove a generalization of
this vanishing result for Steenrod squares.

\begin{vanishing lemma} Let $X$ be a simplicial complex, 
and let $\alpha$ be a cohomology class in $H^p(X, \mathbb{Z}_2)$. 
Let $\pi: S^i \times X \rightarrow X$ be the projection onto the
second factor.  Suppose that $V \subset S^i \times X$ is
an open subset.  For any $\theta \in S^i$, let $V(\theta) 
\subset X$ be the set $\{ x \in X | (\theta, x) \in
V \}$. Let $P[V] = \cap_{\theta \in S^i} [V(\theta) \cup V(-
\theta)]$.  Under these assumptions, if $\pi^*(\alpha)$ vanishes
on $V$, then $Sq_i \alpha$ vanishes on $P[V]$.
\end{vanishing lemma}

We now try to describe the proof of the lower bounds in Theorem 2,
comparing each step to what happened in Theorem 1.  Suppose first that
we have a family of cycles that detects $a(k,n)^2$.  Suppose that $P$
is any hyperplane through the origin, and let $B_1$ and $B_2$ denote
the two components of $B^n(1) - P$.  By using Lusternik-Schnirelmann
theory, we can find a cycle in our
family that meets both $B_1$ and $B_2$ in a substantial volume.  
Now suppose instead
that we have a family of cycles that detects $Sq_1 a(k,n)$.  If we take
a hyperplane $P$ through the origin, the corresponding statement is
false.  Instead, we have to take a 1-parameter family of hyperplanes through the origin,
parametrized by a copy of $\mathbb{RP}^1 \subset \mathbb{RP}^{n-1}$.  For
one of these hyperplanes $P$, we can find a cycle in our family that meets
each half of $B^n(1) - P$ in a substantial volume.

To deal with a class like $Sq_i^P a(k,n)$, we have to iterate the procedure above.
As a warmup, suppose that a family detects $a(k,n)^4$.  Then we cut the ball up into 4 pieces in a
2-step process as follows.  First we pick a hyperplane $P$ that cuts
the ball into two pieces, $B_1$ and $B_2$.  Then we pick a hyperplane $P_1$ that cuts
$B_1$ into two pieces, $B_{11}$ and $B_{12}$.  Similary, we pick a hyperplane $P_2$ that
cuts $B_2$ into two pieces.  We end up with four pieces: $B_{11}, B_{12}, B_{21},$ and $B_{22}$.  Using
Lusternik-Schnirelmann theory, we can find a cycle in
our family that meets each of the four pieces in a substantial volume.  Using this argument,
we get a lower bound for $\mathbb{V}(a(k,n)^4)$.  
The lower bound that we get depends on our choice of planes.  To get
the best lower bound, we want the pieces $B_{ij}$ to be as thick as possible.

Finally suppose a family detects $Sq_1^2 a(k,n)$.  We cut the ball into 4 pieces in
a 2-step process.  First, we pick a 1-parameter family of hyperplanes making a linear copy of
$\mathbb{RP}^1 \subset \mathbb{RP}^{n-1}$.  Then the ball is cut into two pieces, $B_1$ and $B_2$,
along one of these hyperplanes, but we don't get to choose which hyperplane.  Now, for $B_1$, we get
to choose a (possibly different) 1-parameter family of hyperplanes, $\mathbb{RP}^1 
\subset \mathbb{RP}^{n-1}$.  Then the set $B_1$ is cut into two pieces, $B_{11}$ and $B_{12}$, along
one of the hyperplanes in the family, but again we don't get to control which one.  Then we do
the same for $B_2$.  At the end of the cutting process we have divided the ball into four
convex sets.  Using the Steenrod square vanishing lemma, we can prove that our family contains a cycle
that meets each of the four convex sets in a substantial volume.  Notice that we didn't get to choose
the four sets, as we did in the case of $a(k,n)^4$ - we only got to choose the sequence of
1-parameter families of hyperplanes.  Applying the Steenrod vanishing lemma, we get a lower 
bound for $\mathbb{V}(Sq_1^2 a(k,n))$, depending on the shapes of the four pieces.
In order to get the best lower bound, we want to choose our sequence of 1-parameter families
of hyperplanes in order to guarantee that the four pieces are as thick as possible.

The proof of the lower bounds in Theorem 2 combines a topology argument based on the Steenrod
vanishing lemma with a geometry argument estimating the sizes of the pieces that appear in
the above construction.  This proof is the longest and hardest part of the paper.

We now return to the upper bounds in Theorems 1 and 2.  Some of
the upper bounds can be proven using families of real algebraic
cycles, as we discussed above.  For example, the real algebraic
hypersurfaces can be used to prove the upper bounds in Theorem 1
when $k = n-1$.  In general, to prove the upper bounds, we
construct some new families of cycles.

In order to explain the idea, we illustrate the construction in
the case $k=1, n=2$.  We recall Example 2 above, the family of
sets of p vertical lines.  This family of cycles detects the
class $a(1,2)^p$, but it has maximal length on the order of $p$. 
Our
construction is a way to modify the family in Example 2, cutting
out excess length to produce a new family of curves with length
on the order of $p^{1/2}$.

Let $L$ be a lattice in $\mathbb{R}^2$ with side-length
$p^{-1/2}$.  We define a map $\Psi$ which squeezes most of
$\mathbb{R}^2$ into the 1-skeleton of $L$.  Our map $\Psi$ is
periodic, so it suffices to define it on a single square $Q$ of the
lattice $L$.  The boundary of $Q$ is contained in the 1-skeleton
of $L$, and our map is the identity on the boundary of $Q$.  Let
$Q(\epsilon)$ denote a square with the same center as $Q$ but with
side-length $\epsilon p^{-1/2}$.  The map $\Psi$ takes
$Q(\epsilon)$ linearly onto $Q$.  Finally, the map $\Psi$
retracts $Q - Q(\epsilon)$ into the boundary of $Q$.

We now use the map $\Psi$ to ``bend'' the cycles in Example 2. 
First we rotate Example 2, so that each cycle consists of 
$p$ parallel lines at a generic angle.  We then apply the map
$\Psi$ to get a new family.  Each cycle in the new family is
a union of at most $p$ 1-cycles $\Psi(L_1) + ... + \Psi(L_p)$,
where $L_i$ is a line at our generic angle.  The effect of the
map $\Psi$ on a line $L$ is illustrated below.

\includegraphics{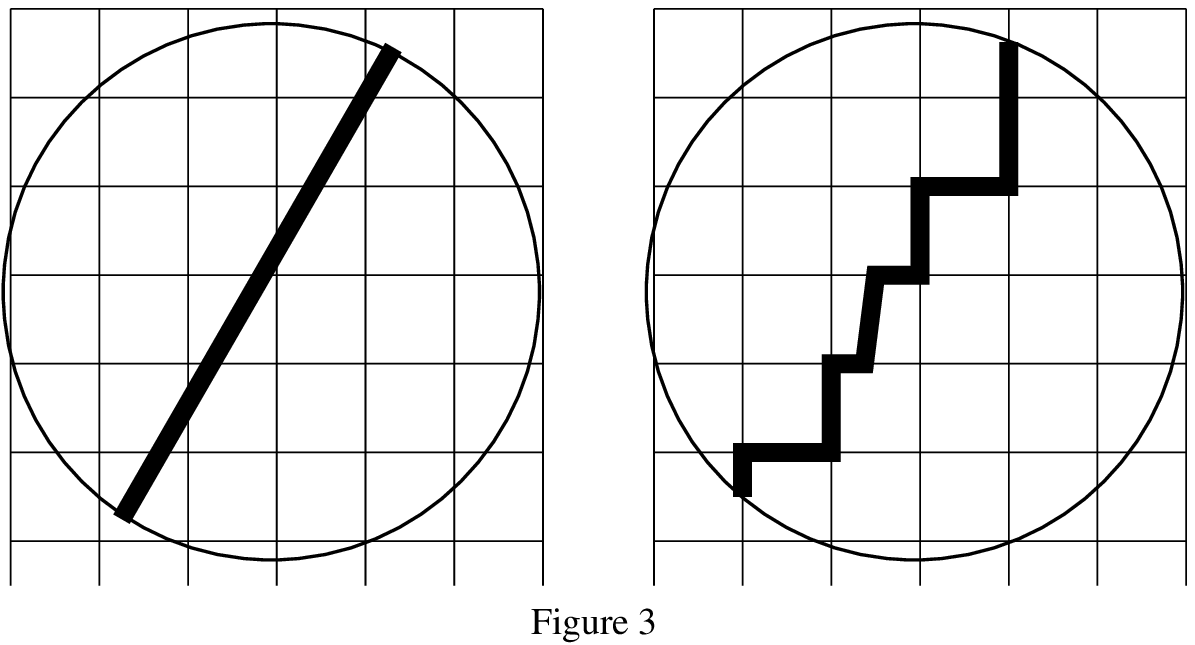}

\noindent On the left side of Figure 3, we see the intersection
of the unit disk with a line $L$.  On the right side, we see the
intersection of the unit disk with $\Psi(L)$.  The scale of the
lattice in Figure 3 corresponds to $p \sim 30$, so the reader
should imagine performing the operation above on 30 parallel lines.

The point of replacing $L_i$ by $\Psi(L_i)$ is that the curve 
$\Psi(L_i)$ is contained mainly in
the 1-skeleton of the lattice $L$.  Each cycle $\Psi(L_i)$ lies in 
the 1-skeleton of $L$ except possibly for one segment of length at 
most $2 p^{-1/2}$.  The cycles $\Psi(L_i)$, ($i = 1, ... ,
p$) overlap a great deal.  Because we are working with mod 2
cycles, we can cancel the overlaps and reduce the total length.

The sum $\sum_{i=1}^p \Psi(L_i)$ lies inside of the 1-skeleton of $L$
except for at most $p$ segments of length at most $2 p^{-1/2}$.  The
total length of all these short segments is at most $2 p^{1/2}$.  On
the other hand, after cancelling all the overlaps, the portion of
our sum inside the 1-skeleton of $L$ has length at most the length
of the 1-skeleton, which is on the order of $p^{1/2}$.

By modifying this bend-and-cancel construction, we prove all
the upper bounds in Theorems 1 and 2.  To prove Theorem 2, we need
to use a sequence of mappings like $\Psi$ with lattices at different
scales.

To finish the introduction, we mention some open questions
connected with families of real algebraic cycles.
There is a general principle in
geometry/topology that algebraic varieties do a good job of
minimizing various things.  A number of examples are given in
Arnold's expository essay \cite{Ar}.  It would be interesting to know
whether the optimal minimax volumes are realized by
families of real algebraic cycles.  

For example, we know that the family $F(d)$ of all
degree d real algebraic planar curves detects the cohomology
class $a(1,2)^{D(d)}$ for $D(d) = (1/2) (d^2 + 3d)$.  Is it
true that the minimax volume $\mathbb{V}(a(1,2)^{D(d)})$
is the maximal length of a degree d curve intersected
with the unit disk?

It would also be interesting to know whether we could prove all
the upper bounds in Theorems 1 and 2 using families of
real algebraic cycles.  By work of Lawson and Lam,
the topology of the space of real algebraic k-cycles in $\mathbb{R}^n$ 
is known.  (The complex case was done by Lawson in \cite{L}.  The real
case was done by Lam in his thesis \cite{La}.  The work is also described
in Lawson's expository article \cite{L2}, page 93.)  The space
of real algebraic k-cycles in $\mathbb{R}^n$ has $\pi_i$ equal to 
$\mathbb{Z}_2$ for $i = n-k$ and $\pi_i$ equal to zero otherwise.  
In particular, their work implies that every cohomology class of the 
form $a(k,n)^p$ or of the form $Sq_0^{Q_0} ... Sq_{n-k-1}^{Q_{n-k-1}} a(k,n)$ can
be detected by a family of real algebraic cycles.  It would be
interesting to know the smallest degree $d$ so that a given cohomology
class can be detected by a family of algebraic cycles of degree
at most $d$.  We call this degree $\mathbb{D}(\alpha)$.

An algebraic k-cycle of degree d meets the unit n-ball in volume
at most $C(k,n) d$.  Therefore we get a lower bound $\mathbb{D}(\alpha)
\ge c(k,n) \mathbb{V}(\alpha)$.  (Applying Theorem 1 and 2, we then
get lower bounds for $\mathbb{D}(a(k,n)^p)$ and $\mathbb{D}(Sq^Q a(k,n))$.)
It would be interesting to know if there is a converse bound
$\mathbb{D}(\alpha) \le C(k,n) \mathbb{V}(\alpha)$.  If it exists,
such a bound would verify the philosophy that algebraic objects
are efficient at solving geometric problems.  Our results show that
this converse inequality holds in a couple of cases including
the case $\alpha = a(n-1,n)^p$.

\vskip3pt The paper is organized as follows.  In Section 1, we
state the problem precisely.  In Section 2, we give Gromov's
lower bound for $\mathbb{V}(a(k,n))$ using the isoperimetric
inequality.  In Section 3, we prove lower bounds for
$\mathbb{V}(a(k,n)^p)$ by combining Lusternik-Schnirelmann theory
with the result in Section 2.  In Section 4, we prove lower
bounds for $\mathbb{V}(Sq^Q a(k,n))$ using our Vanishing Lemma
for Steenrod squares.  In Section 5, we construct families of
cycles using the bend-and-cancel construction, proving all the
upper bounds in Theorems 1 and 2.  In Section 6, for context, we
discuss families of algebraic cycles.  

The paper ends with three appendices.  In Appendix 1 we give a
more standard definition of the space of flat cycles.  In Appendix 2,
we describe the limited known results for families of integral
cycles.  In Appendix 3, we discuss the analogous problem in
Riemannian manifolds.

\vskip3pt

Notation: Unless otherwise indicated, all homology and cohomology
groups have coefficient group $\mathbb{Z}_2$.

\vskip3pt

{\it Acknowledgements.} \hskip3pt I would like to thank Misha Gromov for
suggesting this problem to me.  I believe that the problem first
appeared at the end of his essay \cite{G2}.  I would like to thank
the referee for pointing out to me the paper \cite{G3}.  Finally, I 
would like to thank my
thesis advisor Tom Mrowka for his help and support.  I started
working on this material as a graduate student as a different
approach to my thesis problem on area-contracting maps.  The
connection to area-contracting maps is discussed in Appendix 3.

\section{The space of flat cycles}

In this section, we set up our problem precisely.  In particular,
we define the space of mod 2 flat cycles $Z(k,n)$ and the
fundamental cohomology class $a(k,n)$.

(Our definition is different from the standard definition of the
space of flat cycles.  In Appendix 1, we recall the standard
definition and prove that they agree.)

A mod 2 Lipschitz k-chain in the unit n-ball is a finite sum
$\sum a_i f_i$, where $a_i \in \mathbb{Z}_2$, and $f_i$ is a
Lipschitz map from the standard k-simplex to the closed unit
n-ball.  We let $I_B(k,n)$ denote the space of mod 2 Lipschitz
k-chains in the unit n-ball.  We let $I_{\partial B}(k,n)$ denote
the space of mod 2 Lipschitz k-chains in the boundary of the unit
ball.  (Each map $f_i$ is a Lipschitz map from the k-simplex to
$\partial B$.)  Each of these spaces is a vector space over
$\mathbb{Z}_2$.  We define the space of relative k-chains,
$I_{rel}(k,n)$, to be the quotient $I_B(k,n) / I_{\partial
B}(k,n)$.

Next we define boundaries.  The boundary of a Lipschitz k-chain
is defined in the usual way from singular homology theory.  It gives
a boundary map $\partial: I_B(k,n) \rightarrow I_B(k-1,n)$, and a
boundary map $\partial: I_{\partial B}(k,n) \rightarrow I_{\partial B}
(k-1,n)$.  Hence we get a boundary map between the quotients,
$I_{rel}(k,n) \rightarrow I_{rel}(k-1,n)$.  This map makes the relative
Lipschitz k-chains into a complex.  The cycles in this complex
are called relative k-cycles.  The set of relative k-cycles is
denoted by $Z_{rel}(k,n)$.  A relative k-cycle can be represented by
a k-chain whose boundary lies in $\partial B$.  The homology of the chain
complex $I_{rel}$ is the relative homology of the unit n-ball: 
$H_k(B, \partial B, \mathbb{Z}_2)$.

Next we define volume.  For any k-chain $C \in I_B(k,n)$, we
define the volume of $C = \sum a_i f_i$ to be $\sum |a_i|
\textrm{ Vol }(f_i^* Euc)$.  Here $|a_i|$ is equal to 1 if $a_i =
1$ mod 2 and zero otherwise; and $f_i^*(Euc)$ is the induced
metric on the k-simplex.  Because $f_i$ is Lipschitz, this metric
is well-defined almost everywhere and belongs to $L^\infty$,
which is enough to define its volume.  Now we define the volume
of a relative k-chain $C \in I_{rel}(k,n)$ to be the infimal
volume of any absolute k-chain $C' \in I_B(k,n)$ whose projection
to $I_{rel}(k,n)$ is $C$.

We define a metric on $Z_{rel}(k,n)$ which measures how much
area it takes to span the difference of two cycles.  If $C_1$ and
$C_2$ are two relative cycles, we define the area-distance
between $C_1$ and $C_2$ to be the infimal volume of any relative
(k+1)-chain $D \in I_{rel}(k+1,n)$ with $\partial D = C_1 - C_2$.

The following figure illustrates two 1-cycles that are close
together in terms of the area-distance.

\includegraphics{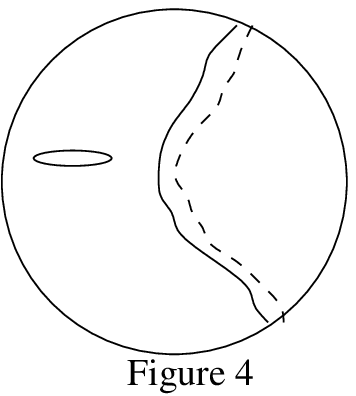}

\noindent The solid line and the dashed line each represent a
relative 1-cycle in the unit disk.  The solid cycle has two
components and the dashed cycle has only one component. One
component of the solid cycle lies far from the dashed cycle, so
the two cycles are far apart in the Hausdorff topology.
Nevertheless, the region between them, consisting of a strip and a
small ellipse, has small area, and so the cycles are close
together in the area-distance.

The area-distance from $C_1$ to
$C_2$ can be zero.  This happens if $C_1$ and $C_2$ are the same
geometric object parametrized in different ways, or if $C_2$ is
obtained from $C_1$ by adding a degenerate map $f_i$ which takes the
k-simplex into a (k-1)-dimensional surface.

We say that two relative cycles $C_1$ and $C_2$ are equivalent if
the area-distance between them is zero.  The set of equivalence
classes of relative cycles is a metric space, where the metric
is the area-distance.  We define $Z(k,n)$ to be the completion
of this metric space.  We say that the volume of a cycle $C \in
Z(k,n)$ is less than $V$ if there is a sequence of relative
Lipschitz cycles $C_i$ with $C_i \rightarrow C$ and the volume of
each $C_i$ less than $V$.  We sometimes denote the volume of $C$
by $|C|$.

By a family $F$ of k-cycles parametrized by a space $X$, we mean
a continuous map $F: X \rightarrow Z(k,n)$.  In this paper, we
will always assume that $X$ is a simplicial complex.  We say that
a family $F$ detects a cohomology class $\alpha \in H^*(Z(k,n))$
if $F^*(\alpha)$ is non-zero in $H^*(X)$.  We let
$\mathbb{F}(\alpha)$ denote the set of all families of cycles
that detect $\alpha$.  We will say that a cycle $C \in Z(k,n)$
belongs to $F$ if $C = F(x)$ for some $x \in X$.  If $C$ belongs
to $F$, we will write $C \in F$.  We define the maximal
volume of the family $F$ to be $\sup_{C \in F} \textrm{ Volume }(C)$.

Now we formally define the minimax volumes.  Let $\alpha$ be a
cohomology class in $H^*(Z(k,n))$.  We define
$\mathbb{V}(\alpha)$ by the following minimax formula.

$$\mathbb{V}(\alpha) = \inf_{F \in \mathbb{F}(\alpha)} \sup_{C
\in F} \textrm{ Volume }(C).$$

Next we construct the fundamental cohomology class $a(k,n) \in
H^{n-k}(Z(k,n), \mathbb{Z}_2)$.  The construction follows
Almgren's original construction in \cite{A}. Informally, our task
is to define what it means for a family of relative k-cycles to
sweep out the unit ball.  Morally, an i-dimensional family of
k-cycles can be glued together to form a (k+i)-cycle, but this is
not literally true.  We now give a construction that takes an
i-dimensional family of cycles and gives a (k+i)-cycle that, in
some sense, is a small perturbation of the family.

The basic object that we introduce to do the construction is called
a complex of cycles.  A complex of cycles is a discrete approximation
to a continuous family of cycles.  A complex of k-cycles is
parametrized by a polyhedral complex $X$.  For each i-face $A$ of $X$,
the complex associates a (k+i)-dimensional mod 2 Lipschitz relative chain
$C(A)$.  These chains have to fit together in the following sense.
If $\partial A = \sum B_i$, then $\partial C(A) = \sum C(B_i)$.
In particular, this formula implies that for each vertex $v$ of $X$,
$C(v)$ is a relative k-cycle.

We include a figure illustrating a complex of 1-cycles in the
unit disk.  In the left half of the figure is
the parameter space of our complex.  It consists of three
0-simplices, p, q, and r, together with two 1-simplices, E and F. 
In the right half of the figure, we see the corresponding cycles
and chains.  For each 0-simplex, there is a corresponding 1-cycle
in the unit disk.  For each 1-simplex there is a corresponding
2-chain in the unit disk.  The figure is supposed to show
geometrically the way these cycles and chains fit together. 

\includegraphics{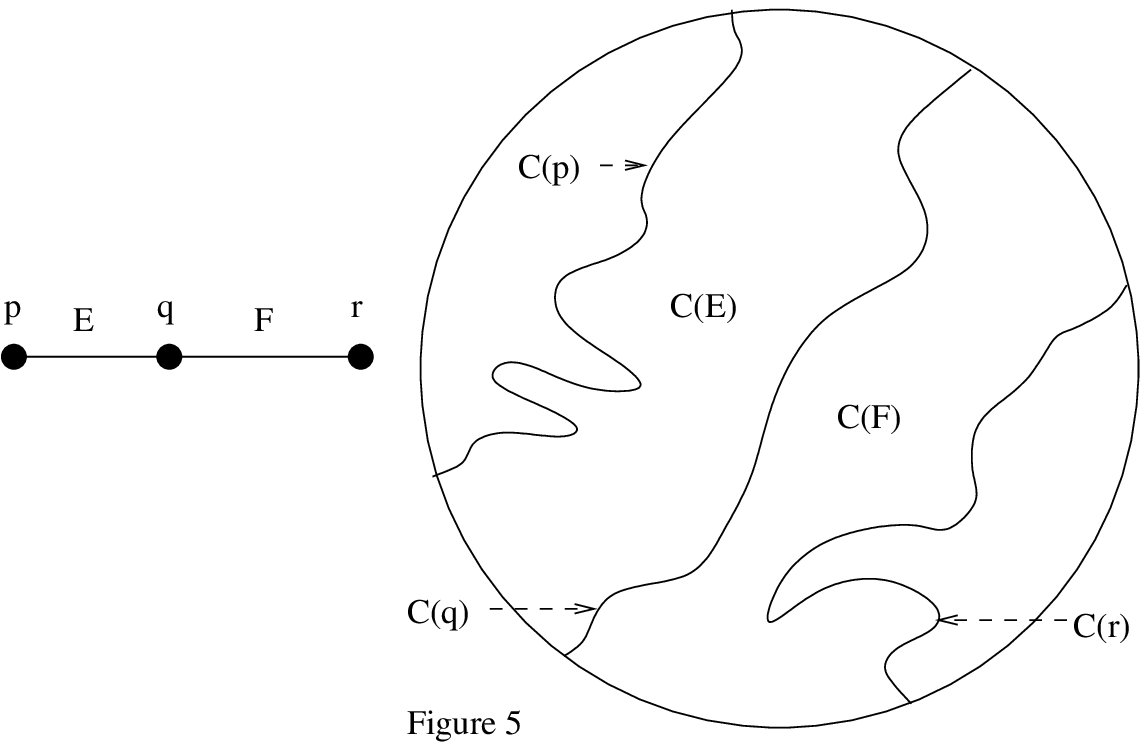}

We can think of $C$ algebraically as a chain map between two
chain complexes. The first chain complex is associated to $X$. 
It has i-chains consisting of sums $c_i A_i$, where $c_i \in
\mathbb{Z}_2$ and $A_i$ is an i-face of $X$. The second chain
complex is the complex of relative Lipschitz chains in the unit
n-ball, $I_{rel}(k,n)$.  A complex of k-cycles $C$ is a chain map
from the first chain complex to the second one with shift k.  (We
remark that $I_{rel}(k,n)$ makes perfectly good sense for $k >
n$.)

The chain map $C$ induces a map from the simplicial homology of
$X$ to the homology of $I_{rel}(k,n)$, which is the relative
homology of the unit n-ball.  Since the chain map has shift $k$,
we get a map from $H_i(X)$ to $H_{k+i}(B, \partial B)$.  We call
this map the gluing homomorphism $G$.

There is a version of homotopy for complexes of cycles.  If $C$
is a complex of cycles parametrized by $X \times [0,1]$ (with the
natural polyhedral structure), then we call $C$ a homotopy.  If
$C$ restricted to $X \times \{ 0 \}$ is $C_0$ and if $C$
restricted to $X \times \{ 1 \}$ is $C_1$, then we say that $C$
is a homotopy from $C_0$ to $C_1$.  The gluing homomorphism is
homotopy invariant.  From the algebraic point of view, $C$ is a
chain homotopy between the chain maps $C_0$ and $C_1$.

Following Almgren, we next explain how to approximate a
continuous family of k-cycles by a discrete complex of cycles. 
The main tool in the construction is the Federer-Fleming
isoperimetric inequality.

\begin{reftheorem} (Federer, Fleming) There is a constant $C(n)$
so that the following holds.  Suppose that $k < n$, and that $C$
is a mod 2 relative Lipschitz k-cycle in the unit n-ball.  Then
$C$ is the boundary of a (k+1)-chain $D$, with $|D| < C(n)
|C|^{\frac{k+1}{k}}$.
\end{reftheorem}

In addition to the isoperimetric inequality, we also use the
following basic facts.  If $C$ is a mod 2 Lipschitz relative
n-cycle in the unit n-ball, and if the volume of $C$ is less than
the volume of the unit ball, then $C$ bounds an (n+1)-chain $D$. 
Any k-cycle of dimension greater than n automatically bounds a
(k+1)-chain.

Let $F$ be a family of k-cycles parametrized by $X$.  We take a
fine triangulation of $X$.  We pick a small number $\delta > 0$. 
For each vertex $v$ of the triangulation, we choose a mod 2
Lipschitz cycle $C(v)$ so that the area-distance from $C(v)$ to
$F(v)$ is at most $\delta$ and so that the volume of $C(v)$ is at
most $|F(v)| + \delta$.

Now, since the triangulation is fine, we may assume that if $v_1$
and $v_2$ are neighboring vertices, then the area-distance
between $C(v_1)$ and $C(v_2)$ is less than $3 \delta$.  By
definition, this means that there is a mod 2 relative Lipschitz
$(k+1)$-chain $D$ with boundary $C_1 - C_2$, and with $|D| <
3 \delta$.

Let E denote the edge from $v_1$ to $v_2$, so that $\partial E =
v_1 - v_2$.  We define $C(E) = D$.  We repeat this operation for
every edge of the triangulation of $X$.  For each edge $E$ with
boundary $v_1 - v_2$, $C(E)$ is a (k+1)-chain with boundary
$C(v_1) - C(v_2)$ and volume at most $3 \delta$.

We continue this procedure inductively.  For each i-dimensional
simplex $\Delta^i$ in $X$, we define a (k+i)-chain $C(\Delta^i)$
with the following properties.  If the boundary of the simplex
$\Delta^i$ is equal to $\sum_j \Delta_j^{i-1}$, then the boundary
of $C(\Delta^i)$ is equal to $\sum_j C(\Delta_j^{i-1})$ in 
$I_{rel}(k,n)$.  Moreover, $C(\Delta^i)$ has volume less than
$C(n) \delta$.  If $i \le n-k$, we can choose such chains by
using the isoperimetric inequality.  Provided that $\delta$ is
sufficiently small, we can also choose such a chain for $i = n-k
+ 1$, because the boundary we are trying to fill has n-volume
less than $C(n) \delta$. For $i > n-k+1$, we can automatically
find such a chain.

The complex of cycles $C$ was not canonical.  If $\delta > 0$ is
sufficiently small, however, the complex of cycles $C$ is
well-defined up to homotopy.  To see this, let $C_0$ and $C_1$ be
two possible choices of chain map following the construction
above. Divide $X \times [0,1]$ into cells given by $\Delta \times
\{0\}, \Delta \times \{1 \},$ and $\Delta \times [0,1]$, where
$\Delta$ varies over the triangulation of $X$.  Now define
$C(\Delta \times
\{0\}) = C_0(\Delta)$ and $C(\Delta \times \{1\}) = C_1(\Delta)$. 
Suppose that $\Delta^i$ is an i-simplex in $X$, and that the
boundary of $\Delta^i \times [0,1]$ is equal to $\Delta^i \times
\{ 1 \} - \Delta^i \times \{ 0 \} +
\sum_j \Delta_j^{i-1} \times [0,1]$. 
We have to define $C(\Delta^i \times [0,1])$ to be a
(k+i+1)-chain with boundary $C(\Delta^i \times \{1\}) -
C(\Delta^i \times \{0\}) + \sum_j C(\Delta_j^{i-1} \times
[0,1])$.  We proceed inductively, beginning with $i=0$.  In this
case, we have to find a (k+1)-chain spanning $C(v,1) - C(v,0)$. 
But both $C(v,1)$ and $C(v,0)$ are within $\delta$ of $F(v)$, and
so the area-distance between them is less than $2 \delta$.  Hence
we can choose a (k+1)-chain with the given boundary and with
volume less than $2 \delta$.  If $i <n-k$, we define these chains
using the isoperimetric inequality.  Provided that $\delta$ is
sufficiently small, we can define the chains for $i = n-k$ since
we are trying to fill a boundary with n-volume less than $C(n)
\delta$.  When $i > n-k$, we can define the chains
automatically.

A similar argument shows that if $C_0$ is a complex of cycles
chosen with respect to a triangulation $T_0$ of $X$, and if $C_1$
is a complex of cycles chosen with respect to a refined
triangulation $T_1$, then the gluing maps of $C_0$ and $C_1$ agree.

The gluing map of the complex $C$ is therefore defined canonically.
The gluing map gives a homomorphism $G: H_i(X) \rightarrow H_{k+i}(B,
\partial B)$.  In particular, taking $i = n-k$, we get a homomorphism
$G: H_{n-k}(X) \rightarrow \mathbb{Z}_2$.  The gluing homomorphism
makes precise the idea of a family of cycles ``sweeping out'' the unit
ball: we say that $F$ sweeps out the unit ball if the corresponding
gluing homomorphism $G: H_{n-k} \rightarrow \mathbb{Z}_2$ is non-trivial.

The gluing map gives a homomorphism from $H_{n-k}(Z(k,n))
\rightarrow \mathbb{Z}_2$.  Let $z$ be an (n-k)-cycle in $Z(k,n)$.  
We can think of $z$ as a family of cycles, and then
consider $G([z])$ in $H_n(B, \partial B) =
\mathbb{Z}_2$.  If $z'$ is homologous to $z$, then let $y$ be an
(n-k+1)-chain with boundary $z - z'$.  The gluing map $G$ extends
to a map from $H_i(y) \rightarrow H_{k+i}(B, \partial B)$ and so
$G([z]) = G([z'])$.  Hence $G$ gives a homomorphism from
$H_{n-k}(Z(k,n))$ to $\mathbb{Z}_2$.  By the universal
coefficient theorem, this homomorphism determines a cohomology
class $a(k,n) \in H^{n-k}(Z(k,n))$.

For example, suppose that $F$ is the family of parallel k-planes
parametrized by $t \in [-1, 1]^{n-k}$, with $F(t)$ equal to the
intersection of the unit n-ball with the plane $\mathbb{R}^k \times
\{ t\}$.  For each $t$ in the boundary of $[-1, 1]^{n-k}$, $F(t)$
is equal to the empty cycle in $Z(k,n)$.  Therefore, we can think
of $F$ as a map of pairs $(I^{n-k}, \partial I^{n-k})
\rightarrow (Z(k,n), *)$, where $I^{n-k}$ denotes $[-1,1]^{n-k}$
and $*$ denotes the empty cycle in $Z(k,n)$.  We can pick a fine
triangulation of $[-1, 1]^{n-k}$ and define $C(\Delta)$ to be
$\mathbb{R}^k \times \Delta$ intersected with the unit n-ball.  The
fundamental homology class $h$ of $[-1,1]^{n-k}$ relative to its
boundary is given by the sum of all top-dimensional simplices of
our triangulation.  Following the definition, we see that $G(h)$
is equal to the fundamental homology class of $(B, \partial B)$. 
Therefore, $F^*(a(k,n))$ is equal to the fundamental cohomology
class of $(I^{n-k}, \partial I^{n-k})$.

\section{Lower bounds based on the isoperimetric inequality}

In this section we give a lower bound for $\mathbb{V}(a(k,n))$.
This result is due to Almgren and later Gromov gave a simpler
proof which we copy here.  The lower bound for $\mathbb{V}(a(k,n))$
is the basis for all the other lower bounds in the paper.  In fact,
for technical reasons, the later lower bounds are based on
a slightly modified minimax volume $\mathbb{V}^+(a(k,n))$ which we
introduce below.

First we give Gromov's elementary lower bound for
$\mathbb{V}(a(k,n))$.  (Gromov's argument appears on page 134 of \cite{G}.)

\begin{prop} (Gromov) The minimax volume $\mathbb{V}(a(k,n)) \ge
c(n)$, for a dimensional constant $c(n) > 0$.
\end{prop}

\proof Suppose that $F: X \rightarrow Z(k,n)$ is a family of
k-cycles in the unit n-ball, and that each cycle in $F$ has
volume less than $\epsilon$, a small number that we will choose
later.  We have to prove that $F$ does not detect the class
$a(k,n)$.  Equivalently, we have to prove that the gluing
homomorphism $G: H_{n-k}(X) \rightarrow H_n(B, \partial B)$ is
zero.

The gluing homomorphism is induced by a complex of cycles $C$. 
We first recall the construction of $C$.  We pick a fine
triangulation of $X$.  For each vertex $v$ of the triangulation,
we pick a Lipschitz cycle $C(v)$ with distance at most $\delta$
from $F(v)$ and with volume at most $\epsilon + \delta$.  Next
for each edge $E$ with boundary $v_1 - v_2$, we choose a
(k+1)-chain $C(E)$ with volume at most $3 \delta$ and with boundary
$C(v_1) - C(v_2)$.  Here $\delta$ is a number smaller than
$\epsilon$ which we can make as small as we like by choosing a
sufficiently fine triangulation of $X$.  Then we proceed to
define $C(\Delta^i)$ with volume at most $C(n) \delta$ for each
$i$-simplex $\Delta^i$ in our triangulation.

We will prove that our complex of cycles is homotopic to the zero
complex.  By definition, that means that we will construct a
complex of cycles $\bar C$ defined on the product $X \times
[0,1]$, with $\bar C(\Delta \times \{ 0 \}) = C(\Delta)$ and
$\bar C(\Delta \times \{ 1 \})$ equal to the empty cycle. 
Because of the homotopy, the gluing homomorphism associated to
$C$ is equal to zero.

Now we do the construction.  We have already defined $\bar C$ on
each face of the form $\Delta \times \{ 0 \}$ and $\Delta \times
\{ 1 \}$, so it remains to define our complex $\bar C$ on faces
of the form $\Delta \times (0,1)$.  We do this inductively,
beginning with 1-faces of the form $v \times (0,1)$, where $v$ is
a vertex in our triangulation of $X$.

Using the isoperimetric inequality, each cycle $C(v)$ can be
filled by a (k+1)-chain of volume less than $C(n)
\epsilon$.  For each vertex $v$ of the triangulation, we define
$\bar C(v \times (0,1))$ to be such a filling.  We have $\partial
\bar C(v \times (0,1)) = C(v) = \bar C(v \times \{0 \}) - \bar
C(v \times 1)$, so this choice obeys the boundary equation for
complexes of cycles.

Now we inductively define $\bar C(\Delta^i \times (0,1))$ for $i
\ge 1$, so that it has volume at most $C(n) \epsilon$.  By
induction, we can assume that we have already defined $\bar C$ on
all lower-dimensional skeleta.  In particular, we have already
defined $\bar C$ on the boundary of our cell $\bar C(\Delta^i
\times (0,1))$.  The boundary is associated to a
$(k+i)$-dimensional cycle of total volume at most $C(n)
\epsilon$.  Now by the isoperimetric inequality, we can find a
filling of this cycle with total volume at most $C(n) \epsilon$,
and we define $\bar C(\Delta^i \times (0,1))$ to be this filling.  
The special case that $i = n-k$
deserves a further remark.  In this case, the boundary of
$\Delta^i \times (0,1)$ is associated to an n-cycle in the unit
n-ball of total volume at most $C(n) \epsilon$.  Because
$\epsilon$ is sufficiently small, this n-cycle must have degree
0.  Because it has degree zero, it admits a filling.  The filling
is an (n+1)-chain, which automatically has zero volume. \endproof

The key step in the proof above was to find a chain $\bar C(v
\times (0,1))$ with boundary $C(v)$ and with (k+1)-volume at most
$C(n) \epsilon$.  Therefore, our proof shows that a family of
cycles detecting $a(k,n)$ must contain a cycle of large filling
volume.  We recall that the filling volume of a Lipschitz k-cycle
$C$ is the smallest (k+1)-volume of any relative chain $D$ with
$\partial D = C$.  The filling volume is clearly continuous in
the area-distance, and so it defines a continuous function on $Z(k,n)$.

\begin{prop} Let $F$ be a family of cycles that detects $a(k,n)$. 
Then $F$ contains a cycle with filling volume at least $c(n)$.
\end{prop}

\proof The proof is essentially the same.  We proceed by
contradiction, assuming that each $F(v)$ has filling volume at
most $\epsilon$.  Then each $C(v)$ has filling volume at most
$\epsilon + \delta$.  Then we can construct $\bar C(v \times (0,1))$
with volume at most $\epsilon + \delta$.  The rest of the proof
goes as above. \endproof

For technical reasons, we now introduce a minor variant of the
minimax volume $\mathbb{V}(\alpha)$.  We say that 
$\mathbb{V}^+(\alpha)$ is at least $V$ if, for any
family of cycles $F: X \rightarrow Z(k,n)$, if we let $S \subset
X$ denote the subset of cycles with volume at most $V$, then
$F^*(\alpha)$ vanishes on {\it an open neighborhood} of $S$. 

The difference between $\mathbb{V}$ and $\mathbb{V}^+$ is as
follows.  If $V < \mathbb{V}(\alpha)$, and if $S \subset X$
denotes the subset of cycles with volume at most $V$, then
$F^*(\alpha)$ vanishes on $S$, but it's not obvious whether it
vanishes on a neighborhood of $S$.  The set $S$ is compact.  If
it happens to have a neighborhood that retracts onto it, then
$F^*(\alpha)$ vanishes on that neighborhood, but the set $S$ may
be a very nasty compact set.  We will need to use open sets later
in the paper because we use theorems from algebraic topology that
hold for open covers but which don't hold for covers by arbitrary
compact sets.

It follows from the definition that $\mathbb{V}^+(\alpha) \le
\mathbb{V}(\alpha)$.  We will prove lower bounds for
$\mathbb{V}^+(\alpha)$, and these bounds immediately
imply lower bounds for $\mathbb{V}(\alpha)$.  
Gromov's method also gives a lower bound for $\mathbb{V}^+(a(k,n))$.

\begin{prop} (Gromov) The minimax volume $\mathbb{V}^+(a(k,n)) \ge
c(n)$, for a dimensional constant $c(n) > 0$.
\end{prop}

\proof Let $F: X \rightarrow Z(k,n)$ be a family
of k-cycles.  Let $S \subset X$ be the subset of cycles with
volume at most $V$.  If $V$ is sufficiently small, we need to
prove that $F^*(a(k,n))$ vanishes on a neighborhood of $S$.

Even though every cycle in $S$ has volume at most $V$, we have no
way to bound the volumes of cycles in any neighborhood of $S$,
because the volume function is not continuous.  The situation
improves by looking at the filling volume.  Because of the
Federer-Fleming isoperimetric inequality, every cycle in $S$ has
filling volume at most $C(n) V^{\frac{k+1}{k}}$.  But the filling volume is a
continuous function on $Z(k,n)$.  Therefore, we can choose a
neighborhood of $S$ in which every cycle has filling volume at
most $2 C(n) V^{\frac{k+1}{k}}$.  Now if $V$ is sufficiently small, then we can
apply Proposition 2.2 to conclude that $F^*(a(k,n))$ vanishes on
this neighborhood.  \endproof

To finish this section, we make some historical and expository
remarks about the lower bound for $\mathbb{V}(a(k,n))$.  This
material is not needed in the proofs of the theorems.

The minimax volume $\mathbb{V}(a(n-1,n))$ can be bounded below using
the isoperimetric inequality.  If we have a 1-parameter family of
hypersurfaces sweeping out the unit n-ball, then one of these surfaces
must divide the ball in half by volume.  (This remains true even
if the surfaces are not embedded.)  In other words, we can choose
a surface $C$ in our family so that $B - C = U_1 \cup U_2$, where
$U_i$ is a union of connected components of $B - C$, and $|U_1| =
|U_2| = (1/2) |B|$.  We can then divide the boundary of $B$ into
two pieces $\partial B_1, \partial B_2$, according to which open
set they border.  We reorder the sets so that $|\partial B_1| \le
(1/2) |\partial B|$.  Now the boundary of $U_1$ is contained in
$\partial B_1 \cup C$, and so it has total volume at most $(1/2)
|\partial B| + |C|$.  By the isoperimetric inequality,
$(1/2)^{\frac{n-1}{n}} |B|^{\frac{n-1}{n}} = |U_1|^{\frac{n-1}{n}}
\le C(n) |\partial U_1| \le C(n) [(1/2) |\partial B| + |C|]$.  In
the isoperimetric inequality, the sharp constant is given by $B$,
so $C(n) |\partial B| = |B|^{\frac{n-1}{n}}$.  Rearranging we get
the following inequality.

$$|C| \ge [(1/2)^{\frac{n-1}{n}} - (1/2)] |\partial B|.$$

\noindent This approach can be made sharp by using the
sharp isoperimetric inequality for relative cycles in the
unit n-ball.

To my knowledge, the first person to consider the analogous
problem with $k < n-1$ was Almgren.  In \cite{A2}, he proved the
following sharp theorem for families of cycles sweeping out the
n-sphere.

\begin{reftheorem} (Almgren) Let $F$ be a family of k-cycles
sweeping out the unit n-sphere.  Then the maximal volume of $F$
is at least the volume of the unit k-sphere.
\end{reftheorem}

Taking a bilipschitz embedding of the unit n-ball into the unit
n-sphere, Almgren's theorem implies a non-sharp lower bound
for $\mathbb{V}(a(k,n))$.  Almgren's lower bound is still
better than the one in Proposition 2.1.  Since the bilipschitz
constant of the embedding can be taken independent of $n$ it
follows that $\mathbb{V}(a(k,n)) \ge c(k)$, independent of $n$.

Almgren's proof involves the theory of
varifolds, which he invented more or less for this purpose. 
Unfortunately, the proof was never published, but similar
arguments appear in Pitts's book \cite{P}.  We might expect to
get better estimates for $\mathbb{V}(a(k,n))$ by applying
Almgren's method instead of by applying his theorem about spheres.  In
\cite{G} (page 135), Gromov mentions that Almgren's
techniques imply $\mathbb{V}(a(k,n)) = \omega_k$.

Recently, in \cite{G3}, Gromov gave a sharp estimate for the
k-dimensional waist of the unit n-sphere.  This result is
similar to Almgren's theorem, but the technical details of the
statement are different.  The paper uses topological methods
instead of minimal surfaces.

\section{Lower bounds based on Lusternik-Schnirelmann theory}

We now estimate the minimax volume of the cohomology class $a(k,n)^p$,
proving the lower bounds from Theorem 1.  The results in this
section are due to Gromov, and they appeared (with different notation)
in Section 8 of \cite{G3}.

\begin{Theorem 1} (Lower bounds) The minimax volume
$\mathbb{V}(a(k,n)^p)$ is greater than $c(n) p^{\frac{n-k}{n}}$.
\end{Theorem 1}

\proof The proof is based on Lusternik-Schnirelmann theory.  Let
$B_1, ..., B_p$ be disjoint balls inside the unit ball, each with
radius $r = (1/4) p^{-1/n}$. (It is not hard to find disjoint
balls with this radius.  A detailed argument is given in the
proof of the cup product theorem below.)

Let $F: X \rightarrow Z(k,n)$ be a family of k-cycles that
detects $a(k,n)^p$.  Let $S(i) \subset X$ be the subset of cycles
$x \in X$ so that $F(x) \cap B_i$ has volume at most $(1/2)
\mathbb{V}^+(a(k,n)) r^k$.  By the definition of
$\mathbb{V}^+(a(k,n))$ and a scaling argument, $F^*(a(k,n))$
vanishes on an open neighborhood of $S(i)$.  By
Lusternik-Schnirelmann theory, the cohomology class
$F^*(a(k,n))^p$ vanishes on a neighborhood of $\cup_{i=1}^p
S(i)$.  (For a proof, see pages 2-3 of \cite{LS}.)  
Since $F$ detects $a(k,n)^p$, there must be a cycle $C$
in $F$ that does not lie in $S(i)$ for any $i$.

By the definition of $S(i)$, the intersection $C \cap B_i$ has
volume greater than $(1 / 2) \mathbb{V}^+(a(k,n)) r^k$. 
Since the balls $B_i$ are disjoint, the cycle $C$ has total
volume greater than $(1 / 2) p \mathbb{V}^+(a(k,n)) r^k$.  According to
Proposition 2.3, $\mathbb{V}^+(a(k,n)) \ge c(n)$, and $r$ was
defined to be $(1/4) p^{-1/n}$.  Plugging in, we see that the 
volume of $C$ is at least $c(n) p^{\frac{n-k}{n}}$. \endproof

The proof of Theorem 1 allows us to estimate the minimax volume
$\mathbb{V}^+(\alpha)$ for a cup product $\alpha = \alpha_1 \cup
... \cup \alpha_p$ in terms of the minimax volumes
$\mathbb{V}^+(\alpha_i)$.

\begin{theoremcp} Suppose that $\alpha = \alpha_1 \cup ... \cup
\alpha_P$ is a cohomology class in $H^*(Z(k,n))$.  Then the
minimax volume of $\alpha$ obeys the following inequality.

$$\mathbb{V}^+(\alpha) \ge 4^{-k} [\sum_{i=1}^P
\mathbb{V}^+(\alpha_i)^{\frac{n}{n-k}}]^{\frac{n-k}{n}}.$$
\end{theoremcp}

\proof Suppose that $B_1, ... B_P$ are disjoint balls inside the
unit ball, with radii $r_1$, ..., $r_P$.  We will choose
particular balls later.

Let $F: X \rightarrow Z(k,n)$ be a family of k-cycles.  Let
$\epsilon > 0$ be any small number.  Let $S
\subset X$ be the subset of cycles with volume at most $V = (1 - \epsilon)
4^{-k} [\sum_{i=1}^P
\mathbb{V}^+(\alpha_i)^{\frac{n}{n-k}}]^{\frac{n-k}{n}}$.  We
have to prove that $F^*(\alpha)$ vanishes on a neighborhood of
$S$.

We define $S(i) \subset X$ to be the subset of cycles $x \in X$
so that $F(x) \cap B_i$ has volume at most $(1 - \epsilon) r_i^k
\mathbb{V}^+(\alpha_i)$.  By scaling and the definiton of
$\mathbb{V}^+$, it follows that $F^*(\alpha_i)$ vanishes on a
neighborhood of $S(i)$.  By Lusternik-Schnirelmann theory, the
class $F^*(\alpha)$ vanishes on a neighborhood of $\cup_{i=1}^P
S(i)$.  Since the balls $B_i$ are disjoint, it follows that this
union contains every cycle of volume at most $\sum_{i=1}^P (1 - \epsilon)
r_i^k \mathbb{V}^+(\alpha_i)$.

In order to finish the proof, we have to choose disjoint balls
$B_i$ so that $V = \sum (1 - \epsilon) r_i^k
\mathbb{V}^+(\alpha_i)$.  We claim that we can find disjoint
balls in the unit ball with radii $r_i$ provided that
$\sum_{i=1}^P r_i^n \le 4^{-n}$.  To find the disjoint balls, we
order the radii so that $r_1 \ge ... \ge r_P$. Then we choose a
point $p_1$ inside the ball of radius $1/2$ around the origin. 
Since $\sum r_i^n \le 4^{-n}$, $r_1 \le 1/4$, and so the ball
$B(p_1, r_1)$ is contained in the unit ball.  Now we proceed
inductively.  We suppose that we have chosen $p_1, ..., p_{j-1}$
in $B(0, 1/2)$ so that the balls $B(p_i, r_i)$ are disjoint.  If
$\omega_n$ denotes the volume of the unit n-ball, then the volume
of $B(0,1/2)$ is $2^{-n} \omega_n$. On the other hand, the
volume of $\cup_{i=1}^{j-1} B(p_i, 2 r_i)$ is $\sum_{i=1}^{j-1}
(2 r_i)^n \omega_n < 2^{-n} \omega_n$. Therefore, we can choose
a point $p_j$ in $B(0, 1/2)$ but not in the union
$\cup_{i=1}^{j-1} B(p_i, 2 r_i)$.  Since $r_j \le r_i$ for $i <
j$, the ball $B(p_j, r_j)$ is disjoint from $\cup_{i=1}^{j-1}
B(p_i, r_i)$.  Continuing to choose balls in this way proves the
claim.

Finally, we choose $r_i$ subject to $\sum r_i^n \le 4^{-n}$ in
order to maximize $\sum (1 - \epsilon) r_i^k \mathbb{V}^+(\alpha_i)$.  The
maximum value is $V = (1 - \epsilon) 4^{-k} [\sum_i
\mathbb{V}^+(\alpha_i)^{\frac{n}{n-k}}]^{\frac{n-k}{n}}$. It is obtained by
using the following values of $r_i$:

$$r_i = (1/4) \mathbb{V}^+(\alpha_i)^{\frac{1}{n-k}} [\sum_{i=1}^P
\mathbb{V}^+(\alpha_i)^{\frac{n}{n-k}}]^{-1/n}.$$

\endproof

Remark: The cup product theorem holds for cohomology with any
coefficients - i.e. if we consider cohomology classes in
$H^*(Z(k,n), A)$ for any coefficient group $A$.  It also holds if
we replace $Z(k,n)$ by the space of integral Lipschitz cycles. 
For more comments on integral cycles, see Appendix 2.

\section{Lower bounds based on Steenrod squares}

In this section, we prove the lower bounds in Theorem 2.  These lower
bounds involve Steenrod squares.  For background on
Steenrod squares, consult \cite{H}.  We use one piece of notation
which is not in \cite{H}.  We define $Sq_i: H^p \rightarrow H^{2p-i}$
to be $Sq^{p-i}$.  We define $Sq_i^2 \alpha$ to be $Sq_i [Sq_i \alpha]$.
In a similar way, we inductively define cohomology classes of the form
$Sq_0^{Q_0} ... Sq_{n-k-1}^{Q_{n-k-1}} \alpha$.

The topological input to the proof is a vanishing lemma
for Steenrod squares.  This lemma will generalize the fact that if $V_1, V_2$
are open sets with $\alpha|_{V_1} = 0$ and $\alpha|_{V_2} = 0$, then 
the cup square $\alpha^2$ vanishes on $V_1 \cup V_2$.

In order to state our lemma, we set up some vocabulary.  
Let $X$ be a simplicial complex.  
Let $\pi: S^i \times X \rightarrow X$ denote
the projection onto the second factor.  Let $V$ be a subset of
$S^i \times X$.  We let $V(\theta)$ denote the set $\{ x \in X |
(\theta, x) \in V \}$.  The main definition that we need to 
state our vanishing lemma is the following.

$$P[V] := \cap_{\theta \in S^i} [V(\theta) \cup V(- \theta)].$$

\noindent We say a few words to describe $P[V]$.  First of all,
$P[V]$ is a subset of $X$.  For each $x \in X$, we examine the
intesection of $V$ with $S^i \times \{ x \}$. 
Roughly, we include $x$ in $P[V]$ if this intersection is sufficiently
large.   Let $V(x)$ denote
the set $\{ \theta \in S^i | (\theta, x) \in V \}$.  
We include $x$ in $P[V]$ if, for each $\theta$,
either $\theta$ or $- \theta$ belongs to $V(x)$.  In other
words, we map $V(x)$ to $\mathbb{RP}^i$ using the 
standard covering $S^i \rightarrow \mathbb{RP}^i$, and we include $x$
in $P[V]$ if the map is surjective.

\begin{vanishing lemma} Suppose $V$ is an open subset of $S^i \times X$, and
that $\alpha$ is a cohomology class in $H^*(X)$.
If $\pi^* \alpha$ vanishes on $V$, then
$Sq_i \alpha $ vanishes on $P[V]$.  Incidentally, $P[V]$ is open.
\end{vanishing lemma}

\proof The main part of the proof concerns Steenrod squares and
the construction of a homotopy.  The final result also requires
a little point-set topology, which we put at the end.  Therefore,
we begin by proving a slightly weaker statement.  Suppose that 
$K$ is any compact subset of $V$.  We will first prove that $Sq_i \alpha$
vanishes on a neighborhood of $P[K]$.

We begin by recalling a construction of Steenrod squares
given in \cite{H} on pages 502-504.

We can assume that $X$ is connected.  Let $x_0$ be a basepoint of
$X$.  Let $X \wedge X$ denote the smash product of $X$ with
itself.  The group $\mathbb{Z}_2$ acts on $S^i \times X \wedge X$
by sending $(\theta, x_1, x_2)$ to $(- \theta, x_2, x_1)$.  We
call the quotient space $\Gamma^i X$.

By abuse of notation, we use $x_0$ to denote the basepoint of $X
\wedge X$ as well as the basepoint of $X$.  The action of
$\mathbb{Z}_2$ on $S^i \times X \wedge X$
sends $(\theta, x_0)$ to $(-\theta, x_0)$.  Therefore, the image
of $S^i \times \{ x_0 \}$ in the quotient $\Gamma^i X$ is a copy
of $\mathbb{RP}^i$.  We let $\Lambda^i X$ denote the quotient of
$\Gamma^i X$ obtained by collapsing this $\mathbb{RP}^i$ to a
point.  

Let $\Delta_0: S^i \times X \rightarrow S^i \times X \wedge X$ be
the diagonal map $\Delta_0(\theta, x) = (\theta, x, x)$.  There
is an action of $\mathbb{Z}_2$ on each of these spaces.  The
action on $S^i \times X$ sends $(\theta, x)$ to $(-\theta, x)$. 
The action on $S^i \times X \wedge X$ is the one we mentioned
above, which sends $(\theta, x_1, x_2)$ to $(- \theta, x_2,
x_1)$.  The map $\Delta_0$ is equivariant with respect to this
action, and so it descends to a map between the quotients.  If we
then collapse $\Gamma^i X$ to $\Lambda^i X$, we get a map
$\Delta: \mathbb{RP}^i \times X \rightarrow \Lambda^i X$.

Let $p$ be the dimension of the cohomology class $\alpha$.
Then the cohomology class $\alpha$ can be represented by a map $f$
from $X$ to the Eilenberg-Maclane space $K(\mathbb{Z}_2, p)$.  We
abbreviate this space as $K(p)$.  The construction of $\Lambda^i
X$ is functorial, so we get a map $\Lambda^i f: \Lambda^i X
\rightarrow \Lambda^i K(p)$.  The definition of Steenrod squares
uses a certain cohomology class $\lambda \in H^{2p}(\Lambda^i
K(p))$.  The class $\lambda$ is defined in \cite{H}, page 504.

The Steenrod squares of $\alpha$ can now be defined in terms of
$\lambda,
\Lambda^i f,$ and $\Delta$.  Recall that by the Kunneth theorem,
$H^*(\mathbb{RP}^i \times X)$ is the tensor product
$H^*(\mathbb{RP}^i) \otimes H^*(X)$.  Let $\omega$ be the
generator of $H^1(\mathbb{RP}^i)$.  Then the following formula
can be taken as a definition of Steenrod squares.

$$\Delta^* (\Lambda^i f^* \lambda) = \sum_{j=0}^i \omega^j
\otimes Sq_j \alpha.$$

To prove our vanishing theorem, we will construct a homotopy of
$\Lambda^i f \circ \Delta$ which maps a neighborhood of
$\mathbb{RP}^i \times P[K]$ to the basepoint of
$\Lambda^i K(p)$.  From this homotopy, it follows that
$\sum_{j=0}^i \omega^j \otimes Sq_j \alpha$ vanishes on
a nieghborhood of $RP^i \times P[K]$.  The Kunneth
theorem then implies that $Sq_i \alpha$ vanishes on a
neighborhood of $P[K]$.

Recall that $\pi: S^i \times X \rightarrow X$ is the projection
onto the second factor.  The map $f \circ \pi$ from 
$S^i \times X$ to $K(p)$ induces the
cohomology class $\pi^* \alpha$.  We know by hypothesis that
$\pi^* \alpha$ vanishes on an open set $V$ containing the compact
set $K \subset
S^i \times X$. Therefore, the map $f \circ \pi$, restricted to
$V$, is null-homotopic.  By the homotopy extension
theorem, we can homotope $f \circ \pi$ to a map that sends a
smaller neighborhood $U \supset K$ to the basepoint of $K(p)$.  More
formally, we have a continuous family $F_t: S^i
\times X \rightarrow K(p)$ where $F_0(\theta, x) =
f(x)$ and where $F_1$ maps $U$ to the
basepoint of $K(p)$.

Using this homotopy of $f$, we now construct a homotopy of
$\Lambda^i f \circ \Delta$.  We define a family of maps $G_t$
from $S^i \times X$ to $S^i \times K(p) \wedge K(p)$ by the
following formula.

$$G_t(\theta, x) = (\theta, F_t(\theta, x), F_t(- \theta, x)).$$

The maps $G_t$ are equivariant with respect to the $\mathbb{Z}_2$
actions defined above.  At time 0, we have $G_0(\theta, x) =
(\theta, f(x), f(x))$.  Let $*$ denote the basepoint of $K(p)
\wedge K(p)$.  At time 1, we have $G_1(\theta, x) \in
S^i \times \{*\}$ if either $(\theta, x)$ or $(- \theta, x)$ is
in $U$.  Because $G_t$ is equivariant, it descends to a family of
maps from $\mathbb{RP}^i \times X$ to $\Gamma^i K(p)$.  Composing
with the quotient map from $\Gamma^i K(p)$ to $\Lambda^i K(p)$,
we get a family of maps $H_t$ from $\mathbb{RP}^i \times X$ to
$\Lambda^i K(p)$.  The map
$H_0$ is equal to $\Lambda^i f \circ \Delta$.  The map $H_1$
sends a point $(\pm \theta, x) \in \mathbb{RP}^i
\times X$ to the basepoint of $\Lambda^i K(p)$ if either $(\theta,
x)$ or $(- \theta, x)$ lies in $U$.  In particular,
$H_1$ maps $\mathbb{RP}^i \times P[U]$ to the basepoint of
$\Lambda^i K(p)$.

We know that $U$ is an open neighborhood of $K$.  
Since $K$ is compact, it follows that
$U$ contains the $\epsilon$-neighborhood of $K$ for some
$\epsilon > 0$.  Therefore, $U(\theta)$ contains the $\epsilon$
neighborhood of $K(\theta)$ for every $\theta$.  Hence
$U(\theta) \cup U(- \theta)$ contains the $\epsilon$-neighborhood
of $K(\theta) \cup K(- \theta)$.  We conclude that $P[U]$ contains
the $\epsilon$-neighborhood of $P[K]$.  In particular, the map
$H_1$ sends a neighborhood of $\mathbb{RP}^i \times P[K]$ to the
basepoint of $\Lambda^i K(p)$.  We conclude that $Sq_i \alpha$
vanishes on a neighborhood of $P[K]$.

Now let $K_1 \subset K_2 \subset ...$ be an exhaustion of $V$
by compact sets.  Let $V_j$ be an open neighborhood of $K_j$
contained in $K_{j+1}$.  Clearly $P[K_1] \subset P[K_2] \subset ...$
We will check that the union of $P[K_j]$ is equal to $P[V]$.
Suppose that $x \in P[V]$.  Let $\tau: S^i \rightarrow \mathbb{RP}^i$ denote
the standard covering, and recall that $V(x) = \{ \theta \in S^i | (\theta, x) \in V \}$.  
Since $x \in P[V]$, $\tau(V(x))
 = \mathbb{RP}^i$.  Since $V$ is equal to the union of $V_j$, $V(x)$ is
equal to the union of the open sets $V_j(x)$.  Hence $\tau(V_j(x))$ is an
open covering of $\mathbb{RP}^i$ by infinitely many sets.  Since
$\mathbb{RP}^i$ is compact, there is a finite subcovering, and we conclude
that one of the sets $\tau(V_j(x))$ covers $\mathbb{RP}^i$.
In other words, $x$ lies in $P[V_j]$, which lies in $P[K_{j+1}]$.  Therefore
$P[V] = \cup P[K_j]$ and so $Sq_i \alpha$ vanishes on $P[V]$.

Incidentally, the argument above implies that $P[V]$ is open.
Suppose that $x$ lies in $P[V]$.  By the last paragraph, we know that $\tau(K_j(x))$
is equal to $\mathbb{RP}^i$ if $j$ is big enough.  Let $B_r(x)$ denote the ball
around $x$ in $X$ of radius $r$.  Since $K_j$ is a compact subset of $V$, the set
$K_j(x) \times B_r(x) \subset V$ for a small positive $r > 0$.  It follows
that $B_r(x)$ is contained in $P[V]$.  Hence $P[V]$ is open. \endproof

Using this vanishing lemma in place of the Lusternik-Schnirelmann
theory, we will prove lower bounds for $\mathbb{V}^+(\beta)$ for
cohomology classes $\beta \in H^*(Z(k,n))$ defined using Steenrod
squares.

\begin{theoremst} For each $\epsilon > 0$, there is a constant
$c(n, \epsilon) > 0$ that makes the following estimate hold.
If $\alpha \in H^*(Z(k,n))$ is a cohomology class, $i$ is an integer
in the range $0 \le i \le n-k-1$, and $P > 0$ is an integer, then
$\mathbb{V}^+(Sq_i^P \alpha)$ obeys the lower bound below.

$$\mathbb{V}^+(Sq_i^P \alpha) \ge c(n,\epsilon) (2 -
\epsilon)^{\frac{n-i-k}{n-i} P} \mathbb{V}^+(\alpha).$$
\end{theoremst}

I believe that this theorem should also hold with $\epsilon = 0$,
but I don't know how to prove it.  When $i=0$, the theorem
follows immediately from the cup product theorem.
We have $\mathbb{V}^+(Sq_0^P \alpha) =
\mathbb{V}^+(\alpha^{2^P})$.  By the cup product theorem, $\mathbb{V}^+
(\alpha^{2^P}) \ge c(n) [2^P]^{\frac{n-k}{n}} \mathbb{V}^+(\alpha) = c(n)
2^{\frac{n-k}{n} P} \mathbb{V}^+(\alpha)$.  
So we see that when $i=0$, the theorem holds even when $\epsilon = 0$.

Given the Steenrod tower theorem, we can quickly prove the lower bounds
in Theorem 2.

\begin{Theorem 2} (Lower bounds) For each $\epsilon > 0$, there is a constant
$c(n, \epsilon) > 0$ so that the following estimate holds.

$$c(n, \epsilon) \prod_{i=0}^{n-k-1} (2 -
\epsilon)^{\frac{n-k-i}{n-i} Q_i} \le \mathbb{V}(Sq_0^{Q_0} ...
Sq_{n-k-1}^{Q_{n-k-1}} a(k,n)).$$

\end{Theorem 2}

\proof By the definition of $\mathbb{V}^+$, we know the following inequality.

$\mathbb{V}(Sq_0^{Q_0} ... Sq_{n-k-1}^{Q_{n-k-1}} a(k,n)) \ge
\mathbb{V}^+(Sq_0^{Q_0} ... Sq_{n-k-1}^{Q_{n-k-1}} a(k,n))$.

As described above, we apply the cup product theorem to deal with the
$Sq_0$ term.

$\mathbb{V}^+(Sq_0^{Q_0} ... Sq_{n-k-1}^{Q_{n-k-1}} a(k,n)) \ge
c(n) [2^{Q_0}]^{\frac{n-k}{n}} \mathbb{V}^+(Sq_1^{Q_1} ... Sq_{n-k-1}^{Q_{n-k-1}} a(k,n)).$

Now we apply the Steenrod tower theorem $n-k-1$ times to deal with the other
Steenrod squares.

$\ge c(n, \epsilon) 2^{\frac{n-k}{n} Q_0} \prod_{i=1}^{n-k-1} (2 - \epsilon)^{\frac{n-i-k}{n-i} Q_i}
\mathbb{V}^+(a(k,n)).$

Finally, we use Proposition 2.3 to bound $\mathbb{V}^+(a(k,n)) \ge c(n)$.

$\ge c(n, \epsilon) 2^{\frac{n-k}{n} Q_0} \prod_{i=1}^{n-k-1} (2 - \epsilon)^{\frac{n-i-k}{n-i} Q_i}.$

This inequality is slightly stronger than the one we had to prove. \endproof

Now we take up the proof of the Steenrod Tower Theorem.

\proof Throughout the proof, we fix $i, k,$ and $n$.  The case $i=0$ was
proven above, so we assume $i \ge 1$.

The main idea of the proof is to inductively use the vanishing lemma
from the beginning of this section.
In the proof of the cup product theorem, in order to lower bound
$\mathbb{V}^+(\alpha^p)$, we cut the unit ball into $p$ disjoint pieces.
To adapt that argument to Steenrod squares, we again need to cut
the ball into pieces.  But this time, instead of cutting the ball
into pieces in one way, we need a high-dimensional family
of different ways of cutting the ball into pieces.  We formalize
this idea as a ``pyramid of subsets of the unit ball''.

We let $\mathbb{O}$ be the set of open sets in
$\mathbb{R}^n$ equipped with the Hausdorff
topology.

An $i$-pyramid of open sets of height $P$ 
consists of the following data.  For each
integer $0 \le p \le P$, we have a map $U_p$ from $(S^i)^{P-p}$
to $\mathbb{O}$.  (By convention, the space $(S^i)^0$ is a single point.)
The maps $U_p$ must obey two rules.

Rule 1. $U_p(\theta_1, ..., \theta_{P-p}) \subset
U_{p+1}(\theta_1, ..., \theta_{P-p-1})$.

Rule 2. $U_p(\theta_1, ..., \theta_{P-p-1}, \theta_{P-p})$ and
$U_p(\theta_1, ..., \theta_{P-p-1}, - \theta_{P-p})$ are
disjoint.

We call $U_P$ the top level of the i-pyramid.  It consists of 
one large open set, $U_P(*)$, where $*$ denotes the one point
in the space $(S^i)^0$.  We say that $U$ is an i-pyramid of open
sets in the unit ball if $U_P(*)$ is contained in the unit ball.
We call the family
of sets $U_0$ the bottom level of the pyramid.

Rule 1 describes how the various levels of an i-pyramid are
related to each other.  Roughly it says that the open sets in the
top levels are the largest and that the open sets in the bottom
levels are the smallest.  Rule 2 forces various open sets in an
i-pyramid to be disjoint.  For example, it implies that the $2^P$
subsets $U_0(\pm \theta_1, ..., \pm \theta_P)$ are disjoint.  In
Figure 6, we illustrate some of the sets in an i-pyramid of open
sets with height $P=2$.

\vskip5pt

\includegraphics{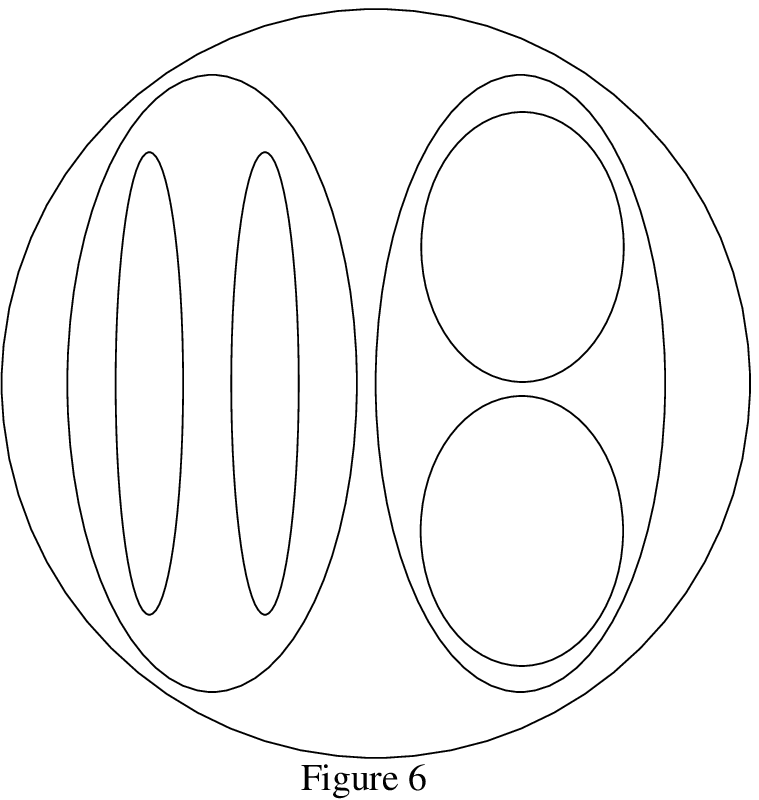}

\noindent The largest circle is $U_2(*)$, 
the unit ball in $\mathbb{R}^2$.  The
medium-sized ellipse on the left is $U_1(\theta)$ and the medium
ellipse on the right is $U_1(- \theta)$, for some point $\theta
\in S^i$. The two small thin ellipses on the left are
$U_0(\theta, \pm \phi)$ for some point $\phi \in S^i$.  The two
small round ellipses on the right are $U_0(- \theta, \pm \psi)$
for some point $\psi \in S^i$.

Let $U$ be an i-pyramid of open sets in the unit ball.  We are going
to define a measure of the thickness of various sets in $U$.  
If $A$ is any open set, we define $Rad[A]$ to be the
largest radius $R$ of any ball $B(x, R) \subset A$.  For any
$\theta \in (S^i)^P$, we define $T_0(\theta) = (1/2) Rad[U_0(\theta)]^k$.
Notice that $T_0$ depends on the value of $k$.  Roughly speaking, $T_0$
is measuring the k-dimensional thickness of the set $U_0(\theta)$.
Our function $T_0$ is defined on the bottom layer
of the pyramid.  We inductively define thickness functions on higher
layers of the pyramid, $T_p: (S^i)^{P-p} \rightarrow \mathbb{R}$,
by the following formula.

$$T_p(\theta_1, ..., \theta_{P-p}) = \inf_{\phi \in S^i}
T_{p-1}(\theta_1, ..., \theta_{P-p}, \phi) + T_{p-1}(\theta_1,
..., \theta_{P-p}, - \phi).$$

\noindent Following the induction, we ultimately define a function
$T_P: (S^i)^0 \rightarrow \mathbb{R}$.  The space $(S^i)^0$ consists
of a single point $*$, and so $T_P$ has a single value $T_P(*)$, 
which we abbreviate as $T$.  We call $T$ the k-thickness of the i-pyramid
$U$.

We can control minimax volumes of Steenrod squares in terms of the
thickness of pyramids according to the following lemma.  This lemma
makes up the first half of the proof of the Steenrod Tower Theorem.

\begin{ML} As above, we fix $n$ and $1 \le k \le n-1$ and $i$ in the range
$1 \le i \le n-k-1$.  Suppose that there is an i-pyramid $U$ of open
sets in the unit ball with height $P$ and k-thickness $T$.  Suppose further
that every set $U_0(\theta)$ is convex.
Then for
every class $\alpha \in H^*(Z(k,n))$, the following inequality holds.

$$\mathbb{V}^+(Sq_i^P \alpha) \ge T \mathbb{V}^+(\alpha).$$
\end{ML}

\proof We will make an inductive argument using the Vanishing Lemma
from the beginning of this section.  The base for our induction
is given by the following lemma.

\begin{lemma} Let $F: X \rightarrow Z(k,n)$ be a family of cycles.  
We define a subset of small cycles $\mathbb{S}_0
\subset (S^i)^P \times X$, by saying that $(\theta, x) \in \mathbb{S}_0$ if
the restriction of $F(x)$ to $U_0(\theta)$ has volume at 
most $T_0(\theta) \mathbb{V}^+(\alpha)$.  Let
$\pi_0$ be the projection from $(S^i)^P \times X$ to the second factor.  Then
$\pi_0^* (F^* \alpha)$ vanishes on a neighborhood of $\mathbb{S}_0$.
\end{lemma}

\proof We can find a ball $B(\theta) \subset U_0(\theta)$ with radius
$Rad[U_0(\theta)]$.  We would like these balls to vary continuously
with $\theta$.  Because the sets $U_0(\theta)$ are convex, we can arrange this
for slightly smaller balls.  If $\delta > 0$, we can find a continuous family of balls
$B(\theta) \subset U_0(\theta)$, where $B(\theta)$ has radius at least $(1 - \delta) Rad[U_0(\theta)]$.  
We pick a fine triangulation of $(S^i)^P$.  For each vertex of the triangulation, we choose a center $c(v) \subset
U_0(v)$, so that the ball of radius $Rad[U_0(v)]$ around $c(v)$ is contained in $U_0(v)$.
Now we extend $c$ to a piecewise-linear function on $(S^i)^P$.  For a sufficiently fine triangulation,
we claim that the ball around $c(\theta)$ of radius $(1 - \delta) Rad[U_0(\theta)]$ lies in
$U_0(\theta)$ for every $\theta$.  This last step uses the convexity of $U_0(\theta)$.
The point $\theta$ lies in some simplex of our triangulation with vertices $v_1, ..., v_N$.
Since the triangulation is fine, we can assume that $U_0(\theta)$ is close to $U_0(v_i)$
in the Hausdorff topology, and so the ball around $c(v_i)$ of radius $(1- \delta) Rad[U_0(\theta)]$
lies in $U_0(\theta)$.  Now, the set of all points $c$ so that the ball of radius $(1 - \delta)
Rad[U_0(\theta)]$ lies in $U_0(\theta)$ is a convex set.  Since $c(\theta)$ is a convex
combination of $c(v_i)$, it follows that the ball around $c(\theta)$ of radius $(1-\delta)
Rad[U_0(\theta)]$ lies in $U_0(\theta)$.

We now define a map $F_0$ from $(S^i)^P \times X$ to $Z(k,n)$.  The cycle $F(\theta, x)$ is
defined by taking the cycle $F(x)$, restricting it to the ball $B(\theta)$, and then rescaling
the ball to get a cycle in the unit ball.  If $(\theta, x) \in \mathbb{S}_0$, then the
rescaled cycle has volume at most $(1- \delta)^{-k} (1/2) \mathbb{V}^+(\alpha) < \mathbb{V}^+(\alpha)$.
By the definition of $\mathbb{V}^+$, the cohomology class $F_0^* \alpha$ vanishes
on a neighborhood of $\mathbb{S}_0$.

Finally, we check that $F_0$ is homotopic to $F \circ \pi_0$.  We construct
a homotopy of our family of balls $B(\theta, t)$ so that $B(\theta,0) = B(\theta)$,
$B(\theta, 1)$ is the unit ball, and each ball $B(\theta, t)$ is contained in
the unit ball.  To construct the family, we first move all the centers $c(\theta)$
to the center of the unit ball, and then we rescale the balls so that all the radii
are $1$.  Now we define $F_t(\theta, x)$ by taking the cycle $F(x)$, restricting to
the ball $B(\theta, t)$, and rescaling to get a cycle in $Z(k,n)$.  At $t=1$, we have
$F_t = F \circ \pi_0$, since the restriction and rescaling are both the identity.  Therefore,
$F_t$ gives a homotopy from $F_0$ to $F \circ \pi_0$.  We conclude that $\pi_0^* [F^* \alpha]$
vanishes on a neighborhood of $\mathbb{S}_0$. \endproof

Now we turn to the inductive step in the proof.  We consider
the spaces $S^{P-p} \times X$ for $p$ from $0$ to $P$.  Inside
each space, we define a subset of "small cycles" $\mathbb{S}_p
\subset S^{P-p} \times X$.  We say that $(\theta, x) \in
\mathbb{S}_p$ if the intersection of $F(x)$ with
$U_p(\theta)$ has volume at most $T_p(\theta) \mathbb{V}^+(\alpha)$.
Let $\pi_p: (S^i)^{P-p} \times X \rightarrow X$ denote the
projection onto the second factor.  We will show inductively
that the cohomology class $\pi_p^* (F^* Sq_i^p \alpha)$ vanishes
on a neighborhood of $\mathbb{S}_p$.  Lemma 4.1 proves
the base case $p = 0$.

We have to set up the problem in such a way that we can apply 
the Vanishing Lemma.  By induction, we assume that $\pi_{p-1}^*
(F^* Sq_i^{p-1} \alpha)$ vanishes on a neighborhood of 
$\mathbb{S}_{p-1}$.  We let $Y = (S^i)^{P-p} \times X$ so
that $S^i \times Y = (S^i)^{P-(p-1)} \times X$.  We let $\pi:
S^i \times Y \rightarrow Y$ denote the projection onto
the second factor.  We let $\beta = \pi_p^* (F^* Sq_i^{p-1} \alpha)$,
a cohomology class in $H^*(Y)$.  Our inductive hypothesis
tells us that $\pi^*\beta$ vanishes on a neighborhood of
$\mathbb{S}_{p-1} \subset S^i \times Y$.  The Vanishing Lemma
then implies that $Sq_i \beta$ vanishes on a neighborhood
of $P[\mathbb{S}_{p-1}]$.  Plugging in the definition of $\beta$,
we see that $Sq_i \beta = \pi_p^*(F^* Sq_i^p \alpha)$.  To
complete the induction, we only have to show that $\mathbb{S}_p$
is contained in $P[\mathbb{S}_{p-1}]$.

Suppose that $(\theta, x)$ is contained in $\mathbb{S}_p$,
where $\theta = (\theta_1, ..., \theta_{P-p})$.  By definition,
the volume of $F(x) \cap U_p(\theta)$ is at most 
$T_p(\theta) \mathbb{V}^+(\alpha)$.  Now let $\phi$ be any
point in $\mathbb{S}^i$.  By the definition of an i-pyramid
of convex sets, we know that $U_{p-1}(\theta, \phi)$ and
$U_{p-1}(\theta, - \phi)$ are disjoint subsets of $U_p(\theta)$.
Therefore, we get the following formula.

$$|F(x) \cap U_{p-1}(\theta, \phi)| + |F(x) \cap U_{p-1}(\theta, - \phi)|
\le T_p(\theta) \mathbb{V}^+(\alpha).$$

But by the definition of the thickness function $T_p$, we know that
$T_p(\theta) \le T_{p-1}(\theta, \phi) + T_{p-1}(\theta, -\phi)$.  

$$ |F(x) \cap U_{p-1}(\theta, \phi)| + |F(x) \cap U_{p-1}(\theta, - \phi)|
\le [T_{p-1}(\theta, \phi) + T_{p-1}(\theta, - \phi)] \mathbb{V}^+(\alpha).$$

We see that either $|F(x) \cap U_{p-1}(\theta,
\phi)| \le T_{p-1}(\theta, \phi) \mathbb{V}^+(\alpha)$ or $|F(x)
\cap U_{p-1}(\theta, - \phi)| \le T_{p-1}(\theta, - \phi)
\mathbb{V}^+(\alpha)$. In other words, either $(\theta, \phi, x)$ is in $\mathbb{S}_{p-1}$ or else
$(\theta, - \phi, x)$ is in $\mathbb{S}_{p-1}$.  Since this
analysis applies to every $\phi \in S^i$, we conclude that 
$\mathbb{S}_p$ is contained in $P[\mathbb{S}_{p-1}]$.  This
argument proves the inductive step.

We conclude that $\pi_P^* (F^* Sq_i^P \alpha)$ vanishes on
a neighborhood of $\mathbb{S}_P \subset (S^i)^{P-P} \times X$.
But the space $(S^i)^{P-P} \times X$ is just $X$, and the
projection $\pi_P$ is just the identity.  In other words,
$F^* Sq_i^P \alpha$ vanishes on a neighborhood of $\mathbb{S}_P$.
The set $\mathbb{S}_P$ is just the subset of $X$ where $F(x)$
has volume at most $T \mathbb{V}^+(\alpha)$.  Since this analysis
applies to any family $F$, we conclude that $\mathbb{V}^+(Sq_i^P \alpha)
\ge T \mathbb{V}^+(\alpha)$. \endproof

In order to use Main Lemma 1, we need to construct an i-pyramid of
open sets with height $P$ and estimate its k-thickness.  
Ideally we would like
to know the largest possible k-thickness for an i-pyramid of open
sets of height $P$ (where the bottom layer consists of convex sets).
This problem is a kind of max-min problem, because for each i-pyramid
the k-thickness is defined by taking a sequence of infima, and we are then
looking for the supremum over all the pyramids.  The best estimate
I know how to prove is contained in the following lemma.

\begin{ML} As above, we fix a dimension $n$, and $1 \le k \le n-1$,
and $1 \le i \le n-k-1$.  For every $P$ we will construct an i-pyramid of open sets in
the unit n-ball of height P.  Each open set in each pyramid will
be convex.  For any $\epsilon > 0$, there is a constant $c(n, \epsilon) > 0$
so that for every $P$, the pyramid of height $P$ has k-thickness
at least $c(n, \epsilon) (2 - \epsilon)^{\frac{n-i-k}{n-i} P}$.
\end{ML}

The Steenrod tower theorem follows immediately by combining the 
two main lemmas.  I believe that this lemma should also hold with
$\epsilon = 0$, but I don't know how to prove it.  If this lemma did hold
with $\epsilon = 0$, then the Steenrod tower theorem and Theorem 2 would
also hold with $\epsilon = 0$, determining $\mathbb{V}(Sq^Q a(k,n))$ up to 
a dimensional constant $C(n)$.  The lemma cannot hold if
we replace $(2 - \epsilon)$ by $(2 + \epsilon)$ because the lower bound
in Theorem 2 would then become larger than the upper bound.

\proof First we will construct an i-pyramid of open sets.  Then we will
estimate its k-thickness.

We are going to talk about some sequences of unit vectors.  First fix
$v_{1-n}, ...,v_{-1}, v_0$ a basis of orthonormal vectors.  Then, for fixed
$p$, 
consider all sequences $v_1, ..., v_{P-p}$ of unit vectors with
the property that $v_a$ is perpendicular to the previous (n-i-1)
vectors: $v_{a-1}, ..., v_{a-n+i+1}$.  In other words, any string
of (n-i) consecutive vectors in the list is orthonormal.  This
condition holds even if $a-n+i+1 < 1$, which is why we 
fixed $v_{1-n},... , v_0$.  Call the set of such
sequences of vectors $O(p)$.  (We define $O(P)$ to be a point.) 
Now there is a map $O(p-1) \rightarrow O(p)$ given by forgetting
the last vector.  The fibers of this map are each i-spheres, and
the map is a fiber bundle.

\begin{lemma} Each of the bundles $O(p-1) \rightarrow O(p)$ is
trivial.
\end{lemma}

\proof We work by induction.  Clearly the map $O(P-1) \rightarrow O(P)$
is trivial, as the base is a point.  We suppose that $O(p)
\rightarrow O(p+1)$ is trival.  To do the induction, we need to show that
the bundle $O(p-1) \rightarrow O(p)$ is trivial.  We accomplish this by
showing that $O(p-1)$ is the pullback of the bundle $O(p)$ by the map $O(p)
\rightarrow O(p+1)$.  That sentence is a bit confusing, so we illustrate it
with a diagram.  Let $\pi$ denote the map $O(p) \rightarrow O(p+1)$.

\[
\begin{CD}
\pi^*(O(p)) @>>> O(p) \\
@VVV @VV\pi V \\
O(p) @>\pi>> O(p+1)
\end{CD}
\]

We prove that the bundle $O(p-1)$ over $O(p)$ is isomorphic to the bundle
$\pi^*(O(p))$ in the diagram above.  
The fiber of $O(p-1) \rightarrow O(p)$ over a point $(v_1, ...,
v_{P-p}) \in O(p)$ is given by the unit vectors perpendicular to
$v_{P-p}, ..., v_{P-p-n+i+2}$.  The map $O(p) \rightarrow O(p+1)$
takes $(v_1, ..., v_{P-p})$ to $(v_1, ..., v_{P-p-1})$.  The
fiber of $O(p) \rightarrow O(p+1)$ over the point $(v_1, ..., v_{P-p-1})$
is the set of
unit vectors perpendicular to $v_{P-p-1}, ..., v_{P-p-n+i+1}$. 
By the definition of $O(p)$, $v_{P-p}$ is perpendicular to the
last list, and $v_{P-p-n+i+1}$ is perpendicular to the first list. 
The two fibers intersect in a cod-1 great circle which is
perpendicular to both $v_{P-p}$ and $v_{P-p-n+i+1}$.
We can map one fiber to the other by mapping $v_{P-p-n+i+1}$ to
$v_{P-p}$, and by using the identity map on the great circle of
intersection.  This defines a bundle isomorphism, showing that
$O(p-1)$ is isomorphic as a bundle to $\pi^*(O(p))$.  Since
$O(p)$ was a trivial bundle by induction, it follows that
$O(p-1)$ is also a trivial bundle. \endproof

As a corollary, we see that $O(p)$ is diffeomorphic to
$(S^i)^{P-p}$.  Moreover, we can inductively define
diffeomorphisms $\psi(p): (S^i)^{P-p} \rightarrow O(p)$, so
that the following diagram is an isomorphism of bundles.

\[
\begin{CD}
(S^i)^{P-p} @>\psi(p)>> O(p) \\
@VVV @VV\pi V \\
(S^i)^{P-p-1} @>\psi(p+1)>> O(p+1)
\end{CD}
\]

The left vertical arrow in this diagram is the map taking
$(\theta_1, ..., \theta_{P-p-1}, \theta_{P-p})$ to $(\theta_1,
..., \theta_{P-p-1})$. Since $\psi$ is an isomorphism of bundles,
if $\psi(p)(\theta_1, ...,
\theta_{P-p-1}, \theta_{P-p}) = (v_1, ..., v_{P-p-1},
v_{P-p})$, then $\psi(p)(\theta_1, ...,
\theta_{P-p-1}, - \theta_{P-p}) = (v_1, ..., v_{P-p-1},
- v_{P-p})$.

We are going to construct maps $CON_p$ associating to each point
of $O(p)$ a convex subset of the unit ball.  In particular we can
think of $CON_p$ as a map from $O(p)$ to $\mathbb{O}$, defined
for each $p$ in the range $0 \le p \le P$.  Our i-pyramid of open
sets will be given by the composition $U_p = CON_p \circ
\psi(p)$.  The two rules for the pyramid $U_p$ are equivalent to
the following two rules for the map $CON_p$.

Rule 1. $CON_p(v_1, ..., v_{P-p}) \subset CON_{p+1}(v_1, ...,
v_{P-p-1})$

Rule 2. $CON_p(v_1, ..., v_{P-p-1}, v_{P-p})$ and $CON_p(v_1,
..., v_{P-p-1}, - v_{P-p})$ are disjoint.

We will define $CON_p$ inductively.  We define $CON_P(*)$ to be
the unit ball.  Suppose we have already defined $CON_{p+1}$, and
we want to define $CON_p$.  Let
$L_0$ be the hyperplane perpendicular to $v_{P-p}$.  Let $L$ be
the translation of $L_0$ that divides $CON_{p+1}(v_1, ...,
v_{P-p-1})$ into two halves of equal volume.  We let $CON_p(v_1,
..., v_{P-p})$ be one of these two halves, chosen as follows. 
Consider a vector $v_{P-p}$ starting at a point on $L$.  The end
of the vector lies on one side of $L$, and $CON_p(v_1, ...,
v_{P-p})$ is the half of $CON_{p+1}(v_1, ..., v_{P-p-1})$ on this
side of $L$.  By definition $CON_p(v_1, ..., v_{P-p})$ is
contained in $CON_{p+1}(v_1, ..., v_{P-p-1})$.  Therefore,
Rule 1 holds.  Now if we change the sign of the last factor,
then the plane L does not change, only we take the half of
$CON_{p+1}(v_1, ..., v_{P-p-1})$ on the other side of this plane.
Hence $CON_p(v_1, ..., v_{P-p-1}, v_{P-p})$ and $CON_p(v_1, ...,
v_{P-p-1}, - v_{P-p})$ are disjoint, and so Rule 2 holds.

We include a figure to help see how this construction works.  Consider
the special case that $i=1$, $n=2$, and $P=2$.  In this
case, $O(2)$ is a point $*$, $O(1)$ is the unit circle $S^1$, 
and $O(0)$ is the torus $(S^1)^2$.  The convex set $CON_2(*)$ is
the unit disk.  We let $v$ be the unit vector $(- 1 / \sqrt2, 1 / \sqrt 2)$,
and we let $w$ be the unit vector $(1,0)$.  In Figure 7, we illustrate
the convex sets $CON_1(v), CON_0(-v,w),$ and $CON_0(-v,-w)$.  The
set $CON_1(-v)$ is the union of $CON_0(-v,w)$ and $CON_0(-v,-w)$.

\includegraphics{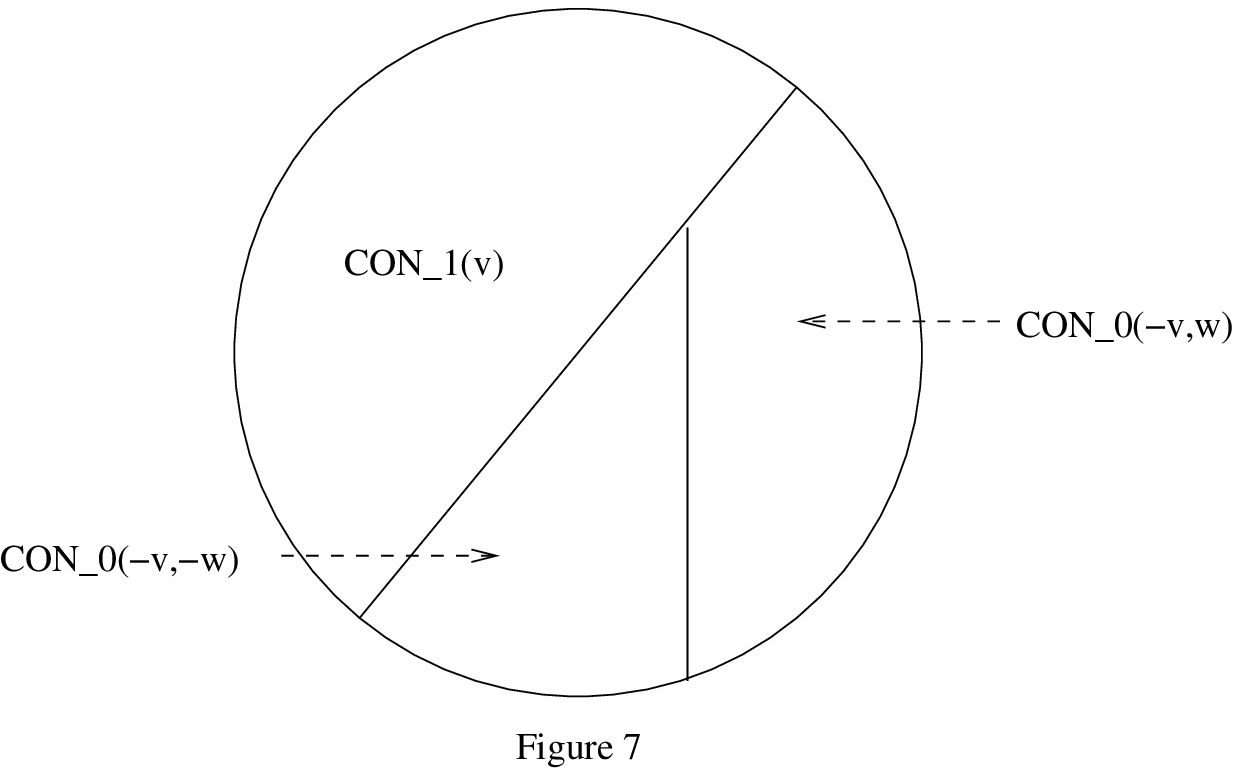}

Now we have to estimate the k-thickness of this pyramid of convex
sets.

We let $\beta$ be a large number.  Our proof will work for each
(sufficiently large) $\beta$, and it will prove the following
estimate for the k-thickness $T$ of our i-pyramid.

$$T \ge c(\beta,n) [2-\epsilon(\beta,n)]^{\frac{n-i-k}{n-i} P}.$$

As $\beta \rightarrow + \infty$, the constant $\epsilon(\beta,n)
\rightarrow 0$, but so does $c(\beta, n)$.  Throughout this section, 
we make the
convention that $\epsilon = \epsilon(\beta, n)$ is a constant
depending only on $\beta$ and $n$ which {\it may change from line
to line}, but always obeying the condition $\lim_{\beta
\rightarrow + \infty} \epsilon(\beta, n) = 0$ for each $n$.

Recall that $T = T_P(*) = \inf_{v \in O(P-1)} T_{P-1}(v) +
T_{P-1} (-v).$ The function $T_{P-1}$ is continuous, so we can
choose $v$ so that $T = T_{P-1}(v) + T_{P-1}(-v)$. Now we proceed
down.  We recall that $T_{P-1}(v) = \inf_{v_1 | (v, v_1) \in
O(P-2)} T_{P-2}(v, v_1) + T_{P-2}(v, - v_1)$.  We choose $v_1$ so
that $T_{P-1}(v) = T_{P-2}(v, v_1) + T_{P-2}(v, - v_1)$. 
Similarly, we choose $v_2$ so that $T_{P-1}(-v) = T_{P-2}(-v,
v_2) + T_{P-2}(-v, - v_2)$.  Continuing in this way, we come to
the following situation.  At the top level, we have a single
convex set $C$ (equal to the unit ball).  At the next level, the
set $C$ is chopped into two sets $C_1$ and $C_2$ along some
plane, perpendicular to the vector $v$.  The sets have equal
volume.  Then each set $C_1$ and $C_2$ is chopped in two.  The
set $C_1$ is chopped in two along a plane perpendicular to $v_1$,
yielding $C_{11}$ and $C_{12}$.  The set $C_2$ is chopped in two
along the plane perpendicular to $v_2$, yielding $C_{21}$ and
$C_{22}$.  This proceeds through P levels, ending up with $2^P$
convex sets.  The chopping obeys two rules.

1. Volume bisection.  Each chop cuts the given set into two sets
of equal volume.  If $C$ has volume 1, then each set at level p has
volume $2^{-p}$.

2. Orthogonality.  Any $(n-i)$ consecutive cuts are orthogonal. 
More precisely, if $I_{p+1}$, ..., $I_{p + n -i}$ are strings of
1's and 2's, with $I_q$ of length q, so that each string is formed
from the previous by adding a 1 or 2 at the end, then the vectors
$v_{I_{p+1}}, ..., v_{I_{p+n-i}}$ are orthogonal.

We let $\{ 1, 2 \}^p$ denote all strings of 1's and 2's of length
$p$.  The k-thickness $T$ of our pyramid of convex sets is equal to
$(1/2) \sum_{I \in \{1, 2 \}^P} Rad[C_I]^k$.  We need to prove
the estimate $T \ge c(\beta, n) (2 -
\epsilon(\beta,n))^{\frac{n-i-k}{n-i} P}$.  In the process, we
will prove the following stronger estimate.

$$ \sum_{I \in \{1, 2 \}^P} Rad[C_I]^{-1} \le C(\beta, n) 2^P (2
+ \epsilon)^{P/(n-i)}. \eqno{(*)}$$

\noindent We check that inequality $(*)$ implies the inequality
that we want to prove.

\begin{lemma} The inequality $(*)$ implies the inequality

$T = (1/2) \sum_{I \in \{1, 2 \}^P} Rad[C_I]^k \ge c(\beta, n) (2
- \epsilon)^{\frac{n-i-k}{n-i} P}.$

\end{lemma}

\proof We begin with the following trivial equation.

$$2^P = \sum_{I \in \{1,2\}^P} Rad[C_I]^{\frac{k}{k+1}}
Rad[C_I]^{- \frac{k}{k+1}}.$$

Applying the Holder inequality to the right-hand side we get the
following.

$$ 2^P \le \bigl[ \sum_{I \in \{1,2\}^P} Rad[C_I]^k \bigr]^{\frac{1}{k+1}}
\bigl[ \sum_{I \in \{1,2\}^P} Rad[C_I]^{-1} \bigr]^{\frac{k}{k+1}}.$$

We use inequality $(*)$ to estimate the last factor.

$$ 2^P \le C(\beta, n)  \bigl[ \sum_{I \in \{1,2\}^P}
Rad[C_I]^k \bigr]^{\frac{1}{k+1}} \bigl[ (2 +
\epsilon)^{P/(n-i)} 2^P \bigr]^{\frac{k}{k+1}}.$$

Rearranging the terms, we get the following.

$$T = (1/2) \sum_{I \in \{1,2\}^P} Rad[C_I]^k \ge c(\beta, n) 2^P (2 +
\epsilon)^{- \frac{k P}{n-i}}.$$

We are assuming that $i \le n-k-1$, and so $n-i-k \ge 1$. 
Therefore, after redefining $\epsilon$, we get the following.

$$T = (1/2) \sum_{I \in \{1, 2\}^P} Rad[C_I]^k \ge c(\beta, n)
(2-\epsilon)^{\frac{n-i-k}{n-i} P}.$$ \endproof

We would like to work inductively by estimating the quantity
$\sum_{I \in \{ 1, 2 \}^p} Rad[C_I]^{-1}$ as $p$ increases from
$0$ to $P$.  Unfortunately, I don't see how to control the
behavior of this quantity inductively.  The main idea of the
proof is to find a related quantity that behaves better from an
inductive point of view and that agrees with $Rad[C_i]^{-1}$ up
to a controlled error.

For each convex set $C$, let $\pi_q(C)$ be the average q-volume
of any orthogonal projection of $C$ onto any q-plane.  (We define
$\pi_0(C) = 1$ for any convex set.)  We now define the quantity
that we will use in our inductive argument.

$$N_{a, \beta}(C) := \beta^a \textrm{ Vol}(C)^{-1/a} \pi_{n-a}(C)^{1/a}.$$

$$N_\beta(C) := \sup_{1 \le a \le n} N_{a, \beta}(C).$$

For large $\beta$, the function $N_\beta$ behaves well in our inductive
arguments.
We will prove below that $N_\beta(C)$ agrees with $Rad[C]^{-1}$ up to a
factor $C(\beta, n)$.  In the proof, we will use the fact that
any convex set can be fairly well approximated by a rectangle, as
described in the following lemma.

\begin{lemma} There is a constant $C(n)$ that makes the following
true.  For any bounded open convex set $A$ in $\mathbb{R}^n$,
there is a rectangle $R$ so that $R \subset A \subset C(n) R$. 
(Here $C(n) R$ denotes the rectangle obtained by magnifying $R$
by a factor $C(n)$ around its center.)
\end{lemma}

\proof Pick a point $a_0 \in A$.

Then pick a point $a_1 \in A$ as far away as possible from $a_0$. 
Let $P_1$ denote the 1-plane containing $a_0$ and $a_1$, and let
$h_1$ denote the distance from $a_0$ to $a_1$.

Then pick a point $a_2 \in A$ as far away as possible from $P_1$. 
Let $P_2$ denote the 2-plane containing $a_0, a_1,$ and $a_2$,
and let $h_2$ denote the distance from $a_2$ to $P_1$.

Continuing inductively we pick points $a_i$ for $0 \le i \le n$,
and define $P_i$ to be the plane spanned by $a_0, ..., a_i$, and
$h_i$ to be the distance from $a_i$ to $P_{i-1}$. 

We now rotate and translate the coordinates of $\mathbb{R}^n$ so
that the point $a_0$ is the origin and so that the plane spanned
by the first $i$ coordinates $x_1, ..., x_i$ is the plane $P_i$. 
If necessary, we make some reflections so that the $x_i$
coordinate of $a_i$ is $h_i > 0$. 

In these coordinates, we can write $a_i = (a_{i,1}, ..., a_{i,
i-1}, h_i, 0, ..., 0)$, where $|a_{i,j}| \le h_j$.  It follows by
induction on the dimension that the convex hull of the $a_i$
contains a rectangle $R$ of the form $[r_1, s_1] \times ...
\times [r_n, s_n]$ with $0 < r_i < s_i < h_i$ and $s_i - r_i \ge
c(n) h_i$.

On the other hand, the set $A$ is contained in the rectangle
$[-h_1 \times h_1] \times ... \times [-h_n , h_n]$.  This
rectangle is in turn contained in $C(n) R$ for an appropriate
dimensional constant $C(n)$. \endproof

Given a convex set $C$ we choose a rectangle $R$ obeying the
conclusion of the lemma, and we let $R_1 \le ... \le R_n$ be its
dimensions.  Then we define $R_i(C)$ to be $R_i$.  Since there
are many rectangles obeying the conclusion of the lemma, $R_i(C)$
is not uniquely defined, but it is defined up to a constant
factor $C(n)$.  The numbers $R_i(C)$ give a rough description of
the shape $C$.

We now describe the sizes $N_{a, \beta}(C)$ in terms of the dimensions
$R_i(C)$. We first note that $\pi_{n-a}(C) \sim R_{a+1} ... R_n$. 
(We write $A \sim B$ to indicate that $A$ approximates $B$ up to
a factor $C(n)$ depending only on the dimension n.)  Using this,
we can describe $N_{a, \beta}(C)$ up to a constant factor $C(n)$.

$$N_{a, \beta} (C) \sim \beta^a (R_1 ... R_a)^{-1/a}.$$

Using this formula, we now check that $N_\beta(C)$ is roughly equal to
$Rad[C]^{-1}$.

\begin{lemma} For any convex set $C$, $N_{\beta}(C)$ agrees with
$Rad[C]^{-1}$ up to a factor $C(\beta, n)$.
\end{lemma}

\proof The invariant $N_{\beta}(C)$ is at least $N_{1, \beta}(C)$, which is at
least $c(n) \beta R_1^{-1}$.  The radius $Rad[C]$ agrees with
$R_1(C)$ up to a constant factor, and so $N_\beta(C) \ge c(n) \beta
Rad[C]^{-1}$.  On the other hand, each $N_{a, \beta}(C)$ is at most $C(n)
\beta^a (R_1 ... R_a)^{-1/a}$, which is at most $C(n) \beta^n
R_1^{-1}$.  Using again the fact that $R_1(C) \sim Rad[C]$, we
conclude that $N_\beta(C) \le C(n) \beta^n Rad[C]^{-1}$. \endproof

We need a few more definitions before we continue.  If $I$ is a
string of 1's and 2's, we denote the length of $I$ by $|I|$. 
Then we define $D_q(I)$ to be the set of all strings of $1$'s and
$2$'s of length $|I| + q$ whose first $|I|$ entries are $I$.  The
expression $D_q(I)$ stands for the descendants of $I$ after $q$
generations.

We have shown that the function $N_\beta(C)$ is roughly $Rad(C)^{-1}$.  The
advantage of $N_\beta(C)$ is that it behaves better in inductive
arguments, especially for large $\beta$.  
The key fact is the following lemma which implies that
$N_{a, \beta}(C_I)$ behaves well inductively whenever $R_{a+1}(C_I) >>
R_a(C_I)$.

\begin{lemma} Suppose that $C_I$ obeys the inequality
$R_{a+1}(C_I) \ge \beta R_a(C_I)$.
Then the following inequality holds.

$$\sum_{J \in D_{n-i}(I)} \pi_{n-a}(C_J) \le (2 + \epsilon)^a
\pi_{n-a}(C_I).$$

\end{lemma}

\proof We first consider the case $a \ge n-i$.  In this case,
the conclusion of the lemma holds regardless of the dimensions of
$R(C_I)$.  Since $C_J \subset C_I$, $\pi_{n-a}(C_J) \le \pi_{n-a}(C_I)$.
The number of $C_J$ is $2^{n-i}$.  Therefore, $\sum_{J \in D_{n-i}(I)}
\pi_{n-a}(C_J) \le 2^{n-i} \pi_{n-a}(C_I) \le 2^a \pi_{n-a}(C_I)$.
We turn to the interesting case $a < n-i$.

We pick a particular rectangle $R$ with $R \subset C_I \subset
C(n) R$ and with side lengths $R_i = R_i(C_I)$, and we let $P$
denote the span of the subrectangle with dimensions $R_1 \times
... \times R_a$.

If $C$ is any convex set, and if $S$ is a piece of hyperplane in
$C$ cutting $C$ into convex sets $A$ and $B$, then we have the
following formula for any $q$.

$$\pi_q(C) + \pi_q(S) = \pi_q(A) + \pi_q(B).$$

This formula is well-known in integral geometry.  It follows
because if $\pi$ is a projection onto a given q-plane, $\pi(\bar
A) \cup \pi(\bar B)$ is exactly $\pi(\bar C)$ and $\pi(\bar A)
\cap \pi(\bar B)$ is exactly $\pi(\bar S)$.  The last formula
follows because if $a \in \pi^{-1}(q) \cap \bar A$ and $b \in
\pi^{-1}(q) \cap \bar B$, then the line segment from $a$ to $b$
lies in $\pi^{-1}(q)$ and intersects $\bar S$.  Since taking
closures does not affect the volume of the image, we see that
$|\pi(A)| + |\pi(B)| = |\pi(C)| + |\pi(S)|$.  Since this formula
holds for each projection $\pi$, it also holds for the average,
proving our formula.

Also, since $S
\subset C$, $\pi_q(S) \le \pi_q(C)$, and so $\pi_q(A) + \pi_q(B)
\le 2 \pi_q(C)$.

Now, we suppose that $C_{L,1}$ and $C_{L,2}$ are immediate
descendants of $C_L$, which is a descendant of $C_I$ of
generation at most $n-i$.  We let $S_L$ denote the hyperplane
in $C_L$ which separates $C_{L,1}$ from $C_{L,2}$.  Since
$S_L$ is a subset of $C_L$, $\pi_q(S_L) \le \pi_q(C_L)$,
which gives us the following inequality.

$$\pi_{n-a} (C_{L,1}) + \pi_{n-a} (C_{L,2}) \le
2 \pi_{n-a}(C_L). \eqno{(1)}$$

If the angle between $v_L$ and $P$ is at least $c(n) > 0$, then
we will prove a stronger estimate.  In this case, $S_L$ is
contained in an (n-1)-dimensional rectangle with dimensions at
most $C(n)$ times bigger than $R_1
\times... \times R_a \times R_{a+2} \times ... \times R_n$. 
Therefore, $\pi_{n-a}(S_L) \le C(n) R_a
R_{a+2} ... R_n \le C(n)
\beta^{-1} R_{a+1} ... R_n$. On the other hand, $\pi_{n-a}(C_L)
\ge c(n) R_{a+1} ... R_n$.  Therefore, we get the following
bound, provided that the angle between $v_L$ and $P$ is at least
$c(n)$.

$$ \pi_{n-a}(C_{L,1}) + \pi_{n-a}(C_{L,2}) \le (1 + C(n)
\beta^{-1}) \pi_{n-a}(C_L). \eqno{(2)}$$

To organize this information, we define a function $h(L)$ that
counts the number of vectors $v_M$ close to $P$ that appear in
the ancestors of $L$.  We define $h(L)$ inductively as follows. 
First, $h(I) = a$.  Now, we suppose we have defined $h(L)$, and
we want to define $h$ on the immediate descendants of $L$.  If the angle
between the vector $v_L$ and the plane $P$ is at most $c(n)$,
then we define $h(L,1) = h(L,2) = h(L) - 1$. On the other hand,
if the angle is greater than $c(n)$, then we define $h(L,1) =
h(L,2) = h(L)$.  With this definition, we can combine equations
(1) and (2) into a single clean inequality.

$$(1 + C(n) \beta^{-1}) 2^{h(L)} \pi_{n-a}(C_L) \ge 2^{h(L,1)}
\pi_{n-a}(C_{L,1}) + 2^{h(L,2)} \pi_{n-a}(C_{L,2}).$$

Applying this inequality repeatedly, we get the following.

$$(1 + C(n) \beta^{-1})^{n-i} 2^{h(I)} \pi_{n-a}(C_I) \ge \sum_{J
\in D_{n-i}(I)} 2^{h(J)} \pi_{n-a}(C_J).$$

Since any sequence of (n-i) vectors is orthonormal, the number of
vectors in any sequence with angle at most $c(n)$ from $P$ is at
most $a$.  Therefore, $h(J) \ge 0$ for every $J \in D_{n-i}(I)$. 
Since $h(I) = a$, we get the following inequality.

$$\sum_{J \in D_{n-i}(I)} \pi_{n-a}(C_J) \le (1 + C(n)
\beta^{-1})^{n-i} 2^a \pi_{n-a}(C_I).$$

This inequality proves our lemma.
\endproof

Using this lemma, we can inductively control the behavior of
$N_{a, \beta}$.

\begin{lemma} Suppose that $C_I$ obeys the inequality
$R_{a+1}(C_I) \ge \beta R_a(C_I)$.  Then the following
inequality holds.

$$\sum_{J \in D_{n-i}(I)} N_{a, \beta}(C_J) \le (2 + \epsilon) 2^{n-i}
N_{a, \beta}(C_I).$$

\end{lemma}

\proof The left-hand side of the inequality we want to prove is

$$\sum_{J \in D_{n-i}(I)} \beta^a \textrm{ Vol}(C_J)^{-1/a}
\pi_{n-a}(C_J)^{1/a}. \eqno{(a)}$$

The volume of $C_J$ is independent of $J$, and it's equal to
$2^{-(n-i)} \textrm{ Vol}(C_I)$.  We apply Holder's inequality to
the term $\sum \pi_{n-a}(C_J)^{1/a}$.

$$\sum_{J \in D_{n-i}(I)} \pi_{n-a}(C_J)^{1/a} \le
2^{(n-i)\frac{a-1}{a}} [\sum_{J \in D_{n-i}(I)} \pi_{n-a}(C_J)]^{1/a}.$$

Applying the last lemma, we get the following inequality.

$$\sum_{J \in D_{n-i}(I)} \pi_{n-a}(C_J)^{1/a} \le 2^{(n-i)
\frac{a-1}{a}} (2 + \epsilon) \pi_{n-a}(C_I)^{1/a}.$$

Putting the last inequality into expression (a), we get the following.

$$\sum_{J \in D_{n-i}(I)} \beta^a \textrm{ Vol}(C_J)^{-1/a}
\pi_{n-a}(C_J)^{1/a} \le \beta^a \textrm{ Vol}(C_I)^{-1/a} 2^{n-i}
(2 + \epsilon) \pi_{n-a}(C_I)^{1/a}.$$

Plugging in the definition of $N_{a, \beta}$ finishes the proof. \endproof

The last key observation is that if $N_{a, \beta}(C) = N_\beta(C)$, then
$R_{a+1}(C) >> R_a(C)$.  This observation is made precise in the
following lemma.

\begin{lemma} If $N_{a, \beta}(C) = N_\beta(C)$, then the dimensions of $C$
obey the following inequality.

$R_{a+1}(C) \ge c(n) \beta^2 R_a(C).$

\end{lemma}

\proof We abbreviate $R_i(C)$ as $R_i$.  Since $N_{a, \beta}(C) = N_\beta(C)$,
we know that $N_{a, \beta}(C) \ge N_{a-1, \beta}(C)$ and $N_{a, \beta}(C) \ge
N_{a+1, \beta}(C)$.  Since $N_{a, \beta}(C) \sim \beta^a (R_1 ... R_a)^{-1/a}$,
we get the following inequalities.

$$ \beta^a R_1^{-1/a} ... R_a^{-1/a} \ge c(n) \beta^{a-1}
R_1^{-1/(a-1)} ... R_{a-1}^{-1/(a-1)}. \eqno{(1)}$$

$$ \beta^a R_1^{-1/a} ... R_a^{-1/a} \ge c(n) \beta^{a+1}
R_1^{-1/(a+1)} ... R_{a+1}^{-1/(a+1)}. \eqno{(2)}$$

In these inequalities and below, we use $c(n)$ to denote a positive
constant, depending only on $n$, whose exact value may change from
line to line.

Taking the second equation and moving $R_{a+1}$ to the left-hand
side and everything else to the right-hand side, we get the
following.

$$ R_{a+1} \ge c(n) \beta^{a+1} R_1^{1/a} ... R_a^{1/a}.
\eqno{(3)}$$

By a similar manipulation, the first equation implies the
following inequality.

$$ R_1^{1/a} ... R_{a-1}^{1/a} \ge c(n) \beta^{-(a-1)}
R_a^{(a-1)/a}. \eqno{(4)}$$

Plugging this inequality into equation 3, we get the following.

$$R_{a+1} \ge c(n) \beta^2 R_a.$$

\endproof

Combining Lemmas 4.7 and 4.8, we can control the
inductive behavior of $N_\beta(C_I)$.  We now pick a number of
generations $G(\beta, n)$.  We choose $G$ just small enough so
that if $J$ is a descendant of $I$ $G$ generations down, then the
dimensions $R_i(C_I)$ and $R_i(C_J)$ agree up to a factor of
$\beta^{1/3}$.  The number of generations $G$ will be on the
order of $\delta(n) \log \beta$ for a small constant $\delta(n)$,
because in each generation the dimensions change by at most a
bounded factor $C(n)$.  We also arrange that $G$ is an integer
multiple of $n-i$.  (We only need to study the cases when $\beta$
is large and so $G$ is large.)

\begin{lemma} The following inequality holds.

$$\sum_{J \in D_G(I)} N_\beta(C_J) \le n 2^G (2 + \epsilon)^{G/(n-i)}
N_\beta(C_I).$$

\end{lemma}

\proof 

Suppose that $N_{a, \beta}(C_J) = N_\beta(C_J)$ for some $J$ in $D_G(I)$.  By
Lemma 4.8, we know that $R_{a+1}(C_J) \ge c(n) \beta^2 R_a(C_J)$.
By the definition of $G$, we know that $R_{a}(C_I) \le
\beta^{1/3} R_a(C_J)$.  On the other hand, $R_{a+1}(C_I) \ge c(n)
R_{a+1}(C_J)$ because $C_J \subset C_I$.  Therefore,
$R_{a+1}(C_I) \ge c(n) \beta^{5/3} R_a(C_I)$.  Similarly, if $K$
is any descendant of $I$ at most $G$ generations down, then
$R_{a+1}(C_K) \ge c(n) \beta^{4/3} R_a(C_K)$.

We are now in a position to apply Lemma 4.7 repeatedly $G/(n-i)$
times.  Assuming that $N_{a, \beta}(C_J) =N_\beta(C_J)$ for some $J$ in $D_G(I)$,
we conclude the following.

$$\sum_{J \in D_G(I)} N_{a, \beta}(C_J) \le 2^G (2+ \epsilon)^{G/(n-i)}
N_{a, \beta}(C_I).$$

Let $A$ denote the set of all $a$ so that for some $J$, $N_{a, \beta}(C_J) =
N_\beta(C_J)$.

$$\sum_{J \in D_G(I)} N_\beta(C_J) \le \sum_{a \in A} \sum_{J \in D_G(I)} N_{a, \beta}(C_J) \le$$

$$\le \sum_{a \in A} 2^G (2 + \epsilon)^{G/(n-i)} N_{a, \beta}(C_I) \le n 2^G (2 + \epsilon)^{G/(n-i)} N_\beta(C_I).$$

\endproof

Finally, we apply the last lemma repeatedly.  If $P$ is a multiple of $G$,
we get the following inequality.

$$\sum_{I \in \{1, 2 \}^P} N_\beta(C_I) \le n^{P/G} (2 +
\epsilon)^{P/(n-i)} 2^P N_\beta(C).$$

In general, if $P$ is not a multiple of $G$, then we can write $P$
as a multiple of $G$ plus a remainder which is at most $G$.  In the last
$G$ generations, the values of $N_\beta(C_I)$ change by at most a constant
$C(\beta, n)$, since $G$ depends only on $\beta$ and $n$.
So for all $P$ we get the following inequality.

$$\sum_{I \in \{1, 2\}^P} N_\beta(C_I) \le C(\beta, n) n^{P/G} (2 + \epsilon)^{P/(n-i)} 2^P N_\beta(C).$$

In this equation, $C$ denotes the unit ball, which is the top of
our pyramid of convex sets.  The term $N_\beta(C)$ is a constant
depending on $\beta$ and $n$.  More importantly, since $G \sim
\delta(n) \log \beta$, $n^{P/G}$ grows at an arbitrarily small
exponential rate in $P$.  Therefore, we can absorb the term
$n^{P/G}$ into the term $(2 + \epsilon)^{P/(n-i)}$, by
changing the definition of $\epsilon$.  We get the following
inequality.

$$\sum_{I \in \{1, 2 \}^P} N_\beta(C_I) \le C(\beta, n) (2 +
\epsilon)^{P/(n-i)} 2^P. $$

According to Lemma 4.8, $N_\beta(C_I)$ agrees with $Rad[C_I]^{-1}$ up
to a factor $C(\beta, n)$.  Making this substitution, we get
inequality $(*)$.

$$\sum_{I \in \{1, 2 \}^P} Rad(C_I)^{-1} \le C(\beta, n) (2 +
\epsilon)^{P/(n-i)} 2^P. \eqno{(*)}$$

According to Lemma 4.3, inequality $(*)$ implies our lower bound
on the k-thickness: $T \ge c(\beta, n) (2 -
\epsilon)^{\frac{n-i-k}{n-i} P}$. \endproof

\section{Families of bent overlapping planes}

In this section, we construct examples of families of cycles with
small volumes, proving the upper bounds in Theorem 1 and Theorem
2. Here is an outline of the approach.

We need to construct a family of k-cycles that detects a given cohomology
class $\alpha \in H^*(Z(k,n))$.  The k-fold suspension map 
$\Sigma: Z(0,n - k) \rightarrow Z(k,
n)$ maps a 0-cycle $C$ in the unit ball $B^{n-k}$ to the product
$C \times \mathbb{R}^k$ restricted to the unit ball $B^n$.

The first step is to construct a family of 0-cycles 
that detects $\Sigma^*(\alpha)$.  Applying $\Sigma$, we
get a family of k-cycles that detects $\alpha$.  Each cycle in
this family is a union of parallel planes.  It turns out that
the cycles in this family have volume much greater than
$\mathbb{V}(\alpha)$.

The second step is the bend-and-cancel construction described in the introduction.
We carefully pick a degree 1 PL map $\Psi$, and we apply the map $\Psi$ to get
a new family of cycles where each cycle consists of a union of bent planes.  By
choosing $\Psi$ carefully, we arrange that the bent planes overlap a great deal.
Since we are working with families of mod 2 cycles, we can cancel the
overlapping parts, reducing the volume.

Before proving the upper bounds in Theorem 1, we gather a few tools
that we will use.  The first tool is the k-fold suspension map 
$\Sigma: Z(0, n-k) \rightarrow Z(k,n)$.  The map is defined in the following
way.  We begin with a 0-cycle $C \in Z(0, n-k)$.  Then we take the product
$C \times \mathbb{R}^k$ which is a (locally finite) cycle in $B^{n-k} \times
\mathbb{R}^k$.  Finally, we restrict this cycle to the unit ball $B^n$.

\begin{lemma} The pullback $\Sigma^* a(k,n)$ is equal to $a(0,n-k)$.
\end{lemma}

\proof Suppose that $F: P \rightarrow Z(0,n-k)$ is a family of 0-cycles, and
suppose that $h$ is a homology class in $H_{n-k}(P)$.  We
need to check that the pairing $<\Sigma F_*(h), a(k,n)>$ is equal to the pairing
$<F_*(h), a(0,n-k)>$.  To compute
the second pairing, we find a complex of cycles $C$ based closely on the family
$F$.  We let $z \subset P$ be a simplicial (n-k)-cycle in the homology class $h$.  Then
the pairing $<F_*(h), a(0,n-k)>$ is equal to the degree of the n-cycle $C(h)$.
Now we define $\Sigma C$, a complex of cycles in $B^n$ parametrized by $P$.  For each
simplex $\Delta$ of $P$, we let $\Sigma C(\Delta)$ be the product $C(\Delta) \times \mathbb{R}^k$
restricted to the unit n-ball.  The point of the proof is that
we can use $\Sigma C$ as the complex of cycles approximating $\Sigma F$.  Therefore,
the pairing $<\Sigma F_*(h), a(k,n)>$ is equal to the degree of $\Sigma C(z)$.  Now
$\Sigma C(z)$ is the product $C(z) \times \mathbb{R}^k$ restricted to the unit n-ball,
so $\Sigma C(z)$ and $C(z)$ have the same degree. \endproof

The first step also 
requires a result of Nakaoka about the cohomology ring of a symmetric
product of spheres.  Nakaoka completely described the cohomology ring, but we
will need only the following facts.

\begin{reftheorem} (Nakaoka, \cite{N}) Let $SP^d S^N$ denote the d-fold symmetric product
of the N-sphere, for $N \ge 2$.  The cohomology group 
$H^N(SP^d S^N)$ is equal to $\mathbb{Z}_2$.
Let $\beta$ denote the generator of this group.  Let $i$ denote the inclusion
of $S^N$ into $SP^d S^N$.  (If $*$ denotes the base point of $S^N$, then $i(x)$ is
defined to be the unordered d-tuple $<x, *, ..., *>$.)  Then $i^* \beta$ is the
generator of $H^N(S^N)$.  The top-dimensional cohomology group $H^{dN} (SP^d S^N)$ is
also isomorphic to $\mathbb{Z}_2$, and it is generated by $\beta^d$.
\end{reftheorem}

The second step is based on the ``bending planes around a
skeleton'' construction in
\cite{Gu1}.  In that paper, the following mappings were constructed.

\begin{reflemma} (\cite{Gu1}) Let $l$ be a dimension in the range
$0 \le l <n$, and let $s > 0$ be a scale.  Let $S$ denote the 
l-skeleton of the lattice with sidelength $s$.  Let $T$ denote
the dual (n-l-1)-skeleton.  

For each dimension $l$, each scale $s$ and each $\epsilon > 0$, there is a
piecewise-linear map $\Psi$ from $\mathbb{R}^n$ to itself with
the following properties.  The map $\Psi$ is linear on each
simplex of a certain triangulation of $\mathbb{R}^n$.  Each
top-dimensional simplex of this triangulation is labelled good or
bad.  For each good simplex $\Delta$, $\Psi(\Delta)$ lies in $S$. 
Each bad simplex lies in the $\epsilon$-neighborhood of $T$.  
The triangulation and the map obey the following bounds.

1. The number of simplices of our triangulation meeting any ball
of radius $s$ is bounded by $C(n)$.

2. The displacement $|\Psi(x) - x|$ is bounded by $C(n) s$.

3. The diameter of each simplex is bounded by $C(n) s$.
\end{reflemma}

We call the map $\Psi$ a skeleton-squeezing map, since it squeezes most of $\mathbb{R}^n$
into the l-skeleton $S$.

\begin{Theorem 1} (Upper bounds) The minimax volume $\mathbb{V}(a(k,n)^p) \le C(n)
p^{\frac{n-k}{n}}.$
\end{Theorem 1}

\proof The first step is to construct a family of cycles $F(p)$
in $Z(0, n-k)$
that detects the cohomology class $a(0, n-k)^p$.  Roughly
speaking, the family of all p-tuples of points in $B^{n-k}$ does
the job.  Our argument involves two cases depending on whether
$n-k = 1$.

If $n-k = 1$, then we define $F(p)$ using the roots of
polynomials.  Let $V(p)$ be the space of all real polynomials of
one variable with degree at most $p$.  The space $V(p)$ is a
vector space of dimension $p+1$.  To each non-zero polynomial in
$V(p)$, we associate its real roots, taken with multiplicity. 
This association defines a map $R_0$ from $V(d) - \{ 0 \}$ to the
space of integral 0-cycles on the real line, but the map is NOT
continuous.  The reason for the discontinuity is that two real
roots may approach each other, become a double root, and then
become two conjugate complex roots.  Since $R_0$ only records the
real roots, two real roots can come together and disappear.
We correct this problem by considering the roots with
multiplicity modulo 2.  We define a root map $R$ from $V(d) - \{
0 \}$ to $Z(0,1)$ by taking the real roots of a polynomial,
keeping only the roots in the interval $(-1, 1)$, and recording
the multiplicity modulo 2.  The map $R$ is continuous.

For any non-zero real number $\lambda$, the polynomials $P$ and
$\lambda P$ have the same roots, and so $R$ induces a map $F(p)$
from $\mathbb{RP}^p = [V(p) - \{ 0 \} ] / \mathbb{R}^*$ to
$Z(0,1)$.  We call this the family of roots of degree d
polynomials.

For example, if $p = 1$, then the map $F(1)$ sends the
polynomial $a x + b$ to its root $-b/a$.  If we fix
$a=1$, then as $b$ goes from $- \infty$ to $+ \infty$,
the point $-b/a$ goes from $+ \infty$ to $- \infty$.  So
the family $F(1)$ sweeps out the unit ball $(-1, 1)$ with
degree 1 modulo 2.  Hence $F(1)^* (a(0,1))$ is the generator
of $H^1(\mathbb{RP}^1)$.

Next we compute that $F(p)^*(a(0,1))$ is the generator of
$H^1(\mathbb{RP}^p)$.  To check this, we pick a homologically
non-trivial curve $c$ in $\mathbb{RP}^p$ and we check that 
$F(c)$ sweeps out the unit interval.  We can take
the curve $c$ given by the projectivization of the linear
polynomials, $V(1) \subset V(p)$.  The map $F(p)$ restricted to this
copy of $\mathbb{RP}^1$ is just $F(1)$, and so the claim
follows from the last paragraph.  Therefore, the family $F(p)$ 
detects $a(0,1)^p$.

If $n-k \ge 2$, then we define $F(p)$ using symmetric products of
spheres. First we define a family $F(1)$ parametrized by
$S^{n-k}$.  We pick a homeomorphism of the upper hemisphere with
the unit ball.  Then we define $F(1)(x)$ for $x$ in the upper
hemisphere to be the corresponding point of the unit ball with
multiplicity 1.  We define $F(1)$ on the lower hemisphere to be
the empty cycle.  (We assume that the basepoint is in the lower
hemisphere, so it maps to the empty cycle.)  The family $F(1)$
sweeps out the unit ball and so it detects $a(0,n-k)$.  Next we
define a family $F(p)$ parametrized by the symmetric product
$SP^p S^{n-k}$.  We define $F(p)$ of an unordered p-tuple $< x_1,
..., x_p >$ to be the sum $\sum_i F(1)(x_i)$.  Using Nakaoka's
theorem, we can check that $F(p)$ detects $a(0,n-k)^p$.  First we
compute $F(p)^*(a(0,n-k))$. The group $H^{n-k}(SP^p S^{n-k})$ is
equal to $\mathbb{Z}_2$ and is generated by the class $\beta$,
and the inclusion $i: S^{n-k} \rightarrow SP^p S^{n-k}$ induces
an isomorphism in $H^{n-k}$.  Therefore, it suffices to compute
$i^* F(p)^* (a(0,n-k))$.  Unwinding the definitions, the map
$F(p) \circ i$ is just $F(1)$, which detects $a(0,n-k)$. 
Therefore, $F(p)^*(a(0,n-k))$ is equal to $\beta$, and
$F(p)^*(a(0,n-k)^p)$ is equal to $\beta^p$.  By Nakaoka's
theorem, $\beta^p$ is non-zero, and so $F(p)$ detects
$a(0,n-k)^p$.

In either case, we have constructed a family $F(p)$ in $Z(0,n-k)$
detecting $a(0,n-k)^p$ where each cycle in the family consists of
at most $p$ points.

(Remark: The case $n-k=1$ is separate because the symmetric
product $SP^p S^1$ is not a cycle.  We could have used the
truncated symmetric product $TP^p S^1$ as in \cite{M}.)

Next we consider the suspension map $\Sigma: Z(0,n-k) \rightarrow
Z(k,n)$.  Let $P(p)$ denote the parameter space of the map $F(p)$
above.  (If $n-k>1$, then $P(p) = SP^p S^{n-k}$. If $n-k=1$, then
$P(p) = \mathbb{RP}^p$.)  We define $F_k(p)$ to be $\Sigma \circ
F(p): P(p)
\rightarrow Z(k,n)$.  By Lemma 5.1, $\Sigma^*(a(k,n)^p) =
a(0,n-k)^p$.  Therefore,
$F_k(p)$ detects $a(k,n)^p$.  Each cycle in $F_k(p)$ consists of a union
of at most $p$ k-planes parallel to the $(x_{n-k+1}, .., x_n)$-plane.

Now we turn to the second step, which is to bend the planes so
that they overlap and cancel the overlaps.  First we scale the
family $F_k(p)$ to get a family of cycles in the ball $B(R)$ for
a radius $R$ that we will choose later.  Then we rotate the
family slightly, so that all the k-planes are parallel to a plane
$P_0$ at a generic angle with respect to the coordinate axes.  We
call the scaled rotated family $F$.  If we restrict $F$ to the
unit ball, we get a family in $Z(k,n)$ that detects $a(k,n)^p$.

Next we bend the cycles using a skeleton-squeezing map $\Psi$ as
described above. The map $\Psi$ depends on three parameters: a
dimension $l$, a scale $s$, and a small number $\epsilon > 0$. 
We choose the dimension $l$ of the skeleton to be $k$, and we
choose the scale $s$ to be $p^{-1/n}$.  (Later we will choose
$\epsilon > 0$ sufficiently small.) The family we are trying to
construct will be $\Psi \circ F$ restricted to the unit ball.

The displacement of $\Psi$ is at most $C(n) s \le C(n)$.  We
choose $R$ sufficiently large that the displacement is less than
$R-1$.  Therefore, $\Psi$ maps the sphere $S(R)$ to the exterior
of the unit ball.  Hence, if $z$ is a relative cycle in $B(R)$,
then $\Psi(z)$ can be restricted to a relative cycle in the unit
ball. So $\Psi \circ F$ defines a family of cycles in the unit
ball.

Now define $\Psi_t$ to be the family of maps $\Psi_t(x) = (1-t) x
+ t \Psi(x)$.  We have $\Psi_0$ equal to the identity and
$\Psi_1$ equal to $\Psi$.  Each map $\Psi_t$ has displacement at
most equal to that of $\Psi$, and so each one maps the sphere
$S(R)$ to the exterior of the unit ball.  Hence we get a homotopy
of families of cycles.  At time 0, we have the restriction of the
family $F$ to the unit ball, which we know detects $a(k,n)^p$. At
time 1, we have the restriction of $\Psi \circ F$ to the unit
ball.  We conclude that the latter family also detects
$a(k,n)^p$.

We now come to the main point: estimating the volumes of the
cycles in $\Psi \circ F$. Since the parallel k-planes in $F$ were
in general position with respect to the coordinate axes, each
k-plane meets the dual (n-k-1)-skeleton $T$ in at most $C(n)$
points.  By taking $\epsilon > 0$ sufficiently small, we can
guarantee that each k-plane meets at most $C(n)$ bad simplices in
$B(R)$. Therefore, for each k-plane $P$, the image $\Psi(P)$ lies
in the k-skeleton $S$ except for at most $C(n)$ pieces of plane,
each with volume at most $C s^k$.  Therefore, each cycle in the
family of bent planes lies inside the k-skeleton $S$ except for
an exceptional region of volume at most $C(n) p s^k = C(n)
p^{\frac{n-k}{n}}$.  The volume of the k-skeleton $S$ intersected
with the unit ball is bounded by $C(n) s^{-n} s^k = C(n)
p^{\frac{n-k}{n}}$.  We cancel all overlaps of cycles inside $S$,
so that each region of $S$ has multiplicity zero or 1.  The
resulting cycle has total volume at most $C(n)
p^{\frac{n-k}{n}}$. \endproof

To prove the upper bounds in Theorem 2 we need a few more tools. 
We will use two more maps between spaces of cycles. The first map
is the addition map which sends a d-tuple of cycles $(z_1, ...,
z_d)$ to the sum $z_1 + ... + z_d$.  The addition in the space of
cycles is commutative, and so we view the addition operation as a
map $A$ from $SP^d Z(0,m)$ to $Z(0,m)$.

The second map is the translation map $T$.  For $t \in [-1, 1]$
and for any cycle $z \in Z(0,m-1)$, we define $T(t,z)$ to be the
restriction of $\{ t \} \times z$ to unit ball $B^m$.  The cycle
$T(t,z)$ is a translation of the cycle $z$.  We illustrate the
map $T$ in Figure 8 below.

\includegraphics{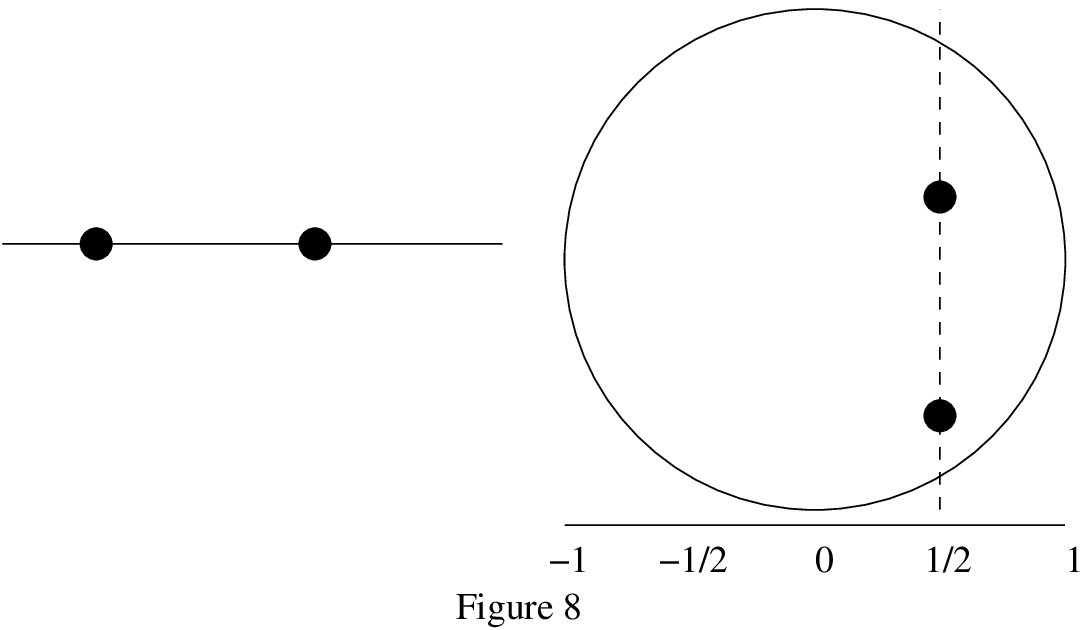}

\noindent The left half of the figure shows a 0-cycle $C \in
Z(0,1)$, consisting of two points.  The right half shows $T(1/2,
z) \in Z(0,2)$, a 0-cycle in the the 2-disk consisting of two
points.

The map $T$ sends $\{-1\} \times Z(0,m-1)$ and $\{1 \} \times
Z(0,m-1)$ to the empty cycle in $Z(0,m)$.  Therefore, it defines
a continuous map $T: S Z(0, m-1) \rightarrow Z(0,m)$, where $S
Z(0,m-1)$ denotes the (non-reduced) suspension of $Z(0,m-1)$.

There is a canonical isomorphism $\sigma: H^{*+1}(S Z(0,m-1))
\rightarrow H^*(Z(0,m-1))$.

\begin{lemma} The cohomology class $\sigma T^* a(0,m) =
a(0,m-1)$.
\end{lemma}

\proof Let $z$ denote a mod 2 (m-1)-cycle in $Z(0,m-1)$.  We need
to check that the pairings $<\sigma S^*(a(0,m)), z>$ and
$<a(0,m-1), z>$ agree.  To compute the second pairing, we pick a
fine triangulation $Tri$ of $z$, and we let $C$ be a complex of
cycles based on this triangulation. Then we compute the degree of
the top-dimensional cycle $C(z)$.  The first pairing is the same
as $<T^*(a(0,m)), Sz>$, where $Sz$ denotes the suspension of $z$
inside of $S Z(0,m-1)$.  This pairing is also equal to $<a(0,m),
T(Sz)>$, where we view $T(Sz)$ as a family of cycles in $Z(0,m)$
parametrized by $Sz$, the suspension of $z$.  We divide the
segment $[-1, 1]$ into segments of length $1/N$.  Then we define
a polyhedral decomposition of $Sz$ with faces of the form $\Delta
\times \{m/N \}$ or $\Delta \times [m/N, (m+1)/N]$, where
$\Delta$ is a simplex of $Tri$, our triangulation of $z$.  Now
$T(Sz)$ is a family of cycles in $Z(0,m)$ parametrized by $Sz$. 
We can build a complex of cycles $\tilde C$ approximating this
family by taking $\tilde C (\Delta \times \{ m/N \} = C(\Delta)
\times \{m/N\}$ and $\tilde C(\Delta \times [m/N, (m+1)/N] =
C(\Delta) \times [m/N, (m+1)/N]$.  It then follows that $\tilde
C(Sz)$ is equal to the restriction to the unit ball of $C(z)
\times [-1,1]$.  Therefore, $C(z)$ and $\tilde C(Sz)$ have the
same degree, and so $<a(0,m-1), z>$ and $<a(0,m), T(Sz)>$ are
equal. \endproof

The last tool concerns symmetric products of pseudomanifolds.  
We use the following vocabulary.
Let $X$ be a finite simplicial complex
of dimension $N$.  Let $S_K(X)$ denote the K-skeleton of $X$.  We
say that $X$ is an N-dimensional pseudomanifold if each (N-1)-simplex of $X$ is 
contained in exactly 2 N-simplices and if $X - S_{N-2}(X)$ is connected.
An equivalent definition is that $X - S_{N-2}(X)$ is a connected N-manifold.
An N-dimensional pseudomanifold $X$ has $H_N(X, \mathbb{Z}_2) = \mathbb{Z}_2$.  We call
the non-trivial element the fundamental homology class of $X$ and denote
it by $[X]$.  If $X$ and $Y$ are N-dimensional pseudomanifolds and $F: X \rightarrow Y$
is a continuous map, then we define the degree of $F$ by $F_*([X]) = (deg F) [Y]$.  
(The degree lies in $\mathbb{Z}_2$.)  If $y$ denotes a generic point in $Y$, then
the degree of $F$ is equal to the number of preimages $F^{-1}(y)$ taken
mod 2, just as for manifolds.  We will need the following lemma about
symmetric products of pseudomanifolds.

\begin{lemma} Suppose that $X$ is an N-dimensional pseudomanifold for $N \ge 2$.  Then
the symmetric product $SP^d X$ is a (dN)-dimensional pseudomanifold.  If $X$ and $Y$
are N-dimensional pseudomanifolds with $N \ge 2$, $F: X \rightarrow Y$ is a map, and 
$SP^d F: SP^d X \rightarrow SP^d Y$ is the d-fold symmetric product of $F$,
then the degree of $SP^d F$ is equal to the degree of $F$.
\end{lemma}

\proof Suppose $X$ is a pseudomanifold.  First we check that
$X^d$ is a pseudomanifold.  We have to consider $X^d - S_{dN - 2} X^d$.  This
set is an open subset of the d-fold product of $X - S_{N-2}(X)$, so it
is a manifold.  Moreover, the complement of this set in $(X - S_{N-2}(X))^d$
has codimension 2, and so this set is connected.  Hence $X^d$ is a pseudomanifold.

Now assume that $N \ge 2$ and consider the symmetric product $SP^d X$.  The
diagonal is the subset of $X^d$ where at least two entries are equal.  We can
triangulate $X^d$ so that the diagonal is a
subcomplex of dimension $N(d-1) \le dN - 2$.  Let $\pi$ denote the quotient
map from $X^d$ to $SP^d X$.  We choose a triangulation of $SP^d X$ so that 
$\pi$ maps the K-skeleton of $X^d$ into the K-skeleton of $SP^d X$ for each $K$.
In particular, the image of the diagonal lies in $S_{dN-2} SP^d X$.  
Let $A$ denote $SP^d(X) - S_{dN-2} SP^d X$.  If we restrict the quotient map 
$\pi$ to $\pi^{-1}(A)$ we get a covering map.  Therefore, $\pi^{-1}(A)$ is
an open subset of $X^d - S_{dN-2} X^d$, and hence a manifold.  Moreover,
the complement of $\pi^{-1}(A)$ in $X^d - S_{dN-2} X^d$ has codimension 2, so
$\pi^{-1}(A)$ is a connected manifold.
Therefore $A$ is the quotient of a connected manifold by a finite group acting
freely and properly, and so $A$ is a connected manifold.  Therefore
$SP^d X$ is a pseudomanifold.

Now assume that $X, Y$ are N-dimensional pseudomanifolds, $N \ge 2$, and that $F: X
\rightarrow Y$ is a continuous map.  We let $F^d$ denote the product map
from $X^d$ to $Y^d$.  By the Kunneth theorem, the degree of $F^d$ is equal
to the degree of $F$ raised to the power $d$.  Since we are working mod 2, the degree
of $F$ and the degree of $F^d$
are equal.  Finally we check that the degree of $SP^d F$ is the same as
the degree of $F^d$.  Let $y$ denote a generic point in $SP^d Y$.  To compute
the degree of $SP^d F$, we count the number of preimages in $[SP^d F]^{-1}(y)$.
The point $y$ has $d!$ preimages in $Y^d$.  Let $y_1$ be one of these preimages.
To compute the degree of $F^d$, we count the number of preimages in
$[F^d]^{-1}(y_1)$.  Because $\pi$ is a covering over generic points, there is
one preimage of $y_1$ in $X^d$ lying over each preimage of $y$ in $SP^d X$.
Therefore, the degree of $F^d$ is equal to the degree of $SP^d F$. \endproof

\begin{Theorem 2} (Upper bounds) The minimax volume $\mathbb{V}(Sq_0^{Q_0} ...
Sq_{n-k-1}^{Q_{n-k-1}} a(k,n)) \le C(n) \prod_{i=0}^{n-k-1}
2^{\frac{n-k-i}{n-i} Q_i}$.
\end{Theorem 2}

\proof We begin by considering the space of 0-cycles $Z(0,m)$.  For each
cohomology class of the form $Sq_0^{Q_0} ... Sq_{m-1}^{Q_{m-1}} a(0,m)$, we will construct
a family of 0-cycles $F_m(Q_0, ..., Q_{m-1})$ that detects it.  These families of
0-cycles are based on certain families of subsets of the unit sphere that were
explained to me by David Wilson.

We will abbreviate $Sq_0^{Q_0} ... Sq_{m-1}^{Q_{m-1}} a(0,m)$ by
$Sq^Q a(0,m)$ and $F_m(Q_0, ..., Q_{m-1})$ by $F_m(Q)$.  The parameter space of
$F_m(Q_0, ..., Q_{m-1})$ will be called $P_m(Q_0, ..., Q_{m-1})$.  Each parameter space will
be a pseudomanifold, and we will check that the pairing $<F_m(Q)^*(Sq^Q a(0,m)), [P_m(Q)] > = 1$.

The construction is inductive in the dimension $m$.  The base
case is $m=1$, and here we have to define $F_1(Q_0)$.  The family
$F_1(Q_0)$ needs to detect $Sq_0^{Q_0} a(0,1) =
a(0,1)^{2^{Q_0}}$.  We constructed such a family at the beginning
of the proof of Theorem 1 (upper bounds), by looking at the set
of roots of polynomials of degree $2^{Q_0}$.  This family is
parametrized by real projective space, so we have $P_1(Q_0) =
\mathbb{RP}^{2^{Q_0}}$.  We checked in the proof of Theorem 1
that the pairing $<F_1(Q_0)^* (Sq_0^{Q_0} a(0,1)), [P_1(Q_0)] > =
1$.

We can now inductively define $F_m(Q_0, ..., Q_{m-1})$, using the
addition map $A$ and the translation map $T$, defined above.  By
the inductive hypothesis, we suppose we already have a family
$F_{m-1}(Q_1, ..., Q_{m-1}): P_{m-1}(Q_1, ..., Q_{m-1})
\rightarrow Z(0,m-1)$.  We first define $F_m(0, Q_1, ...,
Q_{m-1})$.  In this case, we define the parameter space $P_m(0,
Q_1, ..., Q_{m-1})$ to be the suspension $S P_{m-1}(Q_1, ...,
Q_{m-1})$.  We define the map $F_m(0, Q_1, ..., Q_{m-1})$
according to the following diagram.

\[
\begin{CD}
S P_{m-1}(Q_1, ..., Q_{m-1}) @>S F_{m-1}(Q_1, ..., Q_{m-1})>> S
Z(0, m-1) @>T>> Z(0,m)
\end{CD}
\]

We still have to define $F_m(Q_0, ..., Q_{m-1})$ for $Q_0 > 0$. 
To save space, we introduce the following notation.  We use $X$
to denote the space $P_m(0, Q_1, ..., Q_{m-1})$, and we use $f$
to denote the map $F_m(0, Q_1, ..., Q_{m-1}): X \rightarrow
Z(0,m)$.  We define $P_m(Q_0, ..., Q_{m-1})$ to be the symmetric
product $SP^{2^{Q_0}} X$.  To save space, we use $d$ to denote
$2^{Q_0}$.  Now we define $F_m(Q_0, ..., Q_{m-1})$ as the
following composition.

\[
\begin{CD}
SP^d X @>SP^d f>> SP^d Z(0,m) @>A>> Z(0,m)
\end{CD}
\]

To begin, we verify that $P_m(Q)$ is a pseudomanifold.
We already checked that $P_1(Q_0)$ is a pseudomanifold.
We assume that $P_{m-1}(Q_1, ..., Q_{m-1})$ is a pseudomanifold.  The suspension of a pseudomanifold 
is a pseudomanifold, so $S P_{m-1}(Q_1,
..., Q_{m-1})$ is a pseudomanifold.  The pseudomanifold $P_{m-1}(Q_1, ..., Q_{m-1})$ has dimension at least
1, and so its suspension has dimension at least 2.  Therefore the symmetric product
$SP^{2^{Q_0}} S P_{m-1} (Q_1, ..., Q_{m-1})$ is a pseudomanifold.

Next we have to check that the pairing $<F_m(Q)^* (Sq^Q a(0,m)), [P_m(Q)]> = 1$.  Again
we proceed inductively.  We have already checked this equation for $m=1$.  By induction,
we assume that it holds for $m-1$.

First we consider the special case that $Q_0 = 0$.  By definition we need to compute
the following pairing.

\vskip3pt

$<F_{m}(0, Q_1, ... Q_{m-1})^* (Sq_1^{Q_1} ... Sq_{m-1}^{Q_{m-1}} a(0,m)), [P_m(0, Q_1, ..., Q_{m-1})] >$

\vskip3pt

Plugging in the definitions of $F_m(0, Q_1, ..., Q_{m-1})$ and $P_m(0, Q_1, ..., Q_{m-1})$, we get
the following expression.

\vskip3pt

$=<SF_{m-1}(Q_1, ..., Q_{m-1})^* T^*(Sq_1 ^{Q_1} ... Sq_{m-1}^{Q_{m-1}} a(0,m)), [S P_{m-1}(Q_1, ..., Q_{m-1})] >$

\vskip3pt

Using the suspension isomorphism and the fact that $T^*$ commutes with Steenrod squares, 
we get the following expression.

\vskip3pt

$=<F_{m-1}(Q_1, ..., Q_{m-1})^* \sigma (Sq_1^{Q_1} ... Sq_{m-1}^{Q_{m-1}} T^* a(0,m)), [P_{m-1}(Q_1, ..., Q_{m-1})]>.$

\vskip3pt

Steenrod squares commute with $\sigma$.  In other words, $\sigma Sq^i \alpha = Sq^i \sigma \alpha$, as described
in \cite{H}.  We are using lower squares.  Rewriting lower squares in terms of upper squares, it follows
that $\sigma Sq_i \alpha = Sq_{i-1} \sigma \alpha$.  Using this equation, we can interchange $\sigma$ with
the Steenrod squares to get the following formula.

\vskip3pt

$=<F_{m-1}(Q_1, ..., Q_{m-1})^* (Sq_0^{Q_1} ...
Sq_{m-2}^{Q_{m-1}} \sigma T^* a(0,m)), [P_{m-1}(Q_1, ...,
Q_{m-1})]>.$

\vskip3pt

According to Lemma 5.2, $\sigma T^* a(0,m) = a(0,m-1)$. 
Substituting $a(0,m-1)$ for $\sigma T^* a(0,m)$ in the expression
above leaves the following.

$$= <F_{m-1}(Q)^* (Sq_0^{Q_1} ... Sq_{m-2}^{Q_{m-1}} a(0, m-1)),
[P_{m-1}(Q_1, ..., Q_{m-1})] >.$$

By our inductive hypothesis, this pairing is equal to $1$.

To finish the inductive step, we have to deal with the case that
$Q_0 > 0$.  Let $\alpha$ denote $Sq_1^{Q_1} ...
Sq_{m-1}^{Q_{m-1}} a(0,m)$.  Using the abbreviations above, we
have just checked that the pairing $<f^* \alpha, [X]>$ is equal
to 1.  To finish our induction, we need to check that $<(A \circ
SP^d f)^* \alpha^d, [SP^d X]> = 1$.

Let $N$ be the dimension of the cohomology class $\alpha$, and so
also the dimension of $X$.  Let $g: X \rightarrow S^N$ be a map
of degree 1 (mod 2).

Recall that $X$ is a suspension.  We let $x_0$ be the vertex of
one of the two cones whose union is $X$, and we consider $x_0$ to
be a basepoint for $X$. Then we choose a basepoint $p_0$ in
$S^N$, arranging that $g(x_0) = p_0$.  Using the basepoints, we
define the embedding $i: X \rightarrow SP^d X$ by taking a point
$x$ to the d-tuple $<x, x_0, ..., x_0>$ and the embedding $i: S^N
\rightarrow SP^d S^N$ by taking a point $p$ to the d-tuple $<p,
p_0, ..., p_0>$.  Since $f = T \circ SF_{m-1}(Q_1, ...,
Q_{m-1})$, $f$ maps the vertex $x_0$ to the empty cycle. 
Therefore, the following diagram commutes.

\[
\begin{CD}
Z(0,m) @<f<< X @>g>> S^N \\
@VV=V @VViV @VViV \\
Z(0,m) @<A \circ SP^d f<< SP^d(X) @>SP^d g>> SP^d S^N
\end{CD}
\]

\vskip4pt

According to Nakaoka (\cite{N}), $H^N(SP^d S^N, \mathbb{Z}_2) =
\mathbb{Z}_2$.  Let $\beta$ be the generator of this group. 
According to Nakaoka (\cite{N}), $i^*(\beta)$ is the non-trivial
cohomology class $\omega$ in $H^N(S^N, \mathbb{Z}_2)$.

Let $\tilde \alpha = (A \circ SP^d f)^*(\alpha)$ and let $\tilde \beta = (SP^d g)^* (\beta)$.

In order to compute our pairing, we will prove that $\tilde \alpha ^d = \tilde \beta ^d$.

Since $g$ is degree $1$, we know that $g^* \omega =f^* \alpha$.  
Because of the commutative diagram, $i^* \tilde \alpha = i^* \tilde \beta$.

We write $\tilde \alpha^d - \tilde \beta^d = (\tilde \alpha -
\tilde \beta) (\tilde \alpha^{d-1} + \tilde \alpha^{d-2} \tilde
\beta + ... + \tilde \beta^{d-1}) = (\tilde \alpha-\tilde \beta)
\tilde \alpha^{d-1} + ... + (\tilde \alpha-\tilde \beta)
 \tilde \beta^{d-1}$.  We will prove
that each summand in this formula vanishes, using Lusternik-Schnirelmann 
theory.

We write $SP^d X$ as a union of $d+1$ contractible open sets.  This construction
uses the fact that $X$ is a suspension.  It generalizes the fact that a suspension
is the union of two contractible open sets.  We write $X = SY$ for some space
$Y$, and we think of $SY$ as $[0,1] \times Y$ with each component of the
boundary contracted to a point.  We let $t: X \rightarrow [0,1]$ denote the
projection from $[0,1] \times Y$ onto the first coordinate.  Then for each
point in $SP^d X$ we get an unordered d-tuple of times $<t_1, ..., t_d>$.
For $0 \le k \le d$, we define an open set $U_k \in SP^d X$ to be the set
of points where $k$ of the times $t_i$ are (strictly) greater than $1/(k+2)$ and
the other $d-k$ times are (strictly) less than $1/(k+2)$.  The set $U_k$
is contractible: homotope all the times less than $1/(k+2)$ to zero and all the times
more than $1/(k+2)$ to $1$.  

We prove by induction on $d$ that the sets $U_k$ cover $SP^d X$.  When $d=1$, this is true,
because either $t_1 < 1/2$ and our point is in $U_0$, or else $t_1 > 1/3$ and our point is in $U_1$.  
Now suppose that the result holds for
$d-1$.  If all $t_i$ are more than $1/(d+2)$, then our point lies in $U_d$.  If
not, there is at least one $t_i \le 1/(d+2)$.  Renumber the points so 
that $t_d \le 1/(d+2)$.  Now look at the 
(d-1)-tuple of remaining points.  By induction,
it lies in $U_k$ for some $0 \le k \le d-1$.  In other words, $k$ of the $d-1$
other times are more than $1/(k+2)$ and $d-1-k$ of the $d-1$ remaining times
are less than $1/(k+2)$.  But $t_d \le 1/(d+2) < 1/(k+2)$.  Hence for the original
d-tuple, $k$ of the $d$ times are more than $1/(k+2)$ and the other $d-k$ are less
than $1/(k+2)$.  In other words, our point lies in $U_k$.

Lusternik-Schnirelmann theory immediately implies that any (d+1)-fold cup product
vanishes on $SP^d X$.  In our application, we have to deal with a d-fold cup
product.  The next step is to show that $\tilde \alpha - \tilde \beta$ vanishes
on the union of $U_0$ and $U_1$.  This union is not contractible.  We will show
instead that it contracts to $i(X) \subset SP^d X$.  Since $i$ is an embedding,
and since we checked above that $i^*(\tilde \alpha - \tilde \beta) = 0$, it
will follow that $\tilde \alpha - \tilde \beta$ vanishes on the union of $U_0$
and $U_1$.

Now we check that the union of $U_0$ and $U_1$ contracts to $i(X)$.  If a point
lies in this union, then either it has $d$ times less than $1/2$, or else it has
$d-1$ times less than a $1/3$.  In either case, it has $d-1$ times less than $1/2$.
Let $h_t: [0,1] \rightarrow [0,1]$ be a homotopy with $h_0$ 
the identity, $h_t(0) = 0$ and $h_t(1) = 1$.  We can choose $h_t$ so that $h_1$ maps
$[0,1/2]$ to $0$.  By taking the product with the identity map, we can think of $h_t$
as a homotopy of maps from $[0,1] \times Y$ to $[0,1] \times Y$.  Since each
boundary component is mapped to itself, our homotopy descends to a homotopy of maps
from $X$ to $X$.  Finally, taking the d-fold symmetric product, we get a homotopy
of maps from $SP^d X$ to itself.  The final map, $SP^d h_1$ maps our union into $i(X)$.
The times of $h_1( <x_1, ..., x_d>)$ are $<h_1(t_1), ..., h_1(t_d)>$.  Since
$d-1$ times $t_i$ are less than $1/2$, $d-1$ of the points $h_1(x_i)$ are the basepoint
$x_0$.  Hence $SP^d h_1$ maps our union into $i(X)$.

Since $\tilde \alpha - \tilde \beta$ vanishes on the union of $U_0$ and $U_1$, and since
any cohomology class vanishes on $U_k$ for $2 \le k \le d$, Lusternik-Schnirelmann
theory implies that each cohomology class $(\tilde \alpha - \tilde \beta) \cup \tilde \alpha^e \cup
\tilde \beta^{d-e-1}$ vanishes on $SP^d X$.  Summing these terms, we conclude
that $\tilde \alpha^d = \tilde \beta^d$.

We want to compute the pairing $<(A \circ SP^d f)^* \alpha^d,
[SP^d X]> = <\tilde \alpha^d, [SP^d X]>$.  By our result above,
this pairing is the same as $<\tilde \beta^d, [SP^d X]>$.  We can
evaluate this last pairing by pushing it over to $SP^d S^N$.  It
is equal to $<\beta^d, (SP^d g)_* [SP^d X]>$.  We know that $g$ has
degree 1.  According to Lemma 5.3, $SP^d g$ also has degree 1.   
Hence our pairing is equal to
$<\beta^d, [SP^d S^N]>$.  According to Nakaoka's Theorem, this
pairing is equal to 1.  This finishes our induction on $m$.  We
have now computed the pairing $<Sq^Q a(0,m), F_m(Q)> = 1$.

Using the k-fold suspension map $\Sigma$ we can construct
families of k-cycles.  We define $P_{k,n}(Q)$ to be $P_{n-k}(Q)$. 
Then we define $F_{k,n}(Q): P_{k,n}(Q) \rightarrow Z(k,n)$ to be
$\Sigma \circ F_{n-k}(Q)$.  Since $\Sigma^* a(k,n) = a(0,n-k)$,
it follows that $F_{k,n}(Q)$ detects $Sq_0^{Q_0} ...
Sq_{n-k-1}^{Q_{n-k-1}} a(k,n)$. This finishes the first step of
the proof.

We take a little time to describe the geometry of the cycles in
$F_{k,n}(Q)$. These geometric facts will be used to bound the
volumes of cycles in our final family.  Of course each cycle in
$F_{n-k}(Q)$ is a union of points.  A simple induction argument
shows that the number of points is at most $2^{Q_0} ...
2^{Q_{n-k-1}}$. Besides the number of points in each cycle, we
will need to use some information about the way the points are
arranged.  For every $l$ in the range $0 \le l \le n-k$, each
cycle in $F_{n-k}(Q)$ lies in a union of at most $2^{Q_0} ...
2^{Q_{n-l-1}}$ l-planes, each parallel to the $(x_{n-l+1}, ..,
x_{n-k})$-plane in $\mathbb{R}^{n-k}$.  We verify this claim by
induction on $n-k$.  When $n-k=1$, $F_1(Q_0)$ is the family of
all roots of a degree $2^{Q_0}$ polynomial.  Each set of roots
has at most $2^{Q_0}$ points, and they trivially lie in one line. 
We proceed by induction on $n-k$, assuming that the result holds
for $n-k = m-1$.  First we consider the family $F_m(0, Q_1, ...,
Q_{m-1})$.  Each cycle in this family has the form $\{ t \}
\times C$, where $C$ is a cycle in $F_{m-1}(Q_1, ..., Q_{m-1})$
and $t \in [-1, 1]$.  By induction on $m$, this cycle lies in
$2^{Q_1} ... 2^{Q_{m-l-1}}$ l-planes for each $0 \le l \le m-1$. 
Now we consider the family $F_m(Q_0, ..., Q_{m-1})$ for $Q_0 >
0$.  This cycle is a union of $2^{Q_0}$ cycles of the kind above. 
Therefore, it lies in a union of $2^{Q_0} ... 2^{Q_{m-l-1}}$
l-planes for each $0 \le l \le m-1$.  Also, any of these cycles
trivially lies in one m-plane.

Each cycle in $F_{k,n}(Q)$ has the form $\Sigma(C)$ where $C$ is
a cycle in $F_{n-k}(Q)$.  Therefore, each cycle in $F_{k,n}(Q)$
is a union of at most $2^{Q_0} ... 2^{Q_{n-k-1}}$ k-planes each
parallel to the $(x_{n-k+1}, ..., x_n)$-plane. Moreover, for
every $l \ge k$, each cycle in $F_{k,n}(Q)$ lies in a union of at
most $2^{Q_0} ... 2^{Q_{n-l-1}}$ l-planes, all parallel to the
$(x_{n-l+1}, ..., x_n)$-plane.

Before turning to the bending map $\Psi$, we rotate the family
$F_{k,n}$ to a generic angle, and we dilate it so that we have a
family of relative k-cycles in the ball $B(R)$.  The map $\Psi$
we will choose has displacement bounded by $C(n)$ independent of
$Q$.  We choose $R$ big enough that the displacement of $\Psi$ is
less than $R-1$.  Imitating the proof of Theorem 1, it follows
that the restriction of $\Psi F_{k,n}(Q)$ to the unit ball
detects the cohomology class $Sq^Q a(k,n)$. Now we come to the
heart of the matter: how to choose the bending map $\Psi$ so that
each cycle in $\Psi F_{k,n}(Q)$ has small volume.

The first approach one might try is to use the skeleton-squeezing map $\Psi$ to
push most of the parallel k-planes into the k-skeleton of a lattice with a given
side length.  This is the strategy that we employed in the proof of Theorem 1.
If we use a lattice of side length $s$, then each k-plane is pushed into the
k-skeleton except for a region of final volume at most $C(n) s^k$.  On the other
hand, the total volume of the k-skeleton of side length $s$ is roughly $C(n) s^{k-n}$.
Cancelling the overlaps and optimizing $s$, we end up with a family of cycles
each having volume at most $C(n) 2^{\frac{n-k}{n} \sum Q_i}$.  For most vectors
$Q$, this volume is still much larger than the upper bound we want to prove.
We can improve this approach, because the cycles in $F_{k,n}(Q)$ are not
arbitrary unions of at most $2^{\sum Q_i}$ parallel k-planes.  They also obey
a second property, which we noted above.  For
every $l \ge k$, each cycle in $F_{k,n}(Q)$ lies in a union of at most 
$2^{Q_0} ... 2^{Q_{n-l-1}}$ l-planes, all parallel to the $(x_{n-l+1}, ..., x_n)$-plane.
We will choose a map $\Psi$ that takes advantage of this structure.

Our map $\Psi$ will be a composition of skeleton-squeezing
maps.  The maps occur at different scales and use skeleta of
different dimensions.  We use $\Psi_{[l,s]}$ to denote a squeezing map
to the l-skeleton at scale $s$.  Next we define a sequence of
scales $s_i = \prod_{j=0}^i 2^{- \frac{Q_j}{n-j}}$.  We have $s_0 \ge s_1
\ge s_2 ...$  Our map $\Psi$ is the composition 
$\Psi_{[k, s_{n-k-1}]} \circ ... \circ \Psi_{[n-1, s_0]}$.
Roughly speaking, we first squeeze most of space into 
the (n-1)-skeleton of the lattice
of side length $s_0$.  Then we squeeze most of space into 
the (n-2)-skeleton of the finer
lattice with sidelength $s_1$, and so on.  

Technically, we need to describe the map $\Psi_{[l,s]}$ in more
detail.  First of all, the map depends on a parameter $\epsilon_l > 0$.
As we go along, we will need to choose these parameters so that
$\epsilon_k << \epsilon_{k+1} << ... << \epsilon_{n-1}$.  There is a
second complication that we have to introduce.  In the original definition
of $\Psi$, we used a lattice centered at the origin.  In the course of
our estimates, however, we will need to use a general position
argument.  Therefore, when we construct the map $\Psi_{[l,s]}$, we use
a lattice centered at a generic point instead of the origin.  (For each
$l$, we use a different lattice centered at a different generic point.)

The volume of cycles in $\Psi F_{k,n}(Q)$
is controlled inductively by the following lemma.

\begin{lemma} For each $l$ in the range $k \le l \le n$, each cycle in the family
$\Psi_{[l, s_{n-l-1}]} \circ ... \circ \Psi_{[n-1, s_0]} F_{k,n}(Q)$ lies in a union
of at most $C(n) \prod_{m=0}^{n-l-1} 2^{Q_m}$ pieces of $l$-plane, 
each of diameter at most $C s_{n-l-1}$.
\end{lemma}

\proof The proof is by backwards induction, beginning with the case $l=n$.
When $l=n$, the claim is basically vacuous.  It just asserts that each cycle
in $F_{k,n}(Q)$ lies in 1 piece of n-plane
with diameter on the order of 1.

We now assume the lemma holds for $l$ and we need to prove it for $l-1$.

Let $C_0$ be a cycle in $F_{k,n}(Q)$, and let $C$
be the image $\Psi_{[l, s_{n-l-1}]} \circ ... \circ \Psi_{[n-1, s_0]} C_0$.  We
have to control the geometry of $\Psi_{[l-1, s_{n-l}]} C$.  To do this,
we divide the cycle $C$ into two parts.  We let $G$ be the union of the good
simplices for the map $\Psi_{[l-1, s_{n-l}]}$ and we let $B$ be the union
of the bad simplices.  Then we consider separately $\Psi_{[l-1, s_{n-l}]} (C \cap B)$
and $\Psi_{[l-1, s_{n-l}]} (C \cap G)$.

First we deal with the bad simplices.  We know that $C_0$ is
contained in a union of parallel $l-1$-planes $H_a$ for $a = 1,
..., 2^{Q_0} ... 2^{Q_{n-l}}$. The bad set $B$ lies in the
$\epsilon_{l-1}$-neighborhood of an (n-l)-skeleton $T$. We let
$\bar T$ denote the inverse image $(\Psi_{[l, s_{n-l-1}]} \circ
... \circ \Psi_{[n-1, s_0]})^{-1} (T)$.  Since $T$ is a skeleton
of a lattice whose center is in general position, its inverse
image $\bar T$ is an (n-l)-complex.  (The complex $T$ has
finitely many faces that are each pieces of $(n-l)$-plane.  The
map we are considering is a PL map. We can cut the domain into
finitely many simplices so that the map is linear on each
simplex.  Now for each simplex and each face, we look at
$Map(simplex) \cap face$.  Because the face has been translated
by a generic vector, the intersection has the expected dimension,
and so the inverse image of the given face intersected with the
given simplex lies in a plane of the expected dimension.) 

We let $\bar B$ denote the inverse image $(\Psi_{[l, s_{n-l-1}]}
\circ ... \circ \Psi_{[n-1, s_0]})^{-1}(B)$. Since all $\Psi$ are
PL, we can choose $\epsilon_{l-1}$ sufficiently small that $\bar
B$ is contained in the $C \epsilon_{l-1}$ neighborhood of $\bar
T$.  (Remark: The constant $C$ here depends on the choice of
$\epsilon_l, ..., \epsilon_{n-1}$.)  Since the angle of the
planes $H_a$ is in general position, the intersection $H_a \cap
\bar T$ consists of at most $C(n)$ points.  Using general
position again, the intersection $H_a \cap \bar B$ is contained
in at most $C(n)$ balls of radius at most $C \epsilon_{l-1}$. 
When we apply $\Psi_{[n-1, s_0]}$ to these pieces, we stretch by
a factor $C$, and we may bend each piece into $C(n)$ different
planes.  As long as $\epsilon_{l-1}$ is sufficiently small, a
ball of radius $C \epsilon_{l-1}$ meets at most $C(n)$ simplices
of the triangulations for any of the maps $\Psi_{[l, s_{n-l-1}]},
...,
\Psi_{[n-1, s_0]}$.  Therefore, the image $\Psi_{[l, s_{n-l-1}]}
\circ ... \circ \Psi_{[n-1, s_0]} (H_a) \cap
B$ is contained in a union of at most $C(n)$ pieces of plane each
of diameter at most $C \epsilon_{l-1}$.  Each of these pieces of
plane meets at most $C(n)$ simplices of the triangulation for
$\Psi_{[l-1, s_{n-l}]}$.  In summary, the intersection
$(\Psi_{[l, s_{n-l-1}]} \circ ... \circ \Psi_{[n-1, s_0]} H_a)
\cap B$ is contained in a union of at most $C(n)$ pieces of
plane, each piece of plane lying in a single bad simplex.  Adding
the contributions from all the planes, we see that $C \cap B$ is
contained in a union of at most $C(n) 2^{Q_0} ... 2^{Q_{n-l}}$
pieces of (l-1)-plane, each piece of plane lying in a single bad
simplex.  If $\Delta$ is any simplex in the triangulation for
$\Psi_{[l-1, s_{n-l}]}$, then the image of $\Delta$ has diameter
at most $C(n) s_{n-l}$.  Therefore, the image $\Psi_{[l-1,
s_{n-l}]} (C \cap B)$ is contained in a union of at most $C(n)
2^{Q_0} ... 2^{Q_{n-l}}$ pieces of (l-1)-plane, each of diameter
at most $C(n) s_{n-l}$.

Second we deal with the good simplices.  At this step we use the
inductive hypothesis, which tells us that $C$ is contained in a
union of $C(n) 2^{Q_0} ... 2^{Q_{n-l-1}}$ pieces of l-plane, each
of diameter at most $C(n) s_{n-l-1}$.  The image $\Psi_{[l-1,
s_{n-l}]}(C \cap G)$ is contained in the (l-1)-skeleton of a
lattice of side length $s_{n-l}$.  We also know that the map
$\Psi_{[l-1, s_{n-l}]}$ moves each point at most $C(n) s_{n-l}$. 
If $P$ is a piece of l-plane of diameter at most $C(n)
s_{n-l-1}$, then $\Psi_{[l-1, s_{n-l}]}(P \cap G)$ lies in the
part of the (l-1)-skeleton of side length $s_{n-l}$ within a
distance $C s_{n-l}$ of $P$.  This portion of skeleton can be
covered by $C(n) [s_{n-l-1}/ s_{n-l}]^l$ (l-1)-faces, each of
diameter at most $C(n) s_{n-l}$.  Plugging in the definition of
$s_i$, we see that $[s_{n-l-1}/ s_{n-l}]^l = 2^{Q_{n-l}}$. 
Therefore, $\Psi_{[l-1, s_{n-l}]} (C \cap G)$ can be covered by
$C(n) 2^{Q_0} ... 2^{Q_{n-l}}$ pieces of (l-1)-plane of diameter
at most $C(n) s_{n-l}$. \endproof

In particular, each cycle in $\Psi F_{k,n}(Q) = 
\Psi_{[k, s_{n-k-1}]} \circ ... \circ \Psi_{[n-1, s_0]} F_{k,n}(Q)$ is contained in a union
of $C(n) 2^{Q_0} ... 2^{Q_{n-k-1}}$ pieces of k-plane,
each with diameter at most $C s_{n-k-1}$.  So each cycle has total volume
at most $C(n) 2^{Q_0} ... 2^{Q_{n-k-1}} (s_{n-k-1})^k$.  Plugging in the
definition of $s_{n-k-1}$, we see that each cycle has volume at most
$C(n) 2^{Q_0} ...2^{Q_{n-k-1}} \prod_{j=0}^{n-k-1} 2^{- \frac{k}{n-j} Q_j} 
= C(n) \prod_{i=0}^{n-k-1} 2^{\frac{n-k-i}{n-i} Q_i}$.  Since $\Psi F_{k,n}(Q)$
detects $Sq^Q a(k,n)$, $\mathbb{V}(Sq_0^{Q_0} ... Sq_{n-k-1}^{Q_{n-k-1}}
a(k,n)) \le C(n) \prod_{i=0}^{n-k-1} 2^{\frac{n-k-i}{n-i} Q_i}$. \endproof

\section{Families of algebraic cycles}

In this section, we give some examples of families of cycles
coming from algebraic geometry.  These examples can be used to
give alternate proofs of the upper bounds in certain cases of
Theorem 1 and Theorem 2.
Although they don't cover all the cases covered in the previous
section, these examples are simpler in some ways than the
examples using bent planes.  I had several reasons for
including them.  First, seeing examples of
families of cycles helps to put the results of the paper into
context.  Second, it looks plausible to me that
these examples give the optimal values of $\mathbb{V(\alpha)}$
for certain $\alpha$.  Third, families of complex algebraic
varieties give the only proofs of certain upper bounds for
families of integral cycles, as described in Appendix 2.

There is one important technicality.  It is not trivial to prove
that algebraic families of algebraic cycles actually form
continuous families in $Z(k,n)$.  At the end of the section, we
give a self-contained proof that families of algebraic
hypersurfaces are continuous in $Z(n-1,n)$.  This argument makes
Examples 2, 3, and 5 completely rigorous. 

Right now, Example 4 is not completely rigorous.

\begin{example} The roots of a polynomial.
\end{example}

Let $V(d)$ be the space of all real polynomials of one variable with
degree at most $d$.  The space $V(d)$ is a vector space of
dimension $d+1$.  To each non-zero polynomial in $V(d)$, we
associate its real roots, taken with multiplicity.  This
association defines a map $R_0$ from $V(d) - \{ 0 \}$ to the space
of integral 0-cycles on the real line, but the map is NOT
continuous.  The reason for the discontinuity is that two real
roots may approach each other, become a double root, and then
become two conjugate complex roots.  Since $R_0$ only records the
real roots, two real roots can come together and disappear.

We correct this problem by considering the roots with
multiplicity modulo 2.  We define a root map $R$ from $V(d) - \{
0 \}$ to $Z(0,1)$ by taking the real roots of a polynomial,
keeping only the roots in the interval $(-1, 1)$, and recording
the multiplicity modulo 2.  The map $R$ is continuous.

For any non-zero real number $\lambda$, the polynomials $P$ and
$\lambda P$ have the same roots, and so $R$ induces a map $F(d)$
from $\mathbb{RP}^d = [V(d) - \{ 0 \} ] / \mathbb{R}^*$ to
$Z(0,1)$.  We call this the family of roots of degree d
polynomials.

For example, if $d = 1$, then the map $F(1)$ sends the
polynomial $a x + b$ to its root $-b/a$.  If we fix
$a=1$, then as $b$ goes from $- \infty$ to $+ \infty$,
the point $-b/a$ goes from $+ \infty$ to $- \infty$.  So
the family $F(1)$ sweeps out the unit ball $(-1, 1)$ with
degree 1 modulo 2.  Hence $F(1)^* (a(0,1))$ is the generator
of $H^1(\mathbb{RP}^1)$.

Next we compute that $F(d)^*(a(0,1))$ is the generator of
$H^1(\mathbb{RP}^d)$.  To check this, we pick a homologically
non-trivial curve $c$ in $\mathbb{RP}^d$ and we check that 
$F(c)$ sweeps out the unit interval.  We can take
the curve $c$ given by the projectivization of the linear
polynomials, $V(1) \subset V(d)$.  The map $F(d)$ restricted to this
copy of $\mathbb{RP}^1$ is just $F(1)$, and so the claim
follows from the last paragraph.  Therefore, the family $F(d)$ 
detects $a(0,1)^d$.

\begin{example} Planar real algebraic curves.
\end{example}

Now let $V(d,2)$ denote the vector space of real polynomials in
two variables with total degree at most d.  The vector space
$V(d,2)$ has dimension ${d+2 \choose 2} = (1/2) (d^2 + 3d + 2)$.
We let $\mathbb{RP}^{D(d,2)}$ be the projectivization of $V(d,2)
- \{ 0 \}$.  The dimension $D(d,2) = (1/2) (d^2 + 3d)$.  To each
equivalence class $[P]$ in $\mathbb{RP}^{D(d,2)}$ we can
associate the real algebraic variety defined by $P(x,y) = 0$.  We
define $F(d,2) ([P])$ to be the restriction of this real
algebraic variety to the unit disk, considered as a mod 2
relative Lipschitz cycle.

As in the first example, $F(d,2)^*(a(1,2))$ is the generator of
$H^1(\mathbb{RP}^{D(d,2)})$.  We can see this by the same
argument.  The polynomials of the form $a x_1 + b$ make up a
linear copy of $\mathbb{RP}^1 \subset \mathbb{RP}^{D(d,2)}$.  The
map $F(d,2)$ restricted to this $\mathbb{RP}^1$ gives a family of
parallel vertical lines sweeping out the unit disk.  In other words,
$a(1,2)$ evaluated on $F(d,2)_*([\mathbb{RP}^1])$ is equal to 1, and
so $F(d,2)^*(a(1,2))$ is the generator of
$H^1(\mathbb{RP}^{D(d,2)})$.  Therefore, the family $F(d,2)$
detects $a(1,2)^p$ for all $p \le D(d,2)$.

By a standard argument, we can bound the length of a real
algebraic curve in terms of the degree using the Crofton formula. 
The Crofton formula expresses the length of a curve $C$ in the
plane as an appropriate average of the number of intersections of
$C$ with all the lines in the plane.  Suppose that $C$ is a
degree d algebraic curve intersected with the unit disk, and let
$S^1$ denote the unit circle.  If a line $L$ does not meet the
open disk, then $L$ has zero intersections with $C$. If $L$ does
intersect the open disk, then it intersects $S^1$ twice. Because
$C$ has degree $d$, it intersects almost every line $L$ at most
$d$ times.  By the Crofton formula, the length of $C$ is at most
$(d/2)$ times the length of $S^1$.  So the length of $C$ is at
most $\pi d$.  (I suspect that the sharp constant is $2d$, given
by a union of $d$ lines through the origin, but I don't know how
to prove it.)

This example proves that $\mathbb{V}(a(1,2)^p) \le \pi d$ for any
$p \le D(d,2) = (1/2) (d^2 - 3d)$.  Hence $\mathbb{V}(a(1,2)^p) \le
C p^{1/2}$, proving the upper bound in Theorem 1 in case $k=1$ and
$n=2$.

\begin{example} Real algebraic hypersurfaces.
\end{example}

There is an analogous family of real algebraic hypersurfaces in
any dimension.  Let $V(d,n)$ denote the vector space of real
polynomials in $n$ variables of total degree at most $d$.  We let
$\mathbb{RP}^{D(d,n)}$ be the projective space given by the
quotient $[V(d,n) -
\{ 0 \}] / \mathbb{R}^*$.  The dimension $D(d,n)$ grows like $d^n$
in the sense that $c(n) d^n \le D(d,n)
\le C(n) d^n$ for positive constants $c(n) < C(n)$.  We define
$F(d,n)$ to be the map sending an equivalence class $[P]$ to the
real algebraic variety $P(x_1, ..., x_n) = 0$ intersected with
the unit n-ball.  As in the case of curves, $F(d,n)$ gives a
continuous map from $\mathbb{RP}^{D(d,n)}$ to $Z(n-1,n)$.

The pullback $F(d,n)^*(a(n-1,n))$ is again the generator of
$H^1(\mathbb{RP}^{D(d,n)})$.  The proof is essentially the same
as for curves.  We consider a linear $\mathbb{RP}^1 \subset
\mathbb{RP}^{D(d,n)}$ given by the polynomials of the form $a x_1
+ b$.  The restriction of $F(d,n)$ to this $\mathbb{RP}^1$ is a
family of parallel hyperplanes $x_1 = constant$ sweeping out the
unit n-ball.  Therefore, the pairing of $a(n-1,n)$ with
$F(d,n)_*([\mathbb{RP}^1])$ is equal to 1, proving the claim. 
Therefore, $F(d,n)$ detects $a(n-1,n)^p$ for any $p \le D(d,n)$. 

As for curves, the Crofton formula allows us to bound the volume
of a degree $d$ real algebraic hypersurface in the unit n-ball by
$C(n) d$.  This example proves that $\mathbb{V}(a(n-1,n)^p) \le
C(n) d$ for any $p \le D(d,n)$.  Since $D(d,n) \ge c(n) d^n$, it
follows that $\mathbb{V}(a(n-1,n)^p) \le C(n) p^{1/n}$.  This
proves the upper bound in Theorem 1 in the case that $k = n-1$.

\begin{example} Complex algebraic hypersurfaces.
\end{example}

If the dimension $n$ is even, we can think of $B^n$ as the unit
ball in $\mathbb{C}^{n/2}$, and we can consider families of
complex hypersurfaces.  If we consider the complex hypersurfaces
as mod 2 cycles, we get a family $F_\mathbb{C}(d,n)$, mapping
$\mathbb{CP}^{D_\mathbb{C}(d,n)}$ to $Z(n-2,n)$.  The dimension
$D_\mathbb{C}(d,n)$ grows like $d^{n/2}$.  The pullback
$F_\mathbb{C}(d,n)^*(a(n-2,n))$ is the generator of
$H^2(\mathbb{CP}^{D_\mathbb{C}(d,n)}, \mathbb{Z}_2)$.  Therefore,
this family detects $a(n-2,n)^p$ for $p \le D_\mathbb{C}(d,n) \le
C(n) d^{n/2}$.  By the complex version of the Crofton formula,
each degree d complex hypersurface in the unit ball has volume at
most $C(n) d$.  This example proves that $\mathbb{V}(a(n-2,n)^p)
\le C(n) p^{2/n}$, giving the upper bound in Theorem 1 in case $k
= n-2$ and $n$ is even.

Complex cycles are canonically oriented, so we can also look at
the degree d hypersurfaces as families of {\it integral} cycles. 
We discuss integral cycles in Appendix 2.

\begin{example} Products of the previous examples.
\end{example}

By taking Cartesian products we can produce a variety of new
families of cycles.  For example, $F(d_1, n_1) \times F(d_2,
n_2)$ is a map from $\mathbb{RP}^{D(d_1,n_1)} \times
\mathbb{RP}^{D(d_2, n_2)}$ to $Z(n_1 + n_2 - 2, n_1 + n_2)$.  The
pullback of $a(n_1 + n_2 -2, n_1 + n_2)$ is given by the tensor
product $\omega_1 \otimes \omega_2$, where $\omega_i$ is the
generator of $H^1(\mathbb{RP}^{D(d_i,n_i)})$.  

Most of these families have unnecessarily large maximal volumes,
and so they are not useful for proving the upper bounds in
Theorem 1 and Theorem 2.  We focus on only a few cases, where the
maximal volume is near to optimal.  The cases involve products
with the trivial family which consists of a point moving across
an interval.  This family is called $F(1,1)$, the family that
gives the root of a linear polynomial in one variable, described
in Example 1 above.  (We called it $F(1)$ at the time.)

Consider the product $F(d,n) \times F(1,1)$, which is a family of
cycles in $Z(n-1, n+1)$ parametrized by $\mathbb{RP}^{D(d,n)}
\times S^1$.  Let $\omega$ denote the generator of
$H^1(\mathbb{RP}^{D(d,n)})$, and let $\beta$ denote the generator
of $H^1(S^1)$.  Then the pullback of $a(n-1,n+1)$ by our family is
$\omega \otimes \beta$.

This family detects certain towers of Steenrod squares applied to
$a(n-1,n+1)$.  More precisely, if $2^Q \le D(d,n)$, then this
family detects $Sq_1^Q a(n-1,n+1)$.  In order to see this, consider
the smash product map $\mathbb{RP}^{D(d,n)} \times S^1
\rightarrow S \mathbb{RP}^{D(d,n)}$.  The cohomology class
$\omega \times \beta$ is the pullback of the generator of $H^2(S
\mathbb{RP}^{D(d,n)})$.  Because Steenrod squares commute with
suspensions, it follows that $Sq_1^Q (\omega \otimes \beta) =
(Sq_0^Q \omega) \otimes \beta = \omega^{2^Q} \otimes \beta$. 
Since $2^Q \le D(d,n)$, $\omega^{2^Q}$ is a non-vanishing class
in $H^*(\mathbb{RP}^{D(d,n)})$ and $\omega^{2^Q} \otimes \beta$
is a non-vanishing class in $H^*(\mathbb{RP}^{D(d,n)} \times
S^1)$.

A cycle in $F(d,n) \times F(1,1)$ is just a product of a cycle in
$F(d,n)$ with a point, and so it has volume at most $C(n) d$. 
This example proves that $\mathbb{V}(Sq_1^Q a(n-1,n+1)) \le C(n)
d$ as long as $2^Q \le D(d,n) \le C(n) d^n$.  In other words, it
shows that $\mathbb{V}(Sq_1^Q a(n-1,n+1)) \le C(n) 2^{Q/n}$.  If
we write $N = n+1$, then we get the inequality $\mathbb{V}(Sq_1^Q
a(N-2, N)) \le C(N) 2^{\frac{1}{N-1} Q}$.  This
inequality is the upper bound in Theorem 2 in the special case
that $k = N-2$ and $Q_i = 0$ for $i \not= 1$.

More generally, we can consider the product $F(d,n) \times F(1,1)
\times ... \times F(1,1) = F(d,n) \times F(1,1)^s$.  In this
case, we get a family of cycles in $Z(n-1, n+s)$ parametrized by
$\mathbb{RP}^{D(d,n)} \times T^s$, where $T^s$ denotes the
s-dimensional torus $(S^1)^s$.  The pullback of $a(n-1, n+s)$ is
given by $\omega \otimes \beta$ for $\omega$ the generator of
$H^1(\mathbb{RP}^{D(d,n)})$ and $\beta$ the generator of
$H^s(T^s)$.  By the same argument as above, this family detects
the class $Sq_s^Q a(n-1, n+s)$ for $2^Q \le D(d,n)$.  Each cycle
in the family has volume at most $C(n) d$.  This proves the
inequality $\mathbb{V}(Sq_s^Q a(n-1, n+s)) \le C(n) 2^{Q/n}$,
giving some more special cases of the upper bounds in Theorem 2.

\vskip5pt

We now return to the technical problem of showing that algebraic
cycles are flat cycles and that the families we mentioned above
are continuous families of flat cycles.  We will prove that the
family of real algebraic hypersurfaces $F(d,n)$ defined in
Example 3 is genuinely a continuous family of mod 2 flat
(n-1)-cycles.  As a special case, it follows that the family
$F(d,2)$ in Example 2 is a continuous family of flat cycles, and
it also follows that the products in Example 5 are continuous
families of flat cycles.  This material is probably old, but I
don't know a reference for it.

Let us start by considering the following special case.  Fix a
polynomial $P \not= 0$, and consider the two polynomials $P -
\delta$ and $P + \delta$ for a small real number $\delta$.  If
$\delta$ is small, the two polynomials are close together in the
space of all polynomials.  The vanishing sets of these two
polynomials are $P^{-1}(\delta)$ and $P^{-1}(-\delta)$.  These
two cycles bound a chain $P^{-1}([-\delta, \delta])$.  We have to
show that if $\delta$ is small, then this chain has small volume. 
The following lemma gives such a bound.

\begin{lemma} Suppose that $P$ is a real polynomial of degree at
most $d$ in $n$ variables.  We can write $P$ in multi-index
notation as $P(x) = \sum_{|I| \le d} c_I x^I$.  Suppose that
$\max |c_I| = M$.  Then the set $\{ x \in B(1) | |P(x)| \le
\delta \}$ has volume at most $C(d,n) [M^{-1} \delta]^{\frac{1}{dn}}$.
\end{lemma}

Remark: This estimate is far from sharp.  For our purposes, it 
doesn't matter what power of $\delta$ appears in the estimate.

\proof By scaling, it suffices to prove the theorem when $M=1$.
We first consider the trivial case that $\max_{I \not= 0} |c_I| \le (1/2) n^{-d}$.
In this case, we must have $|c_0| = 1$.  But then
$|P(x)| \ge 1/2$ everywhere, and so our inequality holds
automatically.  We may assume that $\max_{I \not=0} |c_I| \ge c(d,n)$.
By scaling, it suffices to prove that if $\max_{I \not= 0} |c_I| = 1$, then the
volume of $P^{-1}([-\delta, \delta])$ is at most $C(d,n) \delta^{\frac{1}{dn}}$.

We proceed by induction on $n$.  First we do the case $n=1$.

Divide the 1-dimensional ball $[-1, 1]$ into at most $d$ segments
$S_i$ so that $P$ is monotonic on each segment.  Let $\bar S_i :=
S_i \cap P^{-1}([- \delta, \delta])$. Because $P$ is monotonic on
$S_i$, each $\bar S_i$ is a single segment.  The oscillation of
$P$ on $\bar S_i$ is at most $2 \delta$.  Therefore, it suffices
to check that the oscillation of $P$ on a segment $S \subset [-1,
1]$ of length $|S|$ is at least $C(d) |S|^d$.

We write $P(x) = \sum_{i=0}^d c_i x^d$.  The $d^{th}$ derivative
$P^{(d)}$ is a constant $d! c_d$.  Therefore, the oscillation of
$P^{(d-1)}$ on a segment $S$ is $d! |c_d| |S|$.  We can choose a
subsegment $S_1 \subset S$ with $|S_1| \ge (1/4) |S|$ so that on
$S_1$, $|P^{(d-1)}|$ is constant up to a factor of 2 with size on
the order of $|c_d| |S|$.  (In particular, the sign of
$P^{(d-1)}$ is constant on $S_1$.)  Therefore, the oscillation of
$P^{(d-2)}$ on $S_1$ is on the order of $|c_d| |S|^2 + |c_{d-1}|
|S|$. We can choose a subsegment $S_2 \subset S_1$ with $|S_2|
\ge (1/8) |S_1|$ so that on $S_2$, $|P^{(d-2)}|$ is constant up
to a factor of 2.  Therefore, the oscillation of $P^{(d-3)}$ on
$S_2$ is on the order of $|c_d| |S|^3 + |c_{d-1}| |S|^2 +
|c_{d-2}| |S|$.  Continuing in this way, we eventually produce a
subsegment of $S$ where the oscillation of $P$ is at least $c(d)
[|c_d| |S|^d + ... + |c_1| |S|]$.  In particular the oscillation
is at least $c(d) |S|^d$. This proves our lemma in case $n=1$.

Now we turn to the inductive step.  We assume that the lemma
holds for $n-1$. Let $P$ be a polynomial in $n$ variables of
degree $d$.  As we described at the beginning of the proof, we
can assume that $|c_I| = 1$ for a non-zero index $I$.  By
reordering the variables, we can assume that $x_n$ divides $x^I$.
We sort the monomials of $P$ according to the power of $x_n$,
writing $P(x) = \sum_{i=0}^d P_i(x_1, ..., x_{n-1}) x_n^i$. 
Here, each $P_i$ is a polynomial of degree at most $d$ in the
first $n-1$ variables.  By assumption, one of the $P_i$ with $i
\not= 0$ has a coefficient with norm $1$.  Using induction, we
can apply the lemma to the polynomial $P_i$.  We conclude that
$|P_i| \ge \beta$ except for a set $B$ of measure $C(d,n)
\beta^{\frac{1}{(n-1) d}}$, for a number $\beta$ that we can
choose later. Suppose $(x_1, ..., x_{n-1})$ is not in $B$. 
Fixing this choice of $(x_1, ..., x_{n-1})$, we consider $P$ as a
polynomial in the one variable $x_n$.  This polynomial has
largest coefficient of size at least $\beta$.  By applying the
lemma to this 1-dimensional polynomial, we conclude that
$|P(x_n)| \ge \delta$ except for a subset of values of $x_n$ of
volume at most $C(d,n) (\beta^{-1} \delta)^{1/d}$.  We take
$\beta = \delta^{\frac{n-1}{n}}$, finishing the induction.
\endproof

With this lemma in hand, we can prove the continuity of $F(d,n)$.
Let $P$ be a real degree d polynomial in n variables.  We write
$P = \sum c_I x^I$, and we define $\| P \| = [\sum
|c_I|^2]^{1/2}$.  We consider the family of all $P$ of norm 1. 
This family is parametrized by the sphere $S = \{ c_I | \sum
|c_I|^2 = 1 \}$.  Let $V_P$ denote the set $\{ x \in B^n(1) |
P(x) = 0\}$.  For an open set of full measure in $S$, the variety
$V_P$ is a smooth hypersurface which meets the boundary of $B^n(1)$
transversely.  In this case, $V_P$ clearly defines a mod 2
Lipschitz cycle and hence a mod 2 flat cycle. We let $S_0 \subset
S$ denote this open set of full measure.  We define a map
$F(d,n): S_0 \rightarrow Z(n-1,n)$ by mapping $P$ to $V_P$
considered as a flat cycle.  We note that $V_P = V_{-P}$. 
Therefore, if $P$ is contained in $S_0$, then so is $-P$, and
$F(d,n)(P) = F(d,n)(-P)$.

\begin{lemma} If $P$ and $Q$ are in $S_0$, then the area-distance
from $V_P$ to $V_Q$ is at most $C(d,n) \| P - Q
\|^{\epsilon(d,n)}$, for a constant $\epsilon(d,n) > 0$.
\end{lemma}

In particular, the map $F(d,n)$ extends to a continuous map from
$S$ to $Z(n-1,n)$.  Since $F(d,n)(P) = F(d,n)(-P)$, the same
holds for the extension, and so we can take a quotient by the
action of $\mathbb{Z}_2$, giving a map $F(d,n):
\mathbb{RP}^{D(d,n)} \rightarrow Z(n-1,n)$.  This $F(d,n)$ is
the family of degree d real algebraic hypersurfaces described
in Example 3 above.

\proof Let $F(x,t) = (1-t) P + t Q$.  We let $V_F = \{ (x,t) \in
B^n(1) \times [0,1] | F(x,t) = 0 \}$.  For each $P \in
S_0$, for almost every $Q \in S_0$, the set $V_F$ is a smooth
manifold meeting the boundary of $B^n(1) \times [0,1]$
transversely.  For the time being, we assume that this is true
for $P$ and $Q$.  Then $V_F$ is a relative n-cycle in $B^n(1)
\times [0,1]$.  Let $\pi: B^n(1) \times [0,1] \rightarrow B^n(1)$
be the projection onto the first factor.  Then $\pi(V_F)$ is a
relative Lipschitz n-chain in $B^n(1)$ with boundary $V_P - V_Q$. 
The area-distance from $V_P$ to $V_Q$ is at most the volume of
$\pi(V_F)$.  In the rest of the proof, we give an upper bound for
this volume.

Let $A$ be the image of $V_F$ under the map $\pi$, considered as
a subset of $B^n(1)$.  Notice that for almost every point of $a
\in A$, $\pi^{-1}(a)$ is a single point in $V_F$.  Therefore, it
suffices to bound the volume of $A$.  We know that the volume of
$V_F$ is at most $C(n) d$ by an argument using the Crofton
formula.  Since $\pi: V_F \rightarrow A$ is 1-to-1 over almost
every point of $A$, we can view almost all of $V_F$ as the graph
on a function $T: A \rightarrow [0,1]$.  Therefore we get the
following integral inequality.

$$\int_A (1 + |\nabla T|^2)^{1/2} \le C(n) d.$$ 

At a point $x$, $\nabla T = - \nabla_x F(x, T(x)) / \partial_t
F(x, T(x))$. We now give bounds for $\nabla F$.  Let $\delta :=
\| P - Q \|$.  The time derivative $\partial_t F(x, t) = - P(x) +
Q(x)$, and so $|\partial_t F(x,t)| \le C(d,n) \delta$.  On the
other hand, the space derivative $\nabla_x F(x,t)$ is equal to
$(1-t) \nabla P + t \nabla Q = \nabla P + t \nabla (Q-P)$.  Now,
the coefficients of $\nabla (Q-P)$ are at most $C(d,n) \delta$,
and so $| t \nabla (Q-P)| \le C(d,n) \delta$.  Therefore,
$|\nabla_x F(x,t) - \nabla P(x)| \le C(d,n) \delta$.  Plugging
these bounds into our last equation, we get the following.

$$\int_A \delta^{-1} |\nabla P| \le C(d,n).$$

Our next idea is to apply the previous lemma to control the
volume of the set where $|\nabla P|$ is small. We have assumed
that $\| P \| = 1$.  As a corollary of the last lemma, the
oscillation of $P$ is at least $c(d,n)$.  Therefore, the maximum
of $|\partial_i P|$ must be at least $c(d,n)$ for some $i$.  We
fix this choice of $i$.  Now $\partial_i P$ is a polynomial of
degree d-1 with norm $\| \partial_i P \| \ge c(d,n) > 0$. 
Applying the lemma, we conclude that the set where $|\partial_i
P|$ is at most $\delta^{1/2}$ has volume at most $C(d,n)
\delta^{\frac{1}{2dn}}$.  Using the last equation, we see that
the volume of $A$ is at most $C(d,n) \delta^{\frac{1}{2dn}} +
C(d,n) \delta^{1/2}$.

So far we have assumed that $V_F$ was non-singular and transverse
to the boundary.  For each $P$, this assumption holds
for almost every $R$, but it may not hold for $Q$.  Hence we
bound the area-distance from $P$ to almost every $R$.  Similarly,
we bound the area-distance from $Q$ to almost every $R$.  Now we
may choose a polynomial $R$ which is good for both $P$ and $Q$
and with $\| P - R \| + \|R - Q \| \le 2 \| P - Q \|$. \endproof

\section{Appendix 1: The standard definition of flat cycles}

In this appendix, we recall the standard definition of the space
of flat cycles.  We check that it is equivalent to the definition
given in Section 1.  The standard definition is more complicated
than the definition we gave in Section 1.
The reason that the standard definition is worthwhile is
that one constructs not only a space of flat cycles but a chain
complex of flat chains, containing the flat cycles as a subset. 
In other contexts, especially in the theory of minimal surfaces,
one is interested in the whole chain complex, but in this paper
we were concerned only with the space of cycles.

As in Section 1, we let $I_{rel}(k,n)$ denote the complex of
relative, mod 2 Lipschitz k-chains in the unit n-ball.  We define
a distance function on this space of k-chains as follows.  We
define the flat norm of a k-chain $C$ to be the infimum, over all
(k+1)-chains $D$ of $|D| + |\partial D - C|$.  We need to make
a remark about this formula in the context of relative
chains.  We view $D$ and $C$ as relative chains, and so we also
view $\partial D - C$ as a relative chain, and we measure its
volume as a relative chain.  For example, if $C$ is the line
$x=0$ in the unit disk and $D$ is the relative 2-chain given by
$x \ge 0$ in the unit disk, then $\partial D - C$ is zero and
has volume zero.  If we were working with absolute chains
instead of relative ones, then $\partial D - C$ would be a
semicircle with volume $\pi$.

In the figure below, we illustrate the sum $|D| + |\partial D -
C|$ for a 1-chain $C$ and a 2-chain $D$ in the unit disk.

\includegraphics{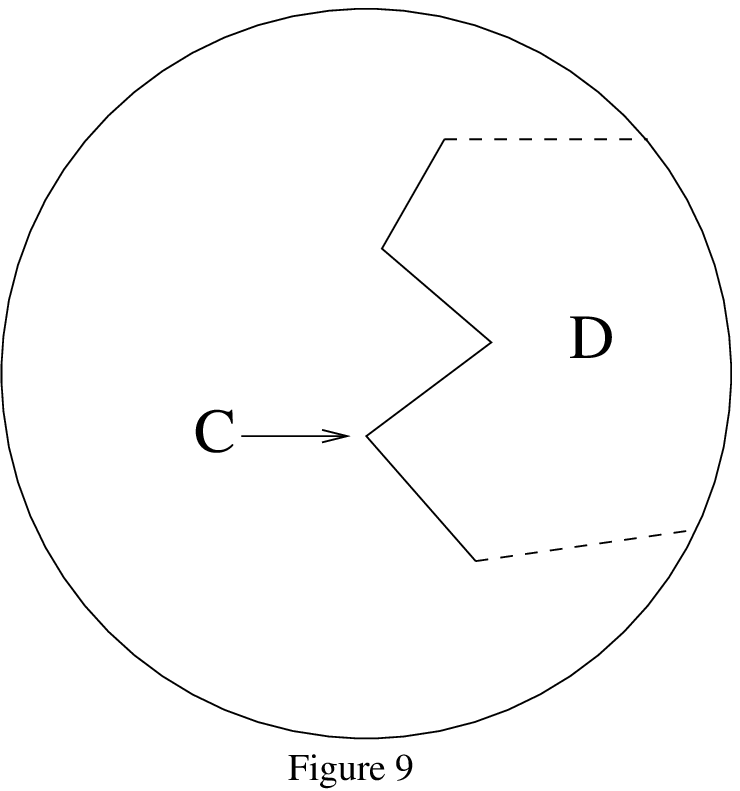}

\noindent The 1-chain $C$ is shown in a solid line.  The 2-chain
$D$ is the region enclosed by $C$ and by the dotted lines.  The
volume of $|\partial D - C|$ is the sum of the lengths of the
two dotted lines.  Because we are working with relative chains, 
it does not include the length of the arc of
the circle bordering $D$.

The flat distance between two chains $C_1$ and $C_2$ is defined
to be the flat norm of $C_1 - C_2$.  We say that $C_1$ and $C_2$
are equivalent if the flat distance between them is zero.  The
set of equivalence classes of relative Lipschitz chains is a
metric space.  We define the completion of this metric space to
be the space of relative flat k-chains in the unit n-ball,
$I_{rel, flat}(k,n)$.

We say that a flat chain $C$ has volume less than $V$ if $C$ is a
limit of relative Lipschitz chains $C_i$ with volume less than
$V$.

If $C$ is a relative Lipschitz k-chain, then the flat norm of its
boundary is at most the flat norm of $C$.  Suppose that $|D| +
|\partial D - C|$ is close to the flat norm of $C$.  We use the
k-chain $C - \partial D$ to bound the flat norm of $\partial C$. 
We get that the flat norm of the boundary of $C$ is at most
$|C - \partial D| + |\partial(C - \partial D) - \partial C|$. 
The second term is zero, and so the flat norm of $\partial C$ is
at most $|C - \partial D|$, which is at most $|C - \partial D| +
|D|$, which is within $\epsilon$ of the flat norm of $C$. 
Therefore, the flat distance between $\partial C_1$ and $\partial
C_2$ is at most the flat distance between $C_1$ and $C_2$. 
Therefore, the boundary operation extends to a continuous map
from $I_{rel, flat}(k+1, n)$ to $I_{rel, flat}(k,n)$.

Now we can define the space of flat cycles $Z_{flat}(k,n)$ to be
the space of cycles in the complex $I_{rel, flat}(k,n)$.  We give
$Z_{flat}(k,n)$ the subspace topology inherited from $I_{rel,
flat}(k,n)$.

The goal of this appendix is to check that the space
$Z_{flat}(k,n)$ is exactly the space $Z(k,n)$ defined in Section 1.

The main point here is that if $C_1$ and $C_2$ are two relative
Lipschitz k-cycles in the unit n-ball, then the flat distance
between them is equal to the area-distance between them.  The
area distance is the infimum over all Lipschitz chains $D$ with
$\partial D = C$ of $|D|$.  On the other hand, the flat distance
is the infimum over all Lipschitz chains $D$ of $|D| + |\partial
D - C|$.  It follows that the flat distance is at most the area
distance.  On the other hand, suppose that $D$ is any relative
Lipschitz k-chain and that $E$ is the cone over $C - \partial D$
with vertex at the origin.  Then $E$ is a relative Lipschitz
(k+1)-chain with volume at most $(k+1)^{-1} |\partial D - C| \le
|\partial D - C|$.  Then $\partial (D + E) = C$ and $|D + E| \le
|D| + |\partial D - C|$.  Therefore, the flat distance is equal
to the area distance.

Now we can check that $Z(k,n)$ is contained in $Z_{flat}(k,n)$. 
Certainly each Lipschitz k-cycle is contained in $I_{rel}(k,n)$.
Hence
$Z(k,n)$ is exactly the subset of $Z_{flat}(k,n)$ which is given
by limits of Lipschitz k-cycles.

We have one more point to check.  If $z$ is a cycle in
$Z_{flat}(k,n)$, it is not apriori clear that $z$ can be written
as a limit of Lipschitz k-cycles.  By definition, a cycle $z \in
Z_{flat}(k,n)$ is a limit of relative Lipschitz {\it chains}
$C_i$.  The fact that $\partial z = 0$ implies $\partial C_i$
converges to zero in the flat topology, but it doesn't imply that
$\partial C_i = 0$ for any finite $i$.  For each $C_i$, we can
choose a k-chain $D_i$ so that $|D_i| + |\partial D_i - \partial
C_i| \rightarrow 0$.  By the isoperimetric inequality, we can
find another k-chain $E_i$ with $\partial E_i = \partial C_i -
\partial D_i$ and with $|E_i| \rightarrow 0$.  We let $\tilde C_i
= C_i - D_i - E_i$.  Each $\tilde C_i$ is a Lipschitz cycle.  On
the other hand, the flat distance from $\tilde C_i$ to $C_i$ is
at most $|D_i| + |E_i| \rightarrow 0$. Therefore, $\tilde C_i$
converges to $z$ in the flat topology.

In summary, we have shown that $Z_{flat}(k,n)$ and $Z(k,n)$ are
the same underlying set equipped with the same metric and hence
the same topology.

\section{Appendix 2: Integral cycles}

The minimax problem described in this paper makes equally good
sense for integral cycles in place of mod 2 cycles.  The lower
bounds in this paper generalize to families of integral cycles. 
Using families of complex algebraic cycles, we can prove matching
upper bounds for some cohomology classes.  For most cohomology
classes, however, there is a large gap between the best upper and
lower bounds that we can prove.

Let $Z_{\mathbb{Z}}(k,n)$ denote the space of integral relative
k-cycles in the unit n-ball.  For a cohomology class $\alpha$ in
$H^*(Z_{\mathbb{Z}}(k,n))$, we define $\mathbb{F}(\alpha)$ to be
the set of all families of integral cycles, $F: P
\rightarrow Z_{\mathbb{Z}}(k,n)$, that detect the cohomology
class $\alpha$ in the sense that $F^*(\alpha)$ is non-zero in
$H^*(P)$.  (This definition makes sense for any choice of
coefficients in the cohomology groups.)  Then we define the
minimax volume $\mathbb{V}(\alpha)$ by the following formula.

$$\mathbb{V}(\alpha) = \inf_{F \in \mathbb{F}(\alpha)} \sup_{C
\in F} mass(C).$$

Using the methods in Section 1, we can construct a fundamental
cohomology class $a_{\mathbb{Z}}(k,n)$ in
$H^{n-k}(Z_{\mathbb{Z}}(k,n), \mathbb{Z})$.  We can then reduce
this class to either $H^{n-k}(Z_{\mathbb{Z}}(k,n), \mathbb{Q})$
or $H^{n-k}(Z_{\mathbb{Z}}(k,n), \mathbb{Z}_2)$.  We denote these
cohomology classes $a(k,n, \mathbb{Z}), a(k,n, \mathbb{Q}),$ and
$a(k,n, \mathbb{Z}_2)$.

All of the lower bounds in the paper apply without modification
to the space of integral cycles.

\begin{Theorem 1A} (Lower bounds) Let $a(k,n)$ be short for any
of the classes $a(k,n, \mathbb{Z}),$ $a(k,n, \mathbb{Q}),
$ or $a(k,n, \mathbb{Z}_2)$.  Then the following lower bound holds.

$$\mathbb{V}(a(k,n)^p) \ge c(n) p^{\frac{n-k}{n}}.$$

\end{Theorem 1A}

\begin{Theorem 2A} (Lower bounds) For each $\epsilon > 0$, there is a constant
$c(n, \epsilon) > 0$ so that the following estimate holds.

$$\mathbb{V}(Sq_0^{Q_0} ... Sq_{n-k-1}^{Q_{n-k-1}} a(k,n,
\mathbb{Z}_2)) \ge c(n, \epsilon) \prod_{i=0}^{n-k-1} (2 -
\epsilon)^{\frac{n-k-i}{n-i} Q_i} .$$

\end{Theorem 2A}

Remark: If Almgren's arguments apply to the space of flat cycles,
which they probably do, then $Z_{\mathbb{Z}}(k,n)$ is weak
homotopic to an Eilenberg-Maclane space $K(\mathbb{Z}, n-k)$. 
The mod 2 cohomology ring of $K(\mathbb{Z}, n-k)$
was determined by Serre in \cite{S}.  It is a free algebra over
$\mathbb{Z}_2$ generated by $Sq_1^{Q_1} ... Sq_{n-k-2}^{Q_{n-k-2}} a(k,n,
\mathbb{Z}_2)$.

We can prove some interesting upper bounds by using families of
complex algebraic cycles.  These upper bounds are not completely
rigorous, because we haven't checked that the families of
algebraic cycles are continuous in the space
$Z_{\mathbb{Z}}(k,n)$.  The continuity sounds true to me. 
Related results are mentioned in \cite{L}, for example.

Example 1. Suppose that $n$ is even and that $a(n-2,n)$ is short
for any of the classes $a(n-2,n, \mathbb{Z}), a(n-2,n,
\mathbb{Z}_2),$ or $a(n-2,n, \mathbb{Q})$.  Then
$\mathbb{V}(a(n-2,n)^p) \le C(n) p^{2/n}$.  This upper bound
matches the lower bound in Theorem 1A up to a constant factor.

\proof We have $n=2m$, and so we can view the unit ball $B^n$ as
the unit ball in complex space $\mathbb{C}^m$.  We consider the
family of complex hypersurfaces of degree at most $d$. This
family is parametrized by $\mathbb{CP}^{N_m(d)}$, where $N_m(d)
\sim d^m$.  It detects the class $a(n-2,n)^p$ for all $p \le
N_m(d)$.  Each complex hypersurface intersected with the unit
n-ball has volume at most $C(n) d$.  Hence
$\mathbb{V}(a(n-2,n)^p) \le C(n) d$ for $p \le c(n) d^m$, and so
$\mathbb{V}(a(n-2,n)^p) \le C(n) p^{2/n}$.
\endproof

Example 2. Suppose that $n$ is even.  Then $\mathbb{V}(Sq_s^Q
a(n-2, n+s, \mathbb{Z}_2)) \le C(n) 2^{\frac{2}{n} Q}$. 
According to Theorem 2A, this minimax volume is at least $c(n,
\epsilon) (2- \epsilon)^{\frac{2}{n} Q}$.  Hence the upper and
lower bounds match up to a factor that grows sub-exponentially in
$Q$.

\proof The proof is essentially the same as in Example 5 of
Section 6.  Let $F(1,1)$ denote a family of 0-cycles in the unit
1-ball consisting of a single point that moves from one end of
the interval to the other.  We can think of $F(1,1)$ as a family
of integral 0-cycles that detects $a(0,1, \mathbb{Z})$ and also
$a(0,1, \mathbb{Z}_2)$.  Take the product of $F(1,1)^s$ with
Example 1 above. \endproof

For many other cases, there is a big gap between the best upper
and lower bounds.  The simplest example concerns 1-cycles in the
unit 3-ball.  Our best lower bound for $\mathbb{V}(a(1,3,
\mathbb{Z})^p)$ is $c p^{2/3}$.  The only upper bound that I know
for $\mathbb{V}(a(1,3,\mathbb{Z})^p)$ is $C p$, which we get by
considering families of p vertical lines.

\section{Appendix 3: Minimax volumes of Riemannian manifolds}

Minimax volumes analogous to those we have studied can also be
defined using a Riemannian manifold $(M^n,g)$ in place of the
unit n-ball.

Let $Z(k,M)$ denote the space of absolute mod 2 k-cycles in $M$. 
If $M$ has a boundary, then let $Z_{rel}(k, M)$ denote the space
of relative mod 2 k-cycles in $(M, \partial M)$. The construction
of the fundamental cohomology class $a(k,n)$ generalizes to the
setting of manifolds.  If $M$ is a closed n-manifold, then we get
a fundamental cohomology class $a(k,M) \in H^{n-k}(Z(k,M))$.  If
$M$ is a compact n-manifold with boundary, then we get a
fundamental cohomology class $a(k,M) \in H^{n-k}(Z_{rel}(k,M))$.

If $\alpha$ is any cohomology class in $H^*(Z(k,M))$, then we can
define a minimax volume associated to $\alpha$.  We let
$\mathbb{F}(\alpha)$ denote the set of all families of cycles $F:
P \rightarrow Z(k,M)$ that detect $\alpha$.  Then we define a
minimax volume $\mathbb{V}_{(M, g)}(\alpha)$ by the usual
formula.

$$\mathbb{V}_{(M, g)}(\alpha) = \inf_{F \in \mathbb{F}(\alpha)} \sup_{C \in F} \textrm{ Vol}(C).$$

\noindent We need the metric $g$ in order to measure the volumes of k-cycles in $M$.  For
a fixed choice of $M$ and $\alpha$, the minimax volume $\mathbb{V}_{(M, g)}(\alpha)$ is
a function of $g$.  Roughly speaking, this function measures how large the manifold $(M,g)$ is.
If $g \ge h$, then it's easy to check that $\mathbb{V}_{(M,g)}(\alpha) \ge \mathbb{V}_{(M,h)}(\alpha)$.

These minimax volumes can be used to control the geometry of degree 1 maps between
Riemannian manifolds.  If $\Phi: M \rightarrow N$ is a Lipschitz map, then it induces a 
map $Z_k(\Phi)$ from $Z(k, M)$ to $Z(k, N)$ for every k.  Similarly, if $\Phi: 
(M, \partial M) \rightarrow (N, \partial N)$ is
a Lipschitz map of pairs, then it induces a map $Z_k(\Phi)$ from $Z_{rel}(k, M)$ to 
$Z_{rel}(k, N)$.  In either case, the pullback $Z_k(\Phi)^* 
[a(k,N)] = (\textrm{deg } \Phi) a(k,M)$.  This equation
follows directly from the construction of $a(k, M)$.

If $\Phi: M \rightarrow N$ is a piecewise $C^1$ map, then we 
say that the k-dilation of $\Phi$ is at most $\Lambda$
if $\Phi$ maps each k-dimensional submanifold of $M$ with 
volume $V$ to an image with volume at most $\Lambda V$.

\begin{prop} Suppose that $\Phi: (M,g) \rightarrow (N,h)$ is a $C^1$ map with degree 1 mod 2 and
k-dilation $\Lambda$.  Suppose that $O$ denotes any natural cohomology operation with mod 2
coefficients.  In particular, $O$ may denote a cup power or any tower of Steenrod squares.
Then the following inequality holds.

$$\Lambda \mathbb{V}_{(M,g)}(O a(k,M)) \ge \mathbb{V}_{(N,h)}(O a(k, N)).$$

\end{prop}

\proof For any $\delta > 0$, let $F_\delta$ be a family of k-cycles in $M$ that detects $O a(k,M)$ 
with volume at most $\mathbb{V}_{(M,g)}(O a(k,M)) + \delta$.  Then $\Phi(F_\delta)$ is a family of k-cycles in 
$(N,h)$ that detects $O a(k,N)$ with volume at most $\Lambda [\mathbb{V}{(M,g)}(O a(k,M)) +
\delta]$. \endproof

This proposition gives a large number of lower bounds for the k-dilation $\Lambda$ of a degree 1
map from $(M,g)$ to $(N,h)$.  If $N$ has a boundary, then we can formulate a slightly more
general version of this proposition.

\begin{prop} Suppose that $(M^n,g)$ is a compact n-manifold and that $(N,h)$ is a compact
n-manifold with boundary.  Suppose that $U \subset M$
is an open set with a piecewise smooth boundary.  Suppose that $\Phi: (U, \partial U) \rightarrow (N, \partial N)$
is a piecewise $C^1$ map with degree 1 mod 2 and k-dilation $\Lambda$.  Let $O$
be a cohomology operation.  Then the following
inequality holds.

$$\Lambda \mathbb{V}_{(M,g)}(O a(k,M)) \ge \mathbb{V}_{(N,h)}(O a(k, N)).$$

\end{prop}

\proof For any $\delta > 0$, let $F_\delta$ be a family of k-cycles in $M$ that detects $O a(k,M)$ 
with volume at most $\mathbb{V}_{(M,g)}(O a(k,M)) + \delta$.  Let $\hat F_\delta$ be the restriction
of this family to $U$.  The family $\hat F_\delta$ detects $O a(k,U)$ with volume at most $\mathbb{V}_{(M,g)}(O a(k, M))$.  Then $\Phi(\hat F_\delta)$ is a family of k-cycles in 
$(N,h)$ that detects $O a(k,N)$ with volume at most $\Lambda [\mathbb{V}_{(M,g)}(O a(k,M)) +
\delta]$. \endproof

In my thesis \cite{Gu2}, I studied the k-dilation of degree 1 maps between various sets in Euclidean
space.  The main theorem of the thesis is the following.

\begin{reftheorem} (\cite{Gu2}, \cite{Gu3}) 
Suppose that $R$ and $S$ are n-dimensional rectangles with dimensions $R_1 \le ...
\le R_n$ and $S_1 \le ... \le S_n$ respectively.  Suppose that $U \subset R$ is an open set with
a piecewise $C^1$ boundary and that $\Phi: U \rightarrow S$ is a degree 1 piecewise $C^1$ map with 
k-dilation $\Lambda$.  Let $Q_i = S_i/ R_i$.  Then for each integer $0 \le j \le k$ and 
each number $k+1 \le l \le n$, the following inequality holds.

$$\Lambda \ge c(n) (Q_1 ... Q_j) (Q_{j+1} ... Q_l)^{\frac{k-j}{l-j}}. \eqno{(*)}$$

Conversely, for each pair of rectangles $R, S$, there is a set $U \subset R$ and a degree 1 map from $U$ to $S$ with k-dilation at most $C(n) \max_{j,l} (Q_1... Q_j) (Q_{j+1} ... Q_l)^{\frac{k-j}{l-j}}$.

\end{reftheorem}

When I began the work on the material in this paper, 
one of my motivations was to find a new proof of $(*)$.  
According to the last proposition, we have the following lower bounds for $\Lambda$.

$$\Lambda \ge \mathbb{V}_S(O a(k,S)) / \mathbb{V}_R(O a(k,R)).$$

Hence, if we could estimate the minimax volumes of a rectangle up to a constant factor, we would get many lower bounds for $\Lambda$.  We end with a conjecture about the minimax volumes that would imply $(*)$.

\begin{conj} For each $k+1 \le l \le n$, the minimax volume $\mathbb{V}_R(Sq_{n-l}^Q a(k,R))$ is equal to the
following expression up to a constant factor $C(n)$.

$$\mathbb{V}_R(Sq_{n-l}^Q a(k,R)) \sim \inf_{0 \le j \le k} R_1... R_j (R_{j+1} ... R_l)^{\frac{k-j}{l-j}} 2^{\frac{l-k}{l-j} Q}.$$

\end{conj}

\end{document}